\documentclass{amsart}
\usepackage{amsmath,amsthm,amssymb,amsfonts}
\usepackage[colorlinks=true]{hyperref}

\hypersetup{urlcolor=blue, citecolor=blue}

\def\doi#1{   {\href{http://dx.doi.org/#1}
   {{\mdseries\ttfamily DOI}}}}

\newcommand{\Hs}[1]{H^{#1}_{\mathrm{rad}}}
\newcommand{\Hsd}[1]{\dot H^{#1}_{\mathrm{rad}}}

\newcommand{\al}{\alpha}    \newcommand{\be}{\beta}
\newcommand{\de}{\delta}    
  \newcommand{\ep}{\varepsilon}
  \newcommand{\La}{\Lambda}  \newcommand{\la}{\lambda}
    
\newcommand{\om}{\omega}    
\newcommand{\ga}{\gamma}    
\newcommand{\R}{\mathbb{R}}\newcommand{\Z}{\mathbb{Z}}
\newcommand{\N}{\mathbb{N}}

\newcommand{\M}{\mathcal{M}}\newcommand{\Sc}{\mathcal{S}}

\newcommand{\pt}{\partial_t}\newcommand{\pa}{\partial}
\newcommand{\D}{{\mathrm D}}
\newcommand{\beeq}{\begin{equation}}\newcommand{\eneq}{\end{equation}}

\newcommand{\ang}{{\not\negmedspace\nabla}}

\newcommand{\ft}{{\mathcal{F}}}

\newcommand{\Sp}{{\mathbb S}}\def\CO{\mathcal {O}}

\newenvironment{prf}{\noindent {\bf Proof.} }{\endprf\par}
\def \endprf{\hfill  {\vrule height6pt width6pt depth0pt}\medskip}

\numberwithin{equation}{section}

\newcommand{\gm}{\mathfrak{g}}

\def\O{{\mathcal{O}}}
\newcommand{\ges}{{\gtrsim}}
\newcommand{\les}{{\lesssim}}

\def\<{\langle}             \def\>{\rangle}
\def\({\left(}                 \def\){\right)}

\newtheorem{thm}{Theorem}[section]

\newtheorem{lem}[thm]{Lemma}
\newtheorem{prop}[thm]{Proposition}

\theoremstyle{definition}

\newtheorem{rem}[thm]{Remark}

\theoremstyle{definition}

\title[Sharp local well-posedness for quasilinear wave equations]
      {Sharp local well-posedness for quasilinear wave equations 
      with spherical symmetry}

\author{Chengbo Wang}\address{School of Mathematical Sciences\\                Zhejiang University\\                Hangzhou 310027, P. R. China}\email{wangcbo@zju.edu.cn}
\urladdr{http://www.math.zju.edu.cn/wang}

\date{\today}
\dedicatory{} \commby{}

\begin{document}

\begin{abstract}
In this paper, we prove a sharp local well-posedness result for 
spherically symmetric solutions to quasilinear wave equations with rough initial data, when the spatial dimension is 
three or higher.
Our approach is based on Morawetz type local energy estimates with fractional regularity for linear wave equations with variable $C^1$ coefficients, which rely on 
multiplier method, weighted Littlewood-Paley theory, duality and interpolation.
Together with weighted linear and nonlinear estimates (including weighted trace estimates, Hardy's inequality, fractional chain rule and  fractional Leibniz rule) which are adapted for the problem, the well-posed result is proved by iteration.
In addition, our argument yields almost global existence for $n=3$ and global existence for $n\ge 4$, when the initial data are small, spherically symmetric with almost critical Sobolev regularity.
\end{abstract}

\keywords{quasilinear wave equation,
local energy estimates, KSS estimates, trace estimates,
fractional chain rule, fractional Leibniz rule, unconditional uniqueness, almost global existence}

\subjclass[2010]{35L72, 35L05, 35L15, 35B30, 35B45, 35B65, 42B25, 42B37, 46B70}

\maketitle 

\section{Introduction}
Let $n\ge 3$, we are interested in the local well-posedness of the spherically symmetric solutions for the Cauchy problem of the quasilinear wave  equations with low regularity
\beeq\label{ea1}
\Box u+ g(u)\Delta u = a(u) u_t^2+b(u) |\nabla u|^2
\ , (t,x)\in (0, T)\times {\mathbb R}^n\ ,\eneq
\beeq\label{ea2}
u(0,x)=u_0(x) \in H^{s}_{{\rm rad}}({\mathbb R}^n), \partial_{t}u(0, x)=u_1(x) \in H^{s-1}_{{\rm rad}}
({\mathbb R}^n)\ ,\eneq
where
$\Box=- \partial_t^2+ \Delta$,
 $g, a, b$ are smooth functions, $g(0) = 0$ and such that
 $\Box+ g(u)\Delta$ satisfy the uniform hyperbolic condition. Here,  $H^{s}_{{\rm rad}}$ stands for the space of spherically symmetric functions lying in the usual Sobolev space $H^s$.

The equation \eqref{ea1} is scale-invariant in the sense that
$u_\la(t,x)=u(t/\lambda, x/\la)$ solves \eqref{ea1}
for every $\la>0$, provided that 
$u(t,x)$ is a solution. This gives us the critical homogeneous Sobolev space $\dot H^{s_c}$ with $$s_{c}=\frac{n}{2}\ ,$$
 which is known to be a lower bound of the regularity for the problem to be well-posed in $H^s$.
On the other hand, another characteristic feature of the  wave equations is that the propagation of singularities along the light cone, which heuristically yields ill posedness for the problem at the regularity level
 $s\le s_l=(n+5)/4$.
 
 For semilinear wave equations, that is, $g\equiv 0$,
it could be shown to be locally well-posed in
$H^s$ for $s>n/2+1-1/q$, with $q=\max(2, 4/(n-1))$,
 by $L^{q+}_t L_x^\infty$ Strichartz estimates.
Moreover, 
 it is known that the problem is
 locally well-posed if $$s>\max(s_{c}, s_{l})\ ,$$ and ill-posed in general, if $s<s_{c}$ or $s\le s_{l}$, see
 Ponce and Sideris \cite{PS} ($n=3$), 
 Tataru \cite{Tataru99} ($n\ge 5$), 
 Zhou \cite{Zh03} ($n=2,4$) for positive results, and 
 Lindblad \cite{Ld93,Ld96}  ($n=3$),
Fang-Wang \cite{FW1} ($n\ge 6$),
Liu-Wang \cite{LiuW21}  ($n\le 5$) for negative results. The critical well-posedness remains open for higher dimensional case ($n\ge 6$).

If the nonlinearity is of the first null-form, that is,
$a(u) u_t^2+b(u) |\nabla u|^2=c(u)(u_{t}^{2}-|\nabla u|^{2})$
for some $c(u)$, improved local well-posed results are also available, which states that $s>s_{c}$ is sufficient for  local well-posedness, see, e.g., Klainerman-Machedon 
\cite{KM95} ($n=3$), Klainerman-Selberg \cite{KS97} ($n\geq 2$), see also Liu-Wang \cite{LiuW21}.

 Furthermore, it is well-known that we could extend the admissible pairs for Strichartz estimates, when the initial data are spherically symmetric  or have certain amount of angular regularity, see
Klainerman-Machedon 
 \cite{KlMa93},
 Sterbenz \cite{Ster07}, Machihara-Nakamura-Nakanishi-Ozawa \cite{MaNaNaOz05},
 Fang-Wang \cite{FW2}. With help of this observation, we could improve the radial results to 
 $s>3/2$ for $n=2$ and $s\ge 2$ for $n=3$. We see that there are still $1/2$ gap of regularity, between the positive results and the scaling regularity.
 When $n=3$ in the case of radial small data, by exploiting the local energy estimates and weighted fractional chain rule, the regularity assumption is improved to the almost critical assumption $s>3/2$, in Hidano-Jiang-Lee-Wang \cite{HJLW20p}, with previous results
of Hidano-Yokoyama \cite{HiYo06} for $s=2$. In view of 
\cite{HJLW20p}, it seems that the critical radial regularity for $n=2$ is $s=3/2$ instead of the scaling critical regularity $s=1$.
 Concerning radial solutions to general semilinear hyperbolic systems in 3D under null condition, 
the global existence for small scaling invariant
$\dot W^{2,1}(\R^3)$ data is known from Yin-Zhou \cite{YinZhou16}.

Turning to the quasilinear problem
 \eqref{ea1}, it is much more delicate. Based on the classical energy argument, it is locally well-posed, as long as $s>n/2+1$ (\cite{MR0420024}).
Similar to the semilinear problem, the approach of using $L^{q+}L^\infty$ Strichartz estimates has been intensively investigated. 
To make the argument work, we need to obtain Strichartz estimates for wave operators with variable coefficients.
It is known that 
we have the full Strichartz estimates, provided that the perturbation is of $C^{1,1}$, see
 Smith \cite{MR1644105},
Tataru \cite{Ta02}, as well as 
Kapitanskii \cite{MR1031987}, Mockenhaupt-Seeger-Sogge \cite{MR1168960} for previous results with smooth perturbation. However, in view of application to
the quasilinear problem
 \eqref{ea1}, $u\in H^s$ with $s<n/2+1$ will only imply $g(u)\in C^{0, s-n/2}$, by Sobolev embedding, which means that it is desirable to obtain the 
Strichartz estimates for wave operators with rough coefficients.

The first breakthrough was achieved through the independent works of Bahouri-Chemin \cite{MR1719798} (Hadamard parametrix) and Tataru \cite{MR1749052}
(FBI transform),
where weaker Strichartz estimates for metric with limited regularity was obtained and so is the local well-posedness for
$$s>\frac{n+2} 2-\left\{\begin{array}{ll }
1/4 ,     &   n\ge 3\ ,\\
1/8,      &   n=2\ .
\end{array}\right.$$
This approach was developed further, to arrive
$$s>\frac{n+2} 2-\left\{\begin{array}{ll }
1/3 ,     &   n\ge 3\ ,\\
1/6,      &   n=2\ ,
\end{array}\right.$$
see Bahouri-Chemin
\cite{MR1728676}, Tataru \cite{Ta02}.
To get further improved result, it is desirable to exploit
the additional geometric information on the metric $g(u)$ and solution $u$.
In the work of
 Klainerman-Rodnianski \cite{MR1962783},  
the nonlinear structure of solutions was exploited to obtain the improved result  $s>3-\sqrt{3}/2$ when $n=3$.
Finally, the well-posedness for
$$s>\frac{n+2} 2-\left\{\begin{array}{ll }
1/2 ,     &   3\le n\le 5\ ,\\
1/4,      &   n=2\ ,
\end{array}\right.$$
 was proven by
Klainerman-Rodnianski \cite{KlRo05} for the 
Einstein-vacuum equations ($n=3$), and
Smith-Tataru \cite{SmTa05} for general quasilinear wave equations
in general spatial dimensions, by constructing a parametrix using wave packets.
Later, Q. Wang \cite{WangQ17} gives an alternative proof for Smith-Tataru's result when $n=3$, by 
commuting vector fields approach.

When $n=3$ and $2$, we know from 
Lindblad \cite{MR1666844} and
Liu-Wang \cite{LiuW21} that the well-posed result in \cite{SmTa05} is sharp in general.
However, concerning the Einstein-vacuum equations, 
the so-called bounded {$L^2$} curvature conjecture (well-posed in $H^2$) was verified in
Klainerman-Rodnianski-Szeftel
\cite{MR3402797}. In contrast,
Ettinger-Lindblad
\cite{MR3583356} 
proved ill-posed result in $H^2$ for Einstein-vacuum equations in the harmonic gauge.

In summary, the quasilinear problem is 
 locally  well-posed in $H^s$ for 
 $$s>\left\{\begin{array}{ll}
(n+5)/4   
      &  n\le 3,\\
   (n+1)/2   &   n=4, 5,\\
n/2+2/3	& n\ge 6\ ,
\end{array}\right.$$
in general.

Comparing with the semilinear problem, we  expect naturally that there should be improved well-posed theory, when
the problem and the initial data are spherically symmetric. Actually,
in  \cite{HWY1}, 
together with Hidano and Yokoyama, we proved that
the $3$-dimensional problem
$$\Box u+g(u)\Delta u=a u_t^2+b |\nabla u|^2$$
 is well-posed for small radial data in $\Hs{2}\times \Hs{1}$, with
 almost global existence of solutions, up to $\exp(c/
\|(\nabla u_0, u_1)\|_{H^1})$ (see also Zhou-Lei \cite{ZhouLei08} for previous work on global existence with $a=b=0$, for compactly supported $\Hs{2}\times \Hs{1}$ data).
On the other hand, when $n=2$, 
the improved local well-posed result for $H^s_{\mathrm{rad}}$ with $s>3/2$ was suggested in Fang-Wang \cite{FWQLW}. 
In addition, as we have mentioned, when $n=2, 3$,
the long time well-posedness with small radial data in $\Hs{s}$ with $s>3/2$ is known from Hidano-Jiang-Lee-Wang \cite{HJLW20p}.

\subsection{Main results}
Let us turn to the first main result of this paper, concerning
 the physically important case, $n=3$.
As we know, it is  well-posed in $\Hs{2}$ (at least for small data) and generally ill-posed in $\Hs{s}$ with $s<s_c=3/2$. 
Heuristically, comparing with the semilinear results, we may expect well-posedness in $\Hs{s}$ for certain $s<2$.
However,
we notice that $\Hs{2}$ is the lowest possible regularity we could obtain by the approach of ($L^2_t L^\infty_x$) Strichartz estimates, even in the radial case.
To break the Sobolev regularity barrier $s=2$, we need to 
circumvent the Strichartz estimates approach.

In the following result, we could prove well-posedness in $\Hs{s}$ for any subcritical regularity, $s>3/2$, which shows that
there are no any other obstacles to well-posedness in the radial case, except scaling. As far as we know, this might be the first well-posed result for three-dimensional quasilinear wave equations, which breaks the Sobolev regularity barrier, $s=2$.
More precisely, we prove the following, with certain low frequency condition on $u_1$.
\begin{thm}\label{th-1}
Let $n= 3$,
$s\in (3/2, 2]$ and
$s_0\in [2-s, s-1]$.
 Considering \eqref{ea1}-\eqref{ea2} with 
 $u_0\in H^{s}_{{\rm rad}}$ and
 $u_1\in 
\dot H^{s-1}_{{\rm rad}}\cap
 \dot H^{s_0-1}_{{\rm rad}}$.
There exists $T_0>0$ such that 
the problem is (unconditionally) locally well-posed in the function space 
 \beeq\label{eq-1-3} u\in L^{\infty}_t \Hs{s}  \cap C^{0,1}_t\Hs{s-1}\cap C_t \dot H^{s_0}([0, T_0]\times \R^3), \pt u \in C_t\dot H^{s_0-1} \ . \eneq
 More precisely, 
 \begin{enumerate}
  \item {\rm (Existence)}
There exists
a universal constant $C>0$, 
  so that
 there exists  a (weak) solution $u$ satisfying \eqref{eq-1-3} and
 $$
\|\pa u\|_{L^\infty \dot H^{\theta}}+
T^{-\frac{\mu}2}\|r^{-\frac{1-\mu}2}D^{\theta}\pa u\|_{L^2([0,T]\times\R^3)}
\le C\|(\nabla u_0, u_1)\|_{\dot H^{\theta}}\ , 
$$
 for all $\theta\in \{s_0-1\}\cup [0, s-1]$ and $T\in (0, T_0]$.
Here $\mu=s-3/2$, $D=\sqrt{-\Delta}$.
  \item {\rm (Uniqueness)} the solution is unconditionally unique in \eqref{eq-1-3}.
  \item {\rm (Persistence of regularity)} Let $T_*$ be the maximal time of existence (lifespan) for the solution in \eqref{eq-1-3}.
If $(u_0, u_1)\in H^{s_1}\times H^{s_1-1}$ for some 
$s_1\ge 3$, 
 then the solution $u\in L^\infty H^{s_1}\times C_t^{0,1} H^{{s_1}-1}$ in $[0,T]\times\R^3$ for any $T<T_*$.
  \item {\rm (Continuous dependence)} We also have 
continuously dependence on the data when $s_0<s-1$, in the following sense:
for any $T\in (0, T_*)$,
$s_1\in (s_c, s)$
 and  $\ep>0$,
there exists $\de>0$, such that
whenever
$\|(\nabla (u_0-v_0), u_1-v_1)\|_{\dot H^{s-1}\cap \dot H^{\max(s_1-2, s_0-1)}}\le \de
$, the corresponding solution $v\in L^\infty H^{s_1}\times C_t^{0,1} H^{{s_1}-1}$ is well-defined in $[0, T]\times \R^3$ and 
 $$
\|\pa (u-v)\|_{L^\infty (\dot H^{s_1-1}
\cap \dot H^{\max(s_1-2, s_0-1)})}\le\ep\ .$$
\end{enumerate}
 \end{thm}

\begin{rem} \label{rem-1.1}
The regularity assumption of the lifespan obtained in Theorem \ref{th-1} is sharp in general. 
More precisely, we could not have well-posedness for data
in some critical space,  $B$, and possibly non-subcritical space 
$\dot H^s$ with $s\le s_c=3/2$.
Actually, let $g=0$, $a=1$, $b=0$, and $\phi, \psi$
be given nonnegative nontrivial, spherically symmetric $C^\infty_0$ functions, then it is well-known, see, e.g., John \cite{Jo81}, that for classical solutions, for any $\ep>0$, the lifespan
$T_*<\infty$ for data $u_0=\ep \phi$, $u_1=\ep \psi$. By persistence of regularity, we know that
$T_*$ is the same as the lifespan for weak solutions.
If the problem is still well-posed, then by continuous dependence for the trivial solution, there exists $\de>0$ such that
 $T_*\ge \de$,
 for any data with critical norm, $\|(u_0, u_1)\|_B\le \de$ and $\ep_s=
 \|(\nabla_x u_0,u_1 )\|_{\dot H^{s-1}}\le \de$.
 Let $\ep\ll 1$ such that the norm $\le \de$ and we obtain a solution $u$ with
 $T_*<\infty$.
For such fixed $\ep>0$,
 by rescaling, 
we know that, for any $0<\la\le 1$, $u_\la(t,x)=u( t/\la,  x/\la)$ solves the equation with rescaled data, for which we have
$$
\|(u_\la(0), \pt u_\la(0))\|_B
=\|(u_0, u_1)\|_B
\le\de, \ 
\ep_{s, \la}=\la^{s_c-s} \ep_{s}\le \ep_{s}\le \ep_{0},\ T_{*,\la}=
\la T_* \ .$$
This gives us
$0<\de\le T_{*,\la}=
\la T_*<\infty$ for any  $0<\la\le 1$, which is clearly a contradiction.
On the other hand, we have an auxiliary low frequency regularity assumption on the initial velocity $u_1\in \dot H^{s_0-1}$, due to the second order feature of the equation and the limited regularity level $s<2$. This assumption plays a key role in our analysis,  to close the iteration,
 and it will be interesting to determine if it is essential for the well-posed result or not. Notice, however, that we do not need to assume $u_1\in  \dot H^{s_0-1}$ when it is compactly supported.
\end{rem}

Next, we present our high dimensional well-posed result.
\begin{thm}\label{th-high}
Let $n\ge 4$, 
$s=n/2+\mu$ with
\beeq\label{eq-mu}
\mu\in
\left\{
\begin{array}{ll }
  (0,  1 /2]  ,  &   n \textrm{ odd} \ ,\\
  (0, 1) ,      &    n \textrm{ even}\ .
\end{array}
\right.
\eneq
The
problem \eqref{ea1}-\eqref{ea2} is (unconditionally) locally well-posed in the function space 
 \beeq\label{eq-1-6} u\in L^{\infty}_t \Hs{s}  \cap C^{0,1}_t\Hs{s-1}\cap CH^1\cap C^1 L^2 \ . \eneq
 More precisely, 
 \begin{enumerate}
  \item {\rm (Existence)}
There exists
a  constant $C>0$
  so that
 for any data $(u_0,u_1)\in
\Hs{s}  \times \Hs{s-1}$, there exist $T
>0$, and a (weak) solution $u$ in \eqref{eq-1-6} in $[0,T]\times \R^n$, satisfying
 $$
\|\pa u\|_{L^\infty \dot H^{\theta}}
+T^{-\frac{\mu}2}\|r^{-\frac{1-\mu}2}D^{\theta}\pa u\|_{L^2([0,T]\times\R^n)}
\le C\|(\nabla u_0, u_1)\|_{\dot H^{\theta}}\ , 
$$
 for all $\theta\in  [0, s-1]$.
  \item{\rm  (Uniqueness)} the solution is unconditionally unique in \eqref{eq-1-6}.
  \item{\rm  (Persistence of regularity)} Let $T_*$ be the 
  lifespan. 
If $(u_0, u_1)\in H^{s_1}\times H^{{s_1}-1}$ for some $s_1\ge [(n+4)/2]$, then the solution $u\in L^\infty H^{s_1}\times C_t^{0,1} H^{{s_1}-1}$ in $[0,T]\times\R^n$ for any $T<T_*$.
  \item{\rm  (Continuous dependence)} We also have 
continuously dependence on the data, in the $H^{s_1}$ topology, for 
$s_1\in (s_c, s)$.
\end{enumerate}
 \end{thm}

When considering the small data problem, it turns out that we could give the following long time existence results.
\begin{thm}[Long time existence for small data]\label{th-2-al-global}
Let $n\ge 3$ and $s>s_c=n/2$. Considering \eqref{ea1}-\eqref{ea2} with 
$(u_0,u_1)\in \Hs{s}  \times \Hs{s-1}$.
When $n=3$,
assuming further that
$u_1\in \dot H^{s-2}$,
there exist $c>0$ and  $\de>0$ such that for any data with
$\ep_1+\ep_s<\de$, the problem
admits an almost global $L^\infty([0,T], H^{\min (s, 2)}(\R^3))$ solution,
where
\beeq\label{eq-1-4}
T=
\exp(c/(\ep_1+\ep_s))\ ,\eneq
\beeq\label{eq-1-7}\ep_s:= \|(\nabla_x u_0,u_1 )\|_{\dot H^{s-1}}\ ,\ 
\ep_c=\|(\nabla u_0, u_1)\|_{\dot B^{s_c-1}_{2,1}}
\ .\eneq
When $n\ge 4$,
for any $s>s_c$,
there exists $\ep>0$ such that 
the problem admits global solutions
whenever $\ep_s+\ep_{1}\le \ep$.
\end{thm}

\begin{rem} 
The lower bounds of the lifespan, \eqref{eq-1-4}, obtained in Theorem \ref{th-2-al-global} for $n=3$ is sharp in general, in terms of the order.
Actually, 
as in Remark \ref{rem-1.1},
for the sample case $a=1$, $g=b=0$, it is well known that there exists data $(\ep \phi, 
\ep \psi)$, so that, for any $\ep\in (0, 1]$, the lifespan of the classical solutions has upper bound $T_*\le \exp(C/\ep)$ for some $C>0$.
By the way, it is clear from the proof of Theorems \ref{th-1}-\ref{th-high}, that,
 when $|g|\ll 1$,
 we could obtain the following lower bound of the lifespan
 $$T_*\ge c(g,a,b,\ep_c)\ep_s^{-1/(s-s_c)}\ .$$
 Moreover, when
 $\ep_c\ll 1$,
  $$T_*\ge c(g,a,b)\ep_s^{-1/(s-s_c)}
  \exp(c(g,a,b)/\ep_c)\ .$$
On the other hand, 
when it satisfies certain  nonlinear structures, such as the null conditions, or many case of the weak null conditions,
the problem admits global solutions with small data. See, e.g., 
\cite{MR2003417},
\cite{Ld08},
\cite{ZhouLei08} for
global results with $a=b=0$.
In such situations, we naturally expect
that global radial results with $s>s_c$ (or even certain critical space like $\dot B^{s_c}_{2,1}$) remain hold, 
which is an interesting further problem.
\end{rem}

\begin{rem} 
Although we state only the results for scalar quasilinear wave equations, as it is clear from the proof, our results apply also for general multi-speed system  of quasilinear wave equations, which permit spherically symmetric solutions.
In particular, the system with multi-speeds ($c_j>0$)
$$\pt^2 u^j-c_j^2 (1+ g_j(u))\Delta u^j= 
Q_{kl}^{j\al\be} (u ) \pa_\al u^k \pa_\be u^l,
1\le j\le N,
$$ is local well-posed in $\Hs{s}(\R^n)
\times(\Hs{s-1}(\R^n)\cap \dot H^{s-2}(\R^n))
$ for $n\ge 3$ and $s>n/2$, 
 as long as the system admits spherically symmetric solutions.
Similar statement holds 
for
$$\pt^2 u^j-c_j^2 (1+ g_j(u,\pa u))\Delta u^j= 
Q_{kl}^{j\al\be} (u ) \pa_\al u^k \pa_\be u^l, 1\le j\le N,
$$
in $\Hs{s}(\R^n)
\times(\Hs{s-1}(\R^n)\cap \dot H^{s-3}(\R^n))
$ 
when $n\ge 3$ and $s>(n+2)/2$.

 In addition, 
the quasilinear part could be replaced by the D'Alembertian
 $\Box_{g(u)}$ with respect to the metric $ds^2=-dt^2+g(u) dx^2$,  as well to $\Box_\gm+g(u)\Delta$
(or $\Box_{g(t,x,u)}$), when $\gm$ is a
small, long range perturbation of the Minkowski metric:
$$ \gm=-K_0(t,r)^2 d t^2
+2K_{01}(t,r)dtdr
+K_1^{2}(t, r)dr^{2}+r^{2}d\omega^{2}\ ,
$$
$$
|(K_0-1, K_{01}, K_1-1)|
\ll 1,
\sum_{j\ge 0}\|r^{|\ga|-\mu}\<r\>^\mu\pa^\ga 
(K_0-1, K_{01}, K_1-1)
\|_{L^\infty_{t,x}(1+|x|\sim 2^j)}\ll 1, 
$$
for $1\le|\ga|\le [n/2]+1$ (and $K_{01}=0$ when $n=3$).
\end{rem}

\subsection{Idea of proof}
Let us discuss the idea of proof. Basically, we rely on the approach of using Morawetz type local energy estimates, instead of Strichartz estimates. As have appeared in many works on dispersive and wave equations, 
 Morawetz type local energy estimates have been proven to be more fundamental and robust than Strichartz estimates, in many nonlinear problems.
 
To make the approach work for quasilinear wave equations, similar to the approach of using Strichartz estimates, we prove a version of Morawetz type local energy estimates,
Theorem \ref{th-LE0}, for linear wave equations with variable $C^1$ coefficients. It is this version of local energy estimates which makes it possible to relax the regularity requirement for quasilinear wave equations. The proof is based on the classical multiplier approach with well-chosen multipliers, which yields such estimates for small perturbation of the Minkowski metric. Furthermore, the property of finite speed of propagation is exploited to handle the general case of large perturbation.

With the help of the well-adapted Morawetz type local energy estimates (weighted space-time $L^2$ estimates), we are naturally led to develop the corresponding linear and nonlinear estimates involving weighted functions. Among others, 
we prove weighted Sobolev type estimates (including 
weighted trace estimates,
Proposition \ref{thm-trace-w},
and
weighted Hardy's inequality,
Lemma \ref{thm-w-Hardy}),
weighted fractional chain rule (Theorem \ref{thm-wLeib0} and Proposition \ref{thm-wchain-high}), as well as the 
weighted Leibniz rule 
(Theorem \ref{thm-wLeib-fromChain}).

The Morawetz type local energy estimates,
Theorem \ref{th-LE0}, is at the regularity level of $\dot H^1$. 
To make the approach work, we need to develop the corresponding version of  local energy estimates, at the regularity level $\dot H^s$ with $s>n/2$.
With the help of interpolation, Littlewood-Paley theory involving weighted functions, together with 
the 
weighted Sobolev type estimates from
Proposition \ref{thm-trace-w},
Lemma \ref{thm-w-Hardy}  and
 Lemma \ref{thm-w-Hardy-trace},
we prove a series of
local energy estimates with positive fractional derivatives,
Propositions \ref{th-LEkey1},
\ref{th-LEkey2} and
\ref{th-LEkey4}.

 Equipped with all these linear and nonlinear estimates, 
 we could then use the standard iteration argument to establish local existence and uniqueness, as well as the long time existence. In particular, for the case of $n=3$, to prove convergence of approximate solutions, we 
develop local energy estimates with negative regularity, and
 need to assume certain low frequency requirement on the initial velocity, due to the second order feature of the equation and the limited regularity level $s<2$.

Recall that, in the approach of using Strichartz estimates, the proof of local existence immediately implies the 
persistence of regularity and 
 continuous dependence on the data, through Gronwall's inequality. Unlike the approach of using Strichartz estimates, in our approach, the proofs for persistence of regularity and 
 continuous dependence on the data are not so direct.
 For example, concerning persistence of regularity, we prove first the result  from regularity index $s$ to any $s\in (s_c, [(n+2)/2])$, and then prove the persistence to $s=[(n+4)/2]$, which is sufficient to conclude  persistence of higher regularity.

 This paper is organized as follows.
In the next section, 
we recall and prove various basic linear and nonlinear estimates, including
weighted Sobolev type estimates,
weighted trace estimates,
weighted Hardy's inequality,
as well as the
weighted fractional chain rule and the 
weighted Leibniz rule.
In Section \ref{sec-3le}, we present a version of
Morawetz type local energy estimates, for linear wave equations with variable $C^1$ coefficients, as well as the estimates with fractional regularity.
In Sections \ref{sec-4} and \ref{sec-5-high}, by iteration argument, we prove local existence and uniqueness, for $n=3$ and $n\ge 4$.
In Section \ref{sec-6persi},
we show that 
persistence of regularity for the weak solutions, when the initial data have higher regularity,
as well as the continuous dependence on the data.
Next, in Sections \ref{sec-7-n=3} and \ref{sec-8}, we present the proof of almost global existence and global existence for $n=3$ and $n\ge 4$, when the initial data are small.
Finally, in the appendix, we present the fundamental Morawetz type estimates, by elementary multiplier approach, with carefully chosen multipliers.

\subsection{Notations}
We close this section by listing the notations.

 \noindent $\bullet$  $\ft(f)$ and $\widehat{f}$  denote
the Fourier transform of $f$. 
$D=\sqrt{-\Delta}:=\ft^{-1} |\xi|\ft$ and $P_j=\phi(2^{-j}D)$ is the
 (homogeneous)
Littlewood-Paley projection on the space-variable, $j\in\Z$.

\noindent  $\bullet$
$r=|x|$, $\<r\>=\sqrt{2+r^2}$, 
$\pa=(\pt, \nabla_x)=(\pt, \nabla)$, 
$\tilde \pa u=(\pa u, u/r)$,
$|\nabla^k u|=\sum_{|\ga|=k} |\nabla^\ga u|$ for multi-index $\ga$.

\smallskip

\noindent $\bullet$ $L^p(\R^n)$ denotes the usual Lebesgue space,
and $L_r^p(\R^+)=L^p(\R^+:r^{n-1}dr)$.

\noindent $\bullet$  $L_r^p L_\omega^q$ is
Banach space defined  by the  following norm
\[\|f\|_{L_{r}^p L_\omega^q}=\big\|{\|f(r\omega)\|_{L_\omega^q}}\big\|_{L_{r}^p}.\]

\noindent $\bullet$  $H^s$, $\dot{H}^s$ ($B^s_{p,q}$,
$\dot{B}^s_{p,q}$) are the usual  inhomogeneous and homogeneous Sobolev (Besov) spaces on $\R^n$.

\noindent $\bullet$ 
With parameters $\mu, \mu_1\in (0,1)$ and $T\in (0,\infty)$, we define
\begin{align}
\label{eb5-LE}
\|u\|_{LE_{T}}=  &\|\pa u\|_{L^\infty_t L^2_x}+   \|r^{-\frac{1-\mu}2}\<r\>^{-\frac{\mu+\mu_1}2} \partial u\|_{L^2_{t,x}}
\\&+\<T\>^{-\frac{\mu}2}\|r^{-\frac{1-\mu}2} \partial  u\|_{L^2_{t,x}}
      +(\ln\<T\>)^{-\frac 12}
\|r^{-\frac{1-\mu}2}
\<r\>^{-\frac{\mu}2}
 \partial u\|_{L^2_{t,x}}\ ,
\nonumber
\end{align}
for functions
 on $[0,T]\times \R^n$.
In the limit case $T=\infty$, 
 we set 
$$\|u\|_{LE}= \|\pa u\|_{L^\infty_t L^2_x}+   \|r^{-\frac{1-\mu}2}\<r\>^{-\frac{\mu+\mu_1}2} \partial u\|_{L^2_{t,x}}+\sup_{T>0}\|u\|_{LE_{T}}\ .$$
In addition, for fixed $\mu\in (0,1)$,
\beeq
\label{eq-XT}\|u\|_{X_T}:=
\|  u\|_{L^\infty_t L^2_x}+
T^{-\frac{\mu}2}\|r^{-\frac{1-\mu}2} u\|_{L^2_{t,x}}\ ,\eneq
\beeq\label{eq-XT*}
\|F\|_{X^*_T}:=\inf_{F=F_1+F_2}
(\| F_1\|_{L^1_t L^2_x}+
T^{\frac{\mu}2}\|r^{\frac{1-\mu}2} F_2\|_{L^2_{t,x}})\ .
\eneq
For $q\in [1,\infty]$, we introduce the Besov version as follows
$$
\| u\|_{X_{T,q}}
 :=
 \| u\|_{L^\infty_t \dot B^{0}_{2,q}}   +
T^{-\frac{\mu}2}\|r^{-\frac{1-\mu}2} P_j u \|_{\ell_{j}^q L^2_{t,x}}\ .$$

\section{Sobolev type and nonlinear estimates}

In this section, we recall and prove various basic estimates to be used.

\subsection{weighted Sobolev type estimates}
We will use the following version of the weighted Sobolev estimates, which essentially are consequences of the well-known trace estimates.
\begin{lem}[Trace estimates]
 \label{thm-trace}
Let $n\geq 2 $ and $s\in [0, n/2)$, then
\begin{equation}
\|r^{(n-1)/2}u\|_{L_{r}^{\infty}L_{\omega}^{2}}\les \|u\|_{\dot B^{1/2}_{2,1}},\ \| r^{n/2-s}u\|_{L_{r}^{\infty}H_{\omega}^{s-1/2}} \les\|u\|_{\dot{H}^{s}}, s>1/2
\ ,
\label{eq-2.1}
\end{equation}
\beeq\label{eq-trace2-sd}
\|r^{n(1/2-1/p)-s} f\|_{L^p_r L^{2}_\omega}\les \|f\|_{\dot H^{s}}, 2\le p<\infty,
1/2-1/p\le s<n/2
\ .
\eneq
\end{lem}

The estimate \eqref{eq-2.1} is well-known, see, e.g., 
\cite[(1.3), (1.7)]{FW11} and references therein.
The inequality \eqref{eq-trace2-sd} with $s=1/2-1/p$ is due to \cite{LiZh95}, see also \cite{HWY3} for alternative proof using real interpolation and \eqref{eq-2.1}.

We shall also use the following weighted variant of the trace estimates.
\begin{prop}[Weighted trace estimates]\label{thm-trace-w}
Let $n\ge 2$, $\al\in (1/2, n/2)$ and $\be \in (\al-n/2, n/2)$. Then we have
\beeq\label{eq-trace-w}
\|r^{n/2-\al+\be}P_j u\|_{l^2_j L^\infty_{r} H^{\al-1/2}_\omega}\les 
\|r^{\be} D^{\al} u\|_{L^{2}}\ ,
\eneq
\beeq\label{eq-trace-w'}
\|r^{n/2-\al+\be} u\|_{L^\infty_{r} H^{\al-1/2-}_\omega}\les 
\|r^{\be} 2^{j\al}P_j u\|_{l^\infty_j L^{2}}\ .
\eneq
In addition, we have
\beeq\label{eq-trace-w2}
\|r^{n(1/2-1/p)-\al+\be} u\|_{L^p_{r} L^{2}_\omega}
\les 
\|r^{\be} D^{\al} u\|_{L^{2}}\ , 
\eneq
for any $p\in [2,\infty]$,  $\al\in (1/2-1/p, n/2)$, and $\be\in (\al-n/2,n/2)$.
\end{prop}
\begin{prf}
We essentially follow \cite[Lemma 4.2]{W17}, where 
\eqref{eq-trace-w} was proven for $\al\in (1/2, 1]$ and $n\ge 3$.
Recall that we have the following weighted Littlewood-Paley square-function estimate 
\beeq\label{eq-LPw}\|w P_j f\|_{L^p \ell^2_j}\simeq \|w f\|_{L^p}, w^{p}\in A_p, f\in L^p(w^p dx), p\in (1,\infty)\ .\eneq
As $r^{2\be}\in A_{2}$ if and only if $|\be|<n/2$,
we get 
\beeq\label{eq-trace-w4}
\|r^{\be} 2^{j\al} P_{j} u\|_{\ell_{j}^{2}L^{2}}
\simeq \|r^{\be} \D^{\al} u\|_{L^{2}}
\ , \be \in (-n/2, n/2)\ .
\eneq

Based on this estimate, we observe that, by rescaling, interpolation and frequency localization, the proof
of \eqref{eq-trace-w} and \eqref{eq-trace-w'} can be reduced to
 the following estimate
\beeq\label{eq-trace-w1}
\|u(\om)\|_{H^{\al-1/2}_\omega(\Sp^{n-1})}\les 
\|r^{\be} \nabla^k u\|_{L^{2}}^{\al/k}\|r^{\be}  u\|_{L^{2}}^{1-\al/k}\ , \al\in [1/2,n/2), 
\eneq
where $\be \in (\al-n/2, n/2)$, $k\in (\al, \al+1]$.

For the proof of \eqref{eq-trace-w1}, we
recall the weighted Hardy-Littlewood-Sobolev estimates of Stein-Weiss 
\beeq\label{eq-SteinWeiss}\|r^{\be-\al} u\|_{2}\les \|r^{\be}\D^{\al} u\|_{2}, \al\in (0,n), \be\in (\al-n/2,n/2)\ .\eneq Then for any $\al\in (0, n)$ and $\be\in (\al-n/2,n/2)$, we have
\beeq\label{eq-SteinWeiss1}\|r^{\be-\al} u\|_{2}\les \|r^{\be}\D^{\al} u\|_{2}
\les \|r^{\be}\nabla^k u\|_{2}^{\al/k}\|r^{\be} u\|_{2}^{1-\al/k}
\ ,\eneq
if $\al<k$.
Moreover, if $k\in (\al, \al+1]\cap \N$, we have
\beeq\label{eq-SteinWeiss2}\|r^{\be-j}  u\|_{L^{2}}\les
\|r^\be \nabla^j u\|_{L^2}\ , \forall j<k\ ,
\eneq
for such $\al, \be$.
Let $\phi$ be a cutoff function of $B_{2}\backslash B_{1/2}$ which equals one for $|x|=1$, we get from
\eqref{eq-2.1} that for $\al \in [1/2, n/2)$ 
 and $\be\in (\al-n/2,n/2)$
\begin{eqnarray*}\|u(\omega)\|_{H^{\al-1/2}_\omega}&\les &
\| \nabla^k (\phi u)\|_{L^{2}}^{\al/k}\|\phi u\|_{L^{2}}^{1-\al/k}\\
&\les& (\sum_{j< k} \|r^{\be-j} \nabla^{k-j} u\|_{L^{2}})^{\al/k}\|r^{\be}  u\|_{L^{2}}^{1-\al/k}
+\|r^{\be-\al} u\|_{L^{2}}\\
&\les& \|r^{\be} \nabla^k u\|_{L^{2}}^{\al/k}\|r^{\be}  u\|_{L^{2}}^{1-\al/k}\ ,
\end{eqnarray*}
where we have used 
\eqref{eq-SteinWeiss1} and
\eqref{eq-SteinWeiss2}.
This gives us \eqref{eq-trace-w1}, and so is
\eqref{eq-trace-w} and \eqref{eq-trace-w'}.

Finally, \eqref{eq-trace-w2} follows directly from interpolation
between
\eqref{eq-SteinWeiss}, \eqref{eq-trace-w}
and
\eqref{eq-trace-w'}.
This completes the proof.
\end{prf}

\subsection{weighted fractional chain rule}
When dealing with the nonlinear problems,
it is natural to introduce the weighted fractional chain rule and Leibniz rule.
We first present the following generalized version of the weighted fractional chain rule of Hidano-Jiang-Lee-Wang \cite{HJLW20p}, which could be viewed as a transition from Sobolev type norm to Besov type norm, as well as a transition from space variables to space-time variables.
For the weight functions, we
recall the Muckenhoupt $A_p$ class,  which by definition,
$$w\in A_1\Leftrightarrow \mathcal{M}w(x)\le Cw(x), \textrm{a.e.}\ x \in\R^n\ ,$$
$$w\in A_p\ (1<p<\infty)\Leftrightarrow \left(\int_Q w(x)dx\right)\left(\int_Q w^{1-p'}(x)dx\right)^{p-1}\le C|Q|^p, \forall \textrm{ cubes }Q\ ,$$ 
with $\M w(x)=\sup_{r>0}r^{-n}\int_{B_r(x)} w(y)dy$ denotes the Hardy-Littlewood maximal function. See, e.g., \cite[\S 2.5.2]{MuSch13}.

\begin{thm}[Weighted fractional chain rule]\label{thm-wLeib0}
Let $s\in (0,1)$,  $\lambda\ge 1$, $q, q_{1}, q_2 \in (1,\infty)$,
 $p, p_{1}, p_2 \in [1,\infty]$
 with \beeq\label{eq-HJLW0}\frac 1q=\frac{1}{q_{1}}+\frac{1}{q_{2}},\ 
 \frac 1p=\frac{1}{p_{1}}+\frac{1}{p_{2}}\ .
\eneq
Assume $F:\R^k\rightarrow \R^l$ is a $C^1$ map, satisfying $F(0)=0$ and 
\beeq\label{eq-wLeib-assu0}|F'(\tau v+(1-\tau)w)|\le \mu(\tau)|G(v)+G(w)|,\eneq
with $G>0$ and $\mu\in L^1([0,1])$. If
$(w_1 w_{2})^{q}\in A_{q}$, $w_{1}^{q_{1}}\in A_{q_{1}}$,
$w_{2}^{q_2}\in A_{q_{2}}$, then we have
\beeq\label{eq-wLeib5}\|w_1w_{2} 2^{js}
P_j F(u)\|_{\ell_{j\in\Z}^\lambda L^p_t  L^q_{x} }\les \|w_1 2^{js}
P_j  u\|_{\ell_{j\in\Z}^\lambda L^{p_{1}}_{t} L^{q_{1}}_{x}} \|w_{2} G(u)\|_{L^{p_2}_t L^{q_{2}}_x}\ ,\eneq
for any $[0,T)\times \R^n\ni (t,x)\rightarrow u(t,x)\in\R^k$.
In addition, when $q_{2}=\infty$ and $q\in (1,\infty)$, if
$w_1^{q}, (w_{1}w_2)^{q}\in A_q$ and $w_2^{-1}\in A_1$, we have
\beeq\label{eq-wLeib6}
\|w_1 w_{2}  2^{js}
P_j  F(u)\|_{
\ell_{j\in\Z}^\lambda L^p_t  L^q_{x}  }\les \|w_1  2^{js}
P_j   u\|_{ \ell_{j\in\Z}^\lambda L^{p_1}_t  L^q_{x}  } \|w_2 G(u)\|_{L^{p_2}_t L^\infty_{x}}\ .\eneq
\end{thm}
As a comparison, we recall here that the estimates obtained from 
 \cite[Theorem 1.2]{HJLW20p}
 state as follows:
\beeq\label{eq-wLeib7}\|w_1w_{2} D^s F(u)\|_{L^q}\les \|w_1 D^s u\|_{L^{q_{1}}} \|w_{2} G(u)\|_{L^{q_{2}}}\ ,\eneq
\beeq\label{eq-wLeib8}\|w_1 w_{2} D^s F(u)\|_{L^q}\les \|w_1 D^s u\|_{L^q} \|w_2 G(u)\|_{L^\infty}\ .\eneq

\begin{prf}
The proof proceeds similar arguments as the estimates obtained from  \cite[Theorem 1.2]{HJLW20p}.
At first, recall that, by repeating essentially the same argument as in the proof of Taylor \cite[(5.6), page 112]{Tay},
 we can obtain 
\beeq\label{eq-HJLW-1}
|P_j  F(u)(x)|\les  \sum_{k\in\Z}\min(1,2^{k-j}) (\M(P_k u)(x)\M(H)(x)+\M(H P_k u)(x))\ ,\eneq
where $H(x)\equiv G(u(x))$.

By  \eqref{eq-HJLW-1}, we know that
\begin{eqnarray*}
&&\|w_1 w_{2}  2^{js}
P_j  F(u)\|_{
\ell_{j\in\Z}^\lambda L^p_t  L^q_{x}
}
\nonumber\\
&\les &
\|w_1w_{2}  2^{js}  \min(1,2^{k-j}) (\M(P_k u)\M(H)+\M(H P_k u)) \|_{ 
\ell_{j\in\Z}^\lambda L^p_t  L^q_{x}\ell_k^1} 
\\
&\les &
\|w_1w_{2} 2^{ks}  \min(2^{(j-k)s},2^{(k-j)(1-s)}) (\M(P_k u)\M(H)+\M(H P_k u)) \|_{ \ell_{j\in\Z}^r \ell_k^1  L^p_t  L^q_{x}} 
\\
&\les &\|w_1w_{2} 2^{ks}  (\M(P_k u)\M(H)+\M(H P_k u)) \|_{ \ell_{k\in\Z}^r   L^p_t  L^q_{x} } 
\ ,\nonumber
\end{eqnarray*}
where we used Young's inequality with the assumption $s\in (0,1)$
in  the last inequality.

By  applying Minkowski's and H\"older's inequalities  to the last expression we have 
\begin{align*}
\|w_1 w_{2}  2^{js}
P_j  F(u) &\|_{
\ell_{j\in\Z}^\lambda L^p_t  L^q_{x}
} \\
 \lesssim \| w_1w_2  2^{ks} \M& (P_k u) \M (H) \|_{ \ell_{k\in\Z}^\lambda L^p_t  L^q_{x}  } +\|w_1w_22^{ks} \M(HP_k u) \big)\|_{ \ell_{k\in\Z}^\lambda L^p_t  L^q_{x} } \\
   \lesssim 
   \|w_2 \M(H)&\|_{L^{p_{2}}_t L^{q_2}_x}
\| w_1 2^{ks}  \M(P_k u)\|_{\ell_{k\in\Z}^\lambda L^{p_1}_t  L^{q_1}_{x}} 
+
\|w_1w_{2} 2^{ks}
H P_k u \|_{\ell_{k\in\Z}^\lambda L^p_t  L^q_{x}}\ ,
\end{align*}
for any
$q\in (1,\infty)$,
 $p, p_{1}, p_2 \in [1,\infty]$,
 and
$q_{1}, q_{2}\in (1,\infty]$ with 
\eqref{eq-HJLW0}.
The last term in the above we used
weighted  Hardy-Littlewood inequality, for $(w_1 w_2)^{q}\in A_q$ with $q\in (1,\infty)$.

If  $q_2<\infty$, 
recall that we have assumed $w_1^{q_1}\in A_{q_1}$, $w_2^{q_2}\in A_{q_2}$,
applying H\"older's inequality and weighted  Hardy-Littlewood inequality
again, we obtain
$$
\|w_1 w_{2}  2^{js}
P_j  F(u) \|_{
\ell_{j\in\Z}^\lambda L^p_t  L^q_{x}
}    \lesssim 
   \|w_2 H\|_{L^{p_{2}}_t L^{q_2}_x}
\| w_1 2^{ks}  P_k u\|_{
\ell_{k\in\Z}^\lambda L^{p_1}_t  L^{q_1}_{x}
 }
\ ,
$$which gives  the desired inequality.

For the
 the remaining case $q_2=\infty$, 
 a similar argument yields \eqref{eq-wLeib6}, if
 we recall the weighted $L^\infty$ estimate of \cite[(2.17)]{HJLW20p}:
\beeq\label{eq-HJLW-2}\|w^{-1} \M(H)\|_{L^\infty}\les 
\|w^{-1} H\|_{L^\infty}\ ,\forall w\in A_1.\eneq
This completes the proof.
\end{prf}

\subsection{weighted fractional Leibniz rule} 
As closely related and useful result is the weighted fractional Leibniz rule. 
\begin{thm}[Weighted fractional Leibniz rule]\label{thm-wLeib-fromChain}
Let $s>0$,
$q_0, q_1, q_2 \in (1,\infty)$,
 $p_1,  p_2 \in (1,\infty]$,
 $s_{j}\in [1,\infty]$ 
 such that
 $$\frac 1{q_0}=\frac{1}{q_{1}}+\frac{1}{p_{1}}=\frac{1}{q_{2}}+\frac{1}{p_{2}} \ ,\ 
 \frac{1}{s_{0}}=\frac{1}{s_{1}}+\frac{1}{s_{2}}=\frac{1}{s_{3}}+\frac{1}{s_{4}}\ . 
 $$ Suppose the time-independent weight functions satisfy
$w_0=w_1 z_1=w_2 z_2>0$,
$w_j^{q_j}\in A_{q_j}$, $z_j^{p_j}\in A_{p_j}$ when $p_j<\infty$
and $z_j^{-1}\in A_{1}$ when $p_j=\infty$,
then we have
\begin{align}
\label{eq-wLeib10}\|w_0 2^{js}
P_j (uv)\|_{\ell_{j}^\lambda L^{s_{0}}_t  L^{q_{0}}_{x} }&\les \|w_1 2^{js}
P_j  u\|_{\ell_{j}^\lambda L^{s_{1}}_{t} L^{q_{1}}_{x}} \|z_{1} v\|_{L^{s_2}_t L^{p_{1}}_x}
 \\
    &  +\|w_2 2^{js}
P_j  v\|_{\ell_{j}^\lambda L^{s_{3}}_{t} L^{q_{2}}_{x}} \|z_{2} u\|_{L^{s_4}_t L^{p_{2}}_x}
\ ,\nonumber
\end{align}
which yields also (for time-independent functions)
\beeq\label{eq-wLeib9}
\|w_0 2^{js}
P_j (uv)\|_{\ell_{j}^\lambda  L^{q_{0}}_{x} }\les \|w_1 2^{js}
P_j  u\|_{\ell_{j}^\lambda   L^{q_{1}}_{x}} \|z_{1} v\|_{  L^{p_{1}}_x}
+\|w_2 2^{js}
P_j  v\|_{\ell_{j}^\lambda  L^{q_{2}}_{x}} \|z_{2} u\|_{ L^{p_{2}}_x}
\ .\eneq
\end{thm}
We remark that the following
weighted fractional Leibniz rule
\beeq\label{eq-wLeib9'}\|w_0 D^s (uv)\|_{L^{q_0}}\les \|w_1 D^s u\|_{L^{q_{1}}} \|z_{1} v\|_{L^{p_{1}}}
+\|w_2 D^s v\|_{L^{q_{2}}} \|z_{2} u\|_{L^{p_{2}}}
\ \eneq
 with $q_j, p_j\in (1,\infty)$ has been obtained, see
Cruz–Uribe and Naibo
\cite{CN16},
D'Ancona
\cite{Da19}. However, as is clear,
in view of application,
 the results with $p_j=\infty$ seem to be more desirable.
 
\begin{prf}
The proof follows from a standard para-product argument and we present only the proof of \eqref{eq-wLeib10}.
In view of $u=\sum P_j u$, $v=\sum P_j v$, we introduce the para-product and decompose $uv$ as follows
$$T_u v=\sum_{j-k> N}P_k u P_j v,\ 
uv=T_u v+T_v u+R(u,v)
$$
where $N$ is chosen such that $P_k u P_j v$ has spectral localization in the annulus of radius $\sim 2^j$.
The estimates for the $T_v u$ and $T_u v$ are easy:
\begin{eqnarray*}
\|w_02^{js} P_j (T_v u)\|_{ L^{s_{0}}_t  L^{q_{0}}_{x} }&\les & \sum_{|l-j|\les N}
 \|w_1 2^{ls}
P_l  u\|_{L^{s_{1}}_{t} L^{q_{1}}_{x}} \|z_{1} P_{<l-N} v\|_{L^{s_2}_t L^{p_{1}}_x}
 \\
&   \les &
\sum_{|l-j|\les N}
 \|w_1 2^{ls}
P_l  u\|_{L^{s_{1}}_{t} L^{q_{1}}_{x}} \|z_{1} v\|_{L^{s_2}_t L^{p_{1}}_x}\ ,
\end{eqnarray*}
and so is 
\eqref{eq-wLeib10} for $T_u v+T_v u$,
 where we have applied \eqref{eq-HJLW-2} 
when $p_1=\infty$, as well as the facts that
\beeq\|P_{<l-N} v\|_{L^p(wdx)}\le C\|v\|_{L^p(wdx)}, w\in A_p, p\in (1,\infty)\ ,\eneq
and
$|P_{<l-N} v|\les \M (v)$.

It remains to control $R(u,v)=\sum_{|l-k|\le N} P_l u P_k v$, for which we have
$$P_j  R(u,v)=P_j(\sum_{|l-k|\le N, j-k\les N} P_k u P_l v)\ .$$
Then it follows that
\begin{eqnarray*}
\|w_02^{js} P_j R(u,v)\|_{l^\la_j L^{s_{0}}_t  L^{q_{0}}_{x} }&\les &
\| \|w_1 2^{js}
P_k  u\|_{L^{s_{1}}_{t} L^{q_{1}}_{x}} \|z_{1} P_{l} v\|_{L^{s_2}_t L^{p_{1}}_x}
\|_{l_j^\la l^1_{k\ge j-CN}l^1_{|l-k|\le N}}
 \\
&   \les &
 \|w_1 2^{js}
P_k  u 
\|_{l_j^\la l^1_{k\ge j-CN}
L^{s_{1}}_{t} L^{q_{1}}_{x}
}
\|z_{1} v\|_{L^{s_2}_t L^{p_{1}}_x}
 \\
&   \les &
 \|
 2^{(j-k)s}
 w_1 
2^{ks}P_k  u 
\|_{l_j^\la l^1_{k\ge j-CN}
L^{s_{1}}_{t} L^{q_{1}}_{x}
}
\|z_{1} v\|_{L^{s_2}_t L^{p_{1}}_x}
 \\
&   \les &
 \| w_1 
2^{ks}P_k  u 
\|_{l_k^\la
L^{s_{1}}_{t} L^{q_{1}}_{x}
}
\|z_{1} v\|_{L^{s_2}_t L^{p_{1}}_x}\ ,
\end{eqnarray*}
where we used Young's inequality in the last inequality, as well as the assumption $s>0$. This 
completes the proof.
\end{prf}

We shall encounter the following weight functions, which are known to be $A_{p}$ weight functions, 
\cite[Lemma 2.5]{HJLW20p}.
 \begin{lem}\label{thm-Ap0}
Let 
 $w(x)=r^{-1+2\de_1}\<r\>^{-2\de_1-2\de_2}$, with
 $0\le 1-2\de_1\le 1+2\de_2<n$.
Then $w\in A_p(\R^n)$, for any $p\in [1,\infty)$.
\end{lem}

As a corollary of the weighted fractional Leibniz rule, Theorem \ref{thm-wLeib-fromChain}, together with a weighted 
 variant of the trace estimates, 
Proposition \ref{thm-trace-w}, we obtain the following inequality which will be frequently used.
 \begin{prop}\label{thm-wtrace}
Let 
$n\ge 3$, $\mu\in (0, 1)$ and 
$|\theta|\le \frac{n-2}{2}+\mu$.
Then
\beeq\label{eq-wtrace}
\|r^{\frac{1-\mu}2}D^{\theta} (fg)\|_{L^{2}}\les
\|r^{-\frac{1-\mu}2}D^{\theta} f\|_{L^{2}}
\|g\|_{\dot H^{\frac{n-2}2+\mu}}\ ,
\eneq
whenever $f, g$ are either spherically symmetric or first-order derivative of spherically symmetric functions. Moreover, for any $q\in [1,\infty]$ and non-endpoint $\theta$, i.e.,
$|\theta|< \frac{n-2}{2}+\mu$
, we could obtain the following estimates by interpolation
\beeq\label{eq-wtrace2}
\|r^{\frac{1-\mu}2}2^{j\theta} P_j(fg)\|_{\ell^q_j L^2_t L^{2}}\les
\|r^{-\frac{1-\mu}2}2^{j\theta} P_j f\|_{\ell^q_j L^2_t L^{2}}
\|g\|_{L^\infty_t \dot H^{\frac{n-2}2+\mu}}\ .
\eneq
\end{prop}
\begin{prf}
At first, we notice that it suffices to prove the result
with $\theta= \frac{n-2}{2}+\mu$, by duality and complex interpolation,
if we recall the well-known fact that
$r^{\al}\in A_{2}$ if and only if $|\al|<n$, and so
$\|r^{\al/2}D^{\theta} f\|_{L^{2}}\simeq
\|r^{\al/2}2^{j\theta} P_{j}f\|_{\ell^{2}_{j}L^{2}}$.

By  Theorem \ref{thm-wLeib-fromChain} \eqref{eq-wLeib9},
with $\la=q_0=2$, we have
$$\|r^{\frac{1-\mu}2}D^{\theta} (fg)\|_{L^{2}}\les
\|r^{-\frac{1-\mu}2}D^{\theta} f\|_{L^{2}}
\|r^{1-\mu}g\|_{L^\infty}+
\|D^{\theta}g\|_{L^2}
\|r^{\frac{1-\mu}2} f\|_{L^{\infty}}\ ,
$$
provided that
$$
r^{-(1-\mu)/2}, r^{\mu-1}\in A_1, \ 
r^{\pm(1-\mu)}\in A_2\ ,
$$
which is true as $\mu\in (1-n, 1)$.
By the symmetric assumption, we have
$$\|r^{1-\mu} g\|_{L^{\infty}}\les
\|r^{1-\mu} g\|_{L^{\infty}_{r}L^{2}_{\omega}}\les
\|D^\theta g\|_{L^{2}}\ ,$$
$$ \|r^{\frac{1-\mu}2} f\|_{L^{\infty}}
\les
\|r^{\frac{1-\mu}2} f\|_{L^{\infty}_{r}L^{2}_{\omega}}
\les \|r^{-\frac{1-\mu}2}D^{\theta} f\|_{L^{2}}\ ,
$$
where we have used
Lemma \ref{thm-trace} and
Proposition \ref{thm-trace-w} \eqref{eq-trace-w2}.
This gives us
\eqref{eq-wtrace}
with $\theta= \frac{n-2}{2}+\mu$.
 This completes the proof.
\end{prf}

\subsection{Inhomogeneous weight}

We will need the following 
weighted Hardy type estimate
with inhomogeneous weight.
\begin{lem}[Weighted Hardy's inequality]\label{thm-w-Hardy}
Let $0\le \al\le\be<n/2-s$ and $s\ge 0$. Then
$$\|r^{-\al-s}\<r\>^{-\be+\al}   u\|_{L^2} \les 
\|r^{-\al}\<r\>^{-\be+\al}   D^s  u\|_{L^2} \ .$$
\end{lem}
\begin{prf} The proof is inspired by \cite[\S 9.3]{MuSch13}.
By Lemma \ref{thm-Ap0},
the conditions are sufficient to ensure
$$(r^{-\al-s}\<r\>^{-\be+\al})^2\in A_2\ ,\ \forall s\in [0, n/2)\ .$$
By Littlewood-Paley theory, we could reduce the proof to the case of $u=P_k u$, with $k\in \Z$, that is, we want to show uniform boundedness of the following operators on $L^2$
$$T_k=2^{-ks}r^{-\al-s}\<r\>^{-\be+\al}P_k
r^{\al}\<r\>^{\be-\al}\ .$$

It is equivalent to the uniform boundedness of
$T_k^* T_k$, with kernel
$$K_k(x,y)=\int
2^{-2ks}w(x) \phi_k(x-z) |z|^{-2s}w^{-2}(z)
  \phi_k(z-y) w(y)dz
$$
where we set
$w(x)=|x|^{\al}\<x\>^{\be-\al}$, 
$\phi_k(x)=2^{kn}\phi(2^k x)$ with $\phi\in\Sc(\R^n)$.
As $K_k(x,y)=K_k(y,x)$, by Schur's test, we need only to prove the uniform boundedness of
\beeq\label{kernel}K_k(x,y)\in L^\infty_x L^1_y\ .\eneq

We will divide the proof into three cases:
i) $|y|\les |z|$, ii) $|y|\gg |z|$ and $|z|\gg 2^{-k}$,  iii) $|y|\gg |z|$ and $|z|\les 2^{-k}$.

{\noindent\bf Case i) $|y|\les |z|$}.
In this case, we have 
$w(y)\les w(z)$ and so
\beeq\label{kernel1}\int |K_k(x,y)|dy\les
\int
2^{-2ks}w(x) |\phi_k(x-z)| |z|^{-2s}w^{-1}(z)
 dz\ .\eneq
 We consider first the sub-case:
 $|z|\ges 2^{-k}$,
for which we have $|2^k z|^{-2s}\les 1$, as $s\ge 0$.
If
$|x|\les |z|$,
for which we have $w(x)\les w(z)$ and so, 
$$\int |K_k(x,y)|dy
 \les \int
|\phi_k(x-z)|
 dz
 \les 1\ .$$
Else, since $|x|\gg |z|$, we know that
$$|\phi_k(x-z)| \les 2^{kn} \<2^k x\>^{-N}\ ,$$
and thus
$$\int |K_k(x,y)|dy\les
\int w(x) |\phi_k(x-z)| w^{-1}(z) dz
\les  \<2^k x\>^{-N}w(x) 2^{kn}\int_{|x|\gg |z|\ges 2^{-k}}  w^{-1}(z) dz
\ .$$
If $|x|\les 1$,  we have $w(x)\simeq |x|^{\al}$ and
$$ \int_{|x|\gg |z|\ges 2^{-k}} w^{-1}(z) dz
\les  |x|^{n-\al}
\les |x|^n w^{-1}(x)\ ,
$$
else, for $|x|\gg 1$,  $w(x)\simeq |x|^{\be}$ and so
$$ \int_{|x|\gg |z|\ges 2^{-k}} w^{-1}(z) dz
\les  |x|^{n-\be}
\les |x|^n w^{-1}(x)\ .
$$
In conclusion, we get \eqref{kernel} for
 $|x|\gg |z|\ges 2^{-k}$:
 $$\int |K_k(x,y)|dy\les  \<2^k x\>^{-N}(2^k |x|)^{n}\les 1
\ .$$

 We consider then another sub-case:
 $|z|\ll 2^{-k}$,
 for which we have
 $$|\phi_k(x-z)| \les 2^{kn} \<2^k x\>^{-N}\ ,$$
 and we need to control
$I=\int |z|^{-2s}w^{-1}(z)
 dz$.
If $k\ge 0$, we have 
$|z|^{-2s}w^{-1}(z)\simeq |z|^{-2s-\al}$ and so
$I\les 2^{-k(n-2s-\al)}$. Then
 $$\int |K_k(x,y)|dy\les
\int
2^{-2ks}w(x) |\phi_k(x-z)| |z|^{-2s}w^{-1}(z)
 dz
 \les
 2^{k\al} \<2^k x\>^{-N}w(x)\les 1\ .$$
For the other case $k<0$, 
$I\les 1+2^{-k(n-2s-\be)}\les 2^{-k(n-2s-\be)}$
and we have
\eqref{kernel} similarly.

{\noindent\bf Case ii) $|y|\gg |z|$ and $|z|\gg 2^{-k}$}.
At first, if $|z|\gg \max(1, 2^{-k})$, we have
$w(y)\simeq |y|^\be$, and
$$\int |\phi_k(z-y)| w(y) dy
\les \int 2^{kn} \<2^k y\>^{-N} |y|^\be dy
\les 2^{-k\be}\les |z|^\be\simeq w(z)\ .
$$
Thus
 we get \eqref{kernel1}, which has been proven to be bounded.
 
 Else, if $2^{-k}\ll |z|\les \max(1, 2^{-k})$, we have $k>0$, $2^{-k}\ll |z|\les 1$  and
 $w(z)\simeq |z|^\al$. Then
 $$\int |\phi_k(z-y)| w(y) dy\les
 \int_{|y|\ge 1}+\int_{|z|\le |y|\le 1} |\phi_k(z-y)| w(y) dy
 \les 2^{k(n-N)}+2^{-k\al}\les |z|^\al\ ,$$
which also gives us \eqref{kernel1}.
 
{\noindent\bf Case iii) $|y|\gg |z|$ and $|z|\les 2^{-k}$}.
In this case, we have $$|\phi_k(x-z)| \les 2^{kn} \<2^k x\>^{-N}\ ,\ 
|\phi_k(z-y)| \les 2^{kn} \<2^k y\>^{-N}\ .
$$
Consider first the case $k\ge 0$, we have
$w(z)\simeq |z|^{\al}$ and
\begin{eqnarray*}
\int |K_k(x,y)|dy & \les &\int_{|z|\ll |y|}
2^{-2k(s-n)}w(x) \<2^k x\>^{-N}
\<2^k y\>^{-N}
 |z|^{-2s}w^{-2}(z)
w(y)dzdy \\
& \les & 2^{-2k(s-n)} w(x) \<2^k x\>^{-N} \int
\<2^k y\>^{-N}
 |y|^{n-2(s+\al)}
w(y)dy\\
& \les &2^{-2k(s-n)} w(x) \<2^k x\>^{-N} 
 2^{k(2s+\al-2n)}\les w(2^k x)\<2^k x\>^{-N} \les 1\ .
\end{eqnarray*}
For the remaining case of $k<0$, we consider three sub-cases:  $|z|\ge 1$, $|y|\le 1$, and $|z|<1<|y|$.

When $|z|\ge 1$, we get
\begin{eqnarray*}
\int |K_k(x,y)|dy & \les &\int_{1\le |z|\ll |y|}
2^{-2k(s-n)}w(x) \<2^k x\>^{-N}
\<2^k y\>^{-N}
 |z|^{-2(s+\be)}
|y|^\be dzdy \\
& \les & 2^{-2k(s-n)} w(x) \<2^k x\>^{-N} \int
\<2^k y\>^{-N}
 |y|^{n-2s-\be}
 dy\\
& \les &2^{-2k(s-n)} w(x) \<2^k x\>^{-N} 
 2^{k(2s+\be-2n)}\les w(2^k x)\<2^k x\>^{-N} \les 1\ ,
\end{eqnarray*}
where we have used the assumption that
$2(s+\be)< n$.

On the other hand, if
 $|y|\le 1$, we have
$\<2^k y\> \simeq 1$ and so
\begin{eqnarray*}
\int |K_k(x,y)|dy & \les &\int_{ |z|\ll |y|\le 1}
2^{-2k(s-n)}w(x) \<2^k x\>^{-N}
 |z|^{-2(s+\al)}
|y|^\al dzdy \\
& \les & 2^{-2k(s-n)} w(x) \<2^k x\>^{-N} \int_{|y|\le 1}
 |y|^{n-2s-\al}
 dy\\
& \les &2^{-2k(s-n)} w(x) \<2^k x\>^{-N} 
\les 1\ ,
\end{eqnarray*}
where we have used the fact that
$-2(s-n)\ge \be\ge \al$.

Finally, when
$|z|<1<|y|$, we see that
\begin{eqnarray*}
\int |K_k(x,y)|dy & \les &\int_{|z|<1<|y|}
2^{-2k(s-n)}w(x) \<2^k x\>^{-N}\<2^k y\>^{-N}
 |z|^{-2(s+\al)}
|y|^\be dzdy \\
& \les & 2^{-2k(s-n)} w(x) \<2^k x\>^{-N} \int_{|y|\ge 1}\<2^k y\>^{-N}
 |y|^{\be}
 dy\\
& \les &
2^{k(n-2s-\be)} w(x) \<2^k x\>^{-N} 
\les 1\ ,
\end{eqnarray*}
where we have used the assumption that
$n-2s-\be\ge \be\ge \al$ in the last inequality.
This completes the proof.
\end{prf}

Based on Theorems \ref{thm-wLeib0}
and  \ref{thm-wLeib-fromChain},
we could obtain
weighted fractional chain rule with higher regularity. For simplicity and future reference, we present the result with the inhomogeneous weight 
$r^{-\al}\<r\>^{-(\be-\al)}$
as in Lemma \ref{thm-Ap0}.

 \begin{prop}[Weighted fractional chain rule, higher regularity]\label{thm-wchain-high}
Let $\theta\in \R_+$, $k=[\theta]\in [0, n/2)$, $0\le \al\le \be< n/ 2-k$, then we have
$$ \|r^{-\al}\<r\>^{-(\be-\al)}D^{\theta}f(u)\|_{ L^2_{x}}\les_f 
C(\max_{j\le k}\|  r^j \nabla^j u\|_{L^\infty_{x}})
\|r^{-\al}\<r\>^{-(\be-\al)}D^{\theta}u\|_{ L^2_{x}}\ ,$$
for any  $f\in C^\infty$.
\end{prop}
\begin{prf}
Let
$w:=r^{- \al}\<r\>^{- (\be-\al)}$, by Lemma \ref{thm-Ap0},
the assumptions on $\al, \be$ ensure
$$w^2,
r^{ -2k}w^2
\in A_2, \ r^{-j}\in A_1, \forall j\in [0, k]\ .$$
The case with $k=0$ follows directly from  Theorem \ref{thm-wLeib0}. In the following, we assume $k\ge 1$.

Let
$\theta=k+\tau$ with $\tau\in [0, 1)$ and $k\ge 1$,
  we have
\begin{eqnarray*}
\|wD^{\theta}f(u)\|_{ L^2_{x}}& \les &
\|w\nabla^k D^{\tau}f(u)\|_{ L^2_{x}}
 \\
& \les & \sum_{|\sum \be_l|=k, |\be_l|\ge 1}
\|w D^{\tau} (f^{(j)}(u)
\Pi_{l=1}^j \nabla^{\be_l} u)
\|_{ L^2_{x}}\ .
\end{eqnarray*}
For each term, we know from  Theorem \ref{thm-wLeib0}
and Theorem \ref{thm-wLeib-fromChain} that
\begin{eqnarray*}
&&\|w D^{\tau} (f^{(j)}(u)
\Pi_{l=1}^j \nabla^{\be_l} u)
\|_{ L^2_{x}}\\&\les &
\|r^{ -k}w D^{\tau} (f^{(j)}(u)-f^{(j)}(0))\|_{ L^2_{x}}
\|r^k\Pi_{l=1}^j \nabla^{\be_l} u
\|_{ L^\infty_{x}}\\
&&+\sum
C(f, \|u\|_{L^\infty_{x}})
\|r^{ -(k-|\be_{l_0}|)}w
D^\tau \nabla^{\be_{l_0}} u\|_{ L^2_{x}}
 \|r^{k-|\be_{l_0}|}
\Pi_{l\neq l_0} \nabla^{\be_l} u
\|_{ L^\infty_{x}}\\
&\les &
C(f, \|u\|_{L^\infty_{x}})\|r^{-k} w D^{\tau}  u\|_{ L^2_{x}}
\Pi_l \|r^{|\be_l|}
\nabla^{\be_l} u
\|_{ L^\infty_{x}}\\
&&+\sum_{l_0}
C(f, \|u\|_{L^\infty_{x}})
\|w
D^{k+\tau } u\|_{ L^2_{x}}
\Pi_{l\neq l_0} \|r^{|\be_l|}
\nabla^{\be_l} u
\|_{ L^\infty_{x}}\\
&\les&
C(f, 
\max_{j\le k}\|  r^j \nabla^j u\|_{L^\infty_{x}})
\|w
D^\theta u\|_{ L^2_{x}}
\ ,
\end{eqnarray*}
where we have also used the
weighted Hardy's inequality, Lemma \ref{thm-w-Hardy}, in the last two inequalities.
\end{prf}

\begin{lem}[Weighted trace estimate]\label{thm-w-Hardy-trace}
Let $n\ge 2$, $0\le \al\le \be\le (n-1)/2$. Then, for any $p\in [2,\infty)$,
$$\|r^{-\al+(n-1)(\frac 1 2-\frac 1p)}\<r\>^{\al-\be}   \phi\|_{L_r^pL_\omega^2} \les 
\|r^{-\al }\<r\>^{\al-\be } D^{\frac 1 2-\frac 1p}  \phi\|_{L^2} \ .$$
\end{lem}

\begin{prf}
Let $w=r^{-\al }\<r\>^{\al-\be }$ with $w^2\in A_2$.
As before, by interpolation, we need only to prove the endpoint case:
$$\|r^{\frac{n-1 }2} w \phi\|_{L_r^\infty L_\omega^2} \les 
\| w P_j 2^{j/2}  \phi\|_{l^1_j L^2} \ ,$$
which follows from 
$$\|r^{\frac{n-1}2}w \phi\|_{L_r^\infty L_\omega^2}^2 \les 
\|w \phi\|_{ L^2}
\|w \nabla \phi\|_{ L^2}
 \ .$$
The proof is elementary, by observing that $r^{n-1}w^2$ is essentially increasing:
\begin{eqnarray*}
\int_{\Sp^{n-1}} w(R)^2 R^{n-1} \phi(R\omega)^2 d\omega
& = & \int_{\Sp^{n-1}}\int_R^\infty w(R)^2 R^{n-1} \pa_r \phi^2(r\omega )dr d\omega \\
 & \les & \int_{\R^{n}} w^2  |\phi \phi'|dx\\
  & \les & \|w\phi\|_{L^2}  \|w\nabla \phi\|_{L^2}\ ,
\end{eqnarray*}
which completes the proof.
\end{prf}

\section{Morawetz type local energy estimates}\label{sec-3le}
In this section, we present a 
 version of Morawetz type local energy estimates, involving fractional derivatives,  for linear wave equations with small, variable $C^1$ coefficients. It is this version of local energy estimates which makes it possible to decrease the regularity requirement for quasilinear wave equations.

The similar estimates for linear wave equations with small $C^2$ coefficients have been well-known, see Metcalfe-Tataru \cite{MeTa12MA}. 
There the authors employ the 
paradifferential calculus and
positive commutator method to obtain a microlocal version of the local energy estimates.
However, it is well-known that for applications to quasilinear problems, $C^2$ requirement is simply too strong to apply for the problem with low regularity. In particular, in the current setting, we are working with the regularity level $s<2$ (in the most physical related case of $n=3$) and the most we could require is a local energy estimate with $C^{1,\al}$ ($\al\le 1/2$) metric, even in the spherically symmetric case.

Here, we present an approach to yield certain weaker but still strong enough variant of the Morawetz type local energy estimates, which apply
for linear wave equations with small $C^1$ coefficients.
The approach is remarkably simple, which is relied basically on the multiplier method with well-chosen multiplier, and interpolation, without consulting 
paradifferential calculus.
The multiplier method has been well-developed in Metcalfe-Sogge
\cite{MetSo06} and Hidano-Wang-Yokoyama \cite[Section 2]{HWY1} (with more general weights which we will mainly follow), for small perturbations of $\Box$.
As we shall see, 
the Morawetz type local energy estimates we shall use are also closely related with the KSS estimates, which appear first in Keel-Smith-Sogge
\cite{KSS}.

Let $T\in (0, \infty]$,
$S_T=[0,T)\times {\mathbb R}^n$,
 $h^{\alpha \beta} \in C^1(S_T)$ with
$h^{\alpha \beta}=h^{\beta \alpha}$, $0\le \al,\be\le n$ satisfying the following uniform hyperbolic condition
\beeq \de_0 (\de^{jk})<(h^{jk}(t,x))< \de_0^{-1} (\de^{jk}), \ 
h^{00}=-1,\ 
|h^{0j}|\le \de_0^{-1}
\ ,\label{eq-3.1}
\eneq
for some $\de_0\in (0,1)$.
Set $\tilde h^{\al\be}= h^{\al\be}-m^{\al\be}$ where
$m^{\al\be}$ is  the flat Minkowski metric component,
$(m^{\al\be})=Diag(-1, 1, 1,\cdots, 1)$.
 Consider the  linear wave equation with variable coefficients (with the
 summation convention for repeated upper and lower indices)
\begin{equation}
\Box_h u:= (h^{\al\be}(t,x)\partial_\al
  \partial_\be) u = F(t,x)
 \quad
\mbox{in}\, \,(0,T)\times{\mathbb R}^n,
\label{eq-3.2}
\end{equation}
with the initial data
\begin{equation}
u(0,\cdot) = u_0 
 , \quad
\partial_t u (0, \cdot) =u_1 
.
\label{eq-3.3}
\end{equation}

\begin{thm}[Morawetz type estimates]
\label{th-LE0}
Let $n \geq 3$, $\mu \in (0,1)$ and consider the initial value problem \eqref{eq-3.2}-\eqref{eq-3.3} with
$h^{0j}=0$, $h^{jk} \in C^\infty(S_{T_0})$
satisfying \eqref{eq-3.1}
and
 \beeq\label{eq-3.4} 
\|r^{1-\mu} \pa h\|_{L^\infty_{t,x}([0,T_0)\times\R^n)}\le \de_0^{-1}\ .
\eneq
 Then, there exist $\de_1\in (0,\min(\de_0, T_0))$ and $C_0\ge 1$, such that for any   $T\in (0, \de_1]$, we have
\beeq\label{eq-3.5} 
\|\tilde \pa u\|_{X_{T}}
\le C_0 (\|(\nabla u_0, u_1)\|_{L^2({\mathbb R}^n)} 
+   \|  F\|_{X_T^*})
\ ,
\eneq
where
$X_T$ and $X_T^*$ are defined in
\eqref{eq-XT} and \eqref{eq-XT*}.
\end{thm}

\subsection{Morawetz type estimates for small perturbations of fixed background}
We begin with a standard energy estimates.
\begin{lem}[Energy estimates]\label{thm-ener}
For any solutions $u\in C^\infty([0,T), C_0^\infty(\R^n))$ to the 
uniformly hyperbolic
equation \eqref{eq-3.2} in $S_T$.
Let $e^0
=(h^{jk} u_j  u_k-h^{00} (u_t)^2 )/2
\simeq |u_t|^2+|\nabla u|^2
$, then there exist a uniform constant $C$ depending only on $\de$ and $n$ such that we have,
for $E(t)=\int_{\R^n} e^0 dx$, 
\beeq\label{eq-3.6}| \frac{d}{dt} E(t)|\le C
\int_{\R^n}( |F| |u_t |+ |\pa h| |e^0|) dx\ .
\eneq
\end{lem}
\begin{prf} The result is classical, which
follows from a multiplier argument:
\begin{eqnarray*}
&&
u_t (h^{\al\be}(t,x)\pa_\al\pa_\be)u  
\\
 & = & \pa_\al (h^{\al\be} u_\be u_t)-\pa_\al (h^{\al\be} )u_\be u_t-(\pa_\al \pt u) h^{\al\be} u_\be \\
 &=& \pa_\al (h^{\al\be} u_\be u_t)
 -\pa_\al (h^{\al\be} )u_\be u_t
  -\pt (\frac{ h^{\al\be} u_\be u_\al}2) + \frac{1}2 (\pt  h^{\al\be}) u_\be u_\al 
 \\
 &=&\pa_\al E^\al+R\ ,
\end{eqnarray*}
where 
$$E^0=-\frac{ h^{\al\be} u_\al  u_\be}2+h^{0\be} u_\be u_t
=\frac{ h^{00} (u_t)^2
-h^{jk} u_j  u_k
}2
,\ E^j=  h^{j\be} u_\be u_t\ ,$$
and
$R= -\pa_\al (h^{\al\be} )u_\be u_t+  (\pt  h^{\al\be}) u_\al  u_\be/2$.
We observe that we have $e^0=-E^0\simeq u_t^2+|\nabla u|^2$, as long as
we have \eqref{eq-3.1}, which gives us
\eqref{eq-3.6} in view of the divergence theorem.
\end{prf}

\begin{lem}[Morawetz type estimates, for small perturbation of Minkowski]
\label{th-3.3}
Let $n \geq 3$, $\mu \in (0,1)$ and consider the initial value problem \eqref{eq-3.2}-\eqref{eq-3.3} for
 $h^{\alpha \beta} \in C^\infty(S_T)$
satisfying \eqref{eq-3.1}.
 Then, there exist $\de\in (0,\de_0)$ and $C\ge 1$, such that for any   $T>0$
 with
  \beeq\label{eq-3.7} 
\|r^{1-\mu} \pa h^{\al\be}\|_{L^\infty_{t,x}([0,T]\times\R^n)}\le\de T^{-\mu}
\ ,\ \|\tilde h^{\al\be}\|_{L^\infty(S_T)}\le \de
\ ,
\eneq
 we have
\beeq\label{eq-3.8} 
\|\tilde \pa u\|_{X_{T}}
\le C  (\|(\nabla u_0, u_1)\|_{L^2({\mathbb R}^n)} 
+   \|  F\|_{X_T^*})
\ .
\eneq
\end{lem}

To prove this result, we need the following fundamental Morawetz type estimates, which follows from the elementary multiplier approach, with carefully chosen multipliers. We leave the tedious proof to the appendix.
\begin{thm}[Morawetz type estimates,  multiplier approach]
\label{th-LE}
Let $n \geq 3$, $\mu \in (0,1)$ and consider the initial value problem \eqref{eq-3.2}-\eqref{eq-3.3} for
 $h^{\alpha \beta} \in C^\infty(S_T)$
satisfying the condition \eqref{eq-3.1}. 
 Then
 there exists $C\ge 1$, which is independent of $T\in (0,\infty)$, such that we have
\begin{align}
\label{eb4-vari}
\|\tilde \pa u\|_{X_{T}}^2
&\le C \|(\nabla u_0, u_1)\|^2_{L^2({\mathbb R}^n)} 
\\&
+ C \int_0^T\int_{{\mathbb R}^n}
|\tilde\partial u|\left(|F|
+|\pa u| \left(|\partial h|+ \frac{|\tilde  h|}{r^{1-\mu}(r+T)^{\mu}}
\right)\right) dxdt
\ ,
\nonumber
\end{align}
for any solutions $u \in
 C^\infty([0,T), C_0^\infty(\R^n))$ to \eqref{eq-3.2}-\eqref{eq-3.3},
with $F \in
 C^\infty([0,T), C_0^\infty)$.
 In addition, for $T\in (0,\infty]$,   we have
\begin{eqnarray}
&\|u\|_{LE_{T}}^2\leq& 
C \|(\nabla u_0, u_1)\|^2_{L^2({\mathbb R}^n)} 
\label{eb4} \\
&  & 
+ C  \int_0^T\int_{{\mathbb R}^n}
|\tilde\partial u|\left(|F|
+|\pa u| \left(|\partial h|+ \frac{|\tilde h|}{r^{1-\mu}\langle r\rangle^{\mu}}
\right)\right) dxdt\ .\nonumber
 \end{eqnarray}
\end{thm}

With the help of \eqref{eb4-vari} and the Cauchy-Schwarz inequality, Lemma \ref{th-3.3} follows directly from the assumption
$$|\partial h|+ \frac{|\tilde  h|}{r^{1-\mu}(r+T)^{\mu}}\ll r^{\mu-1}T^{-\mu}\ .$$

\subsection{Proof of Theorem \ref{th-LE0}  }
With the help of  Lemma \ref{thm-ener} and Lemma \ref{th-3.3}, we are 
ready to present the proof of Theorem \ref{th-LE0}.

Let $\de_1>0$ to be determined.
At first,
without loss of generality, we may assume that
the speed of propagation does not exceed $\de_0^{-1}$, and then
for any $x_0\in \R^n\backslash B_{\de_1}$, the solution $u$ in $$\Lambda_{\de_1}(x_0)=\{(t,x): t\in [0, \de_1],
|x-x_0|<2 \de_1+\de_0^{-1}(\de_1-t)\},
$$ depends only on
$h$, $F$ in $\La_{\de_1}(x_0)$, and
 the data in
$B_{(2+\de_0^{-1})\de_1}(x_0)$.

To apply Lemma \ref{th-3.3}, we need the estimate of perturbation, in $\La_{\de_1}(x_0)$. 
Let $x(s)=x_0+s(x-x_0)$ with $s\in [0,1]$, we have either $\inf |x(s)|\ge |x_0|/2$ or $\inf |x(s)|\le |x_0|/2$. In the second case, there exists $s_0\in [0, 1]$ such that
$|x(s_0)|=\inf |x(s)|\le |x_0|/2$. Then we have
$x-x_0\bot x(s_0)$,
$|x-x_0|\ge |x-x(s_0)|\ge |x_0|/2$
 and $$|x(s)|^2=|x(s_0)|^2+(s-s_0)^2|x-x_0|^2
\ge (s-s_0)^2|x-x_0|^2\ge  \frac{(s-s_0)^2|x_0|^2}4\ .$$
Notice that for the first case, we also have
$|x(s)|\ge |x_0|/2\ge s |x_0|/2$, and we see that in either cases, we have
\beeq\label{eq-lowerofx} |x(s)|\ge   \frac{|s-s_0| |x_0|}2\eneq
for some $s_0\in [0, 1]$.

In view of \eqref{eq-3.4} and \eqref{eq-lowerofx}, the perturbation of $h$ in 
$\La_{\de_1}(x_0)$ could be controlled as follows
\begin{eqnarray*}
|h(t,x)-h(0,x_0)|&\le  &
|\int_0^1 \nabla h(t,x(s))\cdot (x-x_0)ds|+
t  \|\pt h(s, x_0)\|_{L^\infty(s\in [0,t])}
\\&\le&
(|x-x_0|\int_0^1 |x(s)|^{\mu-1}
+
t |x_0|^{\mu-1} )
\|r^{1-\mu} \pa h\|_{L^\infty_{t,x}([0,\de_1]\times\R^n)}
\\&\les&
(|t|+|x-x_0|) |x_0|^{\mu-1}
\|r^{1-\mu} \pa h\|_{L^\infty_{t,x}([0,\de_1]\times\R^n)} \\
&\les& \de_0^{-1}\de_1^\mu\|r^{1-\mu} \pa h\|_{L^\infty_{t,x}([0,\de_1]\times\R^n)} 
\ .
\end{eqnarray*}
Thus, 
$h^{\al\be}$ could be viewed as a small perturbation of $h^{\al\be}(0,x_0)$, in
$\La_{\de_1}(x_0)$, when
 $\de_1\ll 1$. If
 $h^{\al\be}(0,x_0)=m^{\al\be}$,
we could apply Lemma \ref{th-3.3} in $\La_{\de_1}(x_0)$.

In general,
as $h^{jk}$ are uniform elliptic,
there exists a linear transform
$M: \R^n\to  \R^n$ so that in the new coordinates
$h^{jk}(0, x_0)$ reduces to the Euclidean metric.
Suppose that in the new coordinates, $y=M x$,
we have $$H^{\al\be}(t,y)\pa_\al\pa_\be u=F$$
$H^{00}=-1$, $H^{0j}=0$,
$H^{jk}(t,Mx_0)=\de^{jk}$.
Notice that there exists an uniform $C>0$ such that
$$
\|r^{1-\mu} \pa H\|_{L^\infty_{t,y}([0,T]\times\R^n)} \le C
\|r^{1-\mu} \pa h\|_{L^\infty_{t,x}([0,T]\times\R^n)} \ .$$
Thus, when $\de_1\ll 1$, we have
the following variant of 
\eqref{eq-3.7} with $T\le \de_1$,
 \beeq
\|r^{1-\mu} \pa H^{\al\be}\|_{L^\infty_{t,y}([0,\de_1]\times\R^n)}\le\de \de_1^{-\mu}
\ ,\ \|H^{\al\be}-m^{\al\be}\|_{L^\infty(M\La_{\de_1}(x_0))}\le \de
\ ,
\eneq
from which we conclude
 \eqref{eq-3.8} 
in $M\La_{\de_1}(x_0)$,
 with $T\le \de_1$,
 from Lemma \ref{th-3.3}.
Transforming back to the original variable, we obtain for some uniform $C>0$,
\beeq\label{eq-3.12} 
\|\chi_{\La_{\de_1}(x_0)}
\tilde \pa u\|_{X_{T}}
\le C  (\|(\nabla u_0, u_1)\|_{L^2(B_{(2+\de_0^{-1})\de_1}(x_0))} 
+   \|\chi_{\La_{\de_1}(x_0)}  F\|_{X_T^*})
\ ,
\eneq
 for any $|x_0|\ge \de_1\ge T$.

Finally, we choose $x_0\in\{z_j\}_{j=1}^\infty$ so that
  $\cup_j \La_{\de_1}(z_j)=S_{\de_1}$ while $\La_{\de_1}(z_j)$ satisfy finite overlapping property.
Thus we conclude
\eqref{eq-3.5} from 
\eqref{eq-3.12}.

\subsection{Local energy estimates with fractional regularity, $\theta\in [0, 1]$}
Based on  Theorem \ref{th-LE0}, we obtain the following local energy estimates with fractional regularity.
\begin{prop}[Local energy estimates with positive regularity]
\label{th-LEkey1}
Let $n\ge 3$, $\mu\in (0,1)$, $h\in C^1$
with  
$h^{0j}=0$,
 \eqref{eq-3.1}  and
\eqref{eq-3.4}.
Then there exist
$\de_2\in (0,\de_1]$ and a constant $C_1>4 C_0$,
  such that for any $T\in (0,\de_1]$ such that
   \beeq\label{eq-3.4'} 
T^\mu \|r^{1-\mu} \pa h\|_{L^\infty_{t,x}([0,T]\times\R^n)}\le \de_2\ ,
\eneq
   and
  solutions to \eqref{eq-3.2} with data $(u_0,u_1)$, we have
\beeq
\| \tilde \pa D^\theta u\|_{X_T}
\leq
C_1 (\|(\nabla u_0, u_1)\|_{\dot H^{\theta}} +
\|D^\theta F\|_{X_T^*}
\  , \forall \theta\in [0,1],
\label{eb13} 
\eneq
\beeq
\label{eq-LE-crit}
\|\pa D^{1/2} u\|_{X_{T, 1}} 
\le C_1(
\|(\nabla u_0, u_1)\|_{\dot B^{1/2}_{2,1}}     
+T^{\frac{\mu}2}\|r^{\frac{1-\mu}2}2^{j/2}P_j F\|_{\ell_j^1 L^2_{t,x}})
\ .
\eneq
\end{prop}
\begin{prf}
At first, 
by approximation, we could assume $h\in C^\infty$, $u, F\in C^\infty_t C^\infty_0$ so that we could apply  Theorem \ref{th-LE0}.

Let us begin with proving a higher order estimate of \eqref{eq-3.5}.
Applying spatial derivative $\pa_j$ to the equation \eqref{eq-3.2}, we get
\beeq
(- \partial_t^2 + \Delta +  \tilde h^{mk}\partial_m\pa_k )\pa_j u =\pa_j F(t,x)- (\pa_j h^{mk })\partial_m
  \partial_k u \ .
\eneq
By  \eqref{eq-3.5}, we see that
\begin{eqnarray*}
\|\tilde \pa \nabla u\|_{X_T} &\les &
\|u_0\|_{\dot H^2}+\|u_1\|_{\dot H^1}+
\|\nabla F\|_{X_T^*}+
T^{\frac{\mu}2}\|r^{\frac{1-\mu}2}(\nabla h) \nabla^2
u\|_{L^2_{t,x}} \\
  &\les&
  \|u_0\|_{\dot H^2}+\|u_1\|_{\dot H^1}+
\|\nabla F\|_{X_T^*}
+T^{\mu}\|r^{1-\mu}\nabla h\|_{L^\infty}
\|\pa\nabla
u\|_{X_{T}}
 \ .
\end{eqnarray*}
In view of \eqref{eq-3.4'} for some $\de_2\ll 1$, we could absorb the last term and have
\beeq\label{eq-pf1}
\|\tilde \pa \nabla u\|_{X_T} \les
\|u_0\|_{\dot H^2}+\|u_1\|_{\dot H^1}+
\|\nabla F\|_{X_T^*}\ .
\eneq

Notice that
all of the weights occurred in $\|\tilde \pa u\|_{X_T}$ and $X_T^*$
are among the functions
$w=r^{-\frac{1-\mu}2}, r^{-\frac{3-\mu}2}
$
 and their reciprocals,
which share the property that $w^2\in A_2$.
Based on this fact, we know 
\eqref{eq-LPw} holds for $p=2$, 
and so is the complex interpolation satisfied by the weighted Sobolev space of fractional order (see e.g.  \cite[Theorem 6.4.3]{BL}, \cite[Lemma 4.6]{MW17} for similar results)
\beeq\label{eq-LPw2}[\dot H^0_w, \dot H^1_w]_\theta=\dot H^\theta_w, \theta\in [0,1],
\|f\|_{\dot H^\theta_w}:= \|w D^\theta f\|_{L^2}\ .\eneq

With the help of 
\eqref{eq-LPw}, we see that
 \eqref{eq-pf1} gives us
\eqref{eb13}  with $\theta=1$.
As \eqref{eq-3.5} is just \eqref{eb13}  with $\theta=0$,
the general estimate \eqref{eb13} with $\theta\in [0,1]$ follows, in view of
\eqref{eq-LPw2}.
Finally, with the help of the  \eqref{eb13}, \eqref{eq-LPw}, and real interpolation with
$\theta=1/2$, we obtain 
\eqref{eq-LE-crit}.
\end{prf}

From basically the same argument, based on
Theorem \ref{th-LE},
 we could also get the following
local energy estimates with fractional regularity, for small perturbation of Minkowski.
\begin{prop}
\label{th-LEkey2}
Let $n\ge 3$, $\mu\in (0,1)$ and $h\in C^1$
with \eqref{eq-3.1}.
There exists a constant $C>1$ such that if
\beeq\label{eb15} 
\|(r^{1-\mu} \pa h,\tilde  h)\|_{L^\infty_{t,x}([0,T]\times B_1)}
+(\ln\<T\>)
\|(r \pa h, \tilde h)\|_{L^\infty_{t,x}([0,T]\times B_1^c)}
\le \frac{1}{C}\ , 
\eneq
then
for any weak solutions to \eqref{eq-3.2} with data $(u_0, u_1)$, we have
\beeq
\label{eb16} 
\|D^\theta u\|_{LE_{T}}\leq
C  (\|(\nabla u_0, u_1)\|_{\dot H^{\theta}} +
(\ln\<T\>)^{\frac 12}
\|r^{\frac{1-\mu}2}
\<r\>^{\frac{\mu}2}
D^\theta F\|_{L^2_{t,x}} )\  , 
\eneq
for any $\theta\in [0,1]$.
Similarly,
if instead of \eqref{eb15}, we assume  \beeq\label{eb17} 
\|\<r\>^{\mu_1}(r^{1-\mu}\<r\>^{\mu} \pa h, \tilde h)\|_{L^\infty_{t,x}([0,\infty)\times\R^n)}
\le\frac{1}{C}\ ,
\eneq
then we have
\beeq
\|D^\theta \phi\|_{LE}
\leq
C (\|(\nabla u_0, u_1)\|_{\dot H^{\theta}} +
\|r^{\frac{1-\mu}2} \<r\>^{\frac{\mu+\mu_1}2}
D^\theta F\|_{L^2_{t,x}})
\  , \forall \theta\in [0,1]\ .
\label{eb18} 
\eneq
\end{prop}

\subsection{Local energy estimates with negative regularity}
It is well-known that the quasilinear problems endure the issue of loss of regularity, which naturally occurs when we try to prove the convergence  of the Picard iteration series. More precisely, we will need to control some term like $(g(u) - g(v))\Delta v$, for which one standard way to bypass is to prove the convergence in certain weaker topology. One typical choice will be the standard energy norm, for which we are led to the requirement $s\ge 2$ for the regularity. In this sense, to break the regularity barrier $2$ (for dimension three), it is very natural to consider energy type estimates with negative regularity. To obtain such estimates, as we have limited regularity for $h$, it is natural to work for equations in divergence form.
\begin{prop}[Local energy estimates with negative regularity]
\label{th-LEkey3}
Let $n\ge 3$, $\mu\in (0,1/2]$,
 $h^{\al\be}\in C^1$ 
 with  
$h^{0j}=0$,
 \eqref{eq-3.1}
 and
\eqref{eq-3.4}.
Then 
there exist
$\de_3\in (0,\de_2]$ and a constant $C>0$,
  such that for any 
  $T\in (0,\de_1]$
  with
   \beeq\label{eq-3.4''} 
T^\mu \|r^{1-\mu} \pa h\|_{L^\infty_{t,x}([0,T]\times\R^n)}\le \de_3\ ,
\eneq
 we have
 \beeq\label{eb21}
\|\pa D^{-\theta} u\|_{X_T}\leq C \|D^{-\theta} F\|_{X_T^*}\  , \forall \theta\in [0, 1] \  ,  \eneq
for any weak solutions to 
\beeq\label{eq-3.25} (\pa_\al h^{\al\be}
  \partial_\be) u=F \ ,\eneq
 with vanishing data.
Here   $X_T$ and $X_T^*$ are defined in \eqref{eq-XT}-\eqref{eq-XT*}.
In addition, if $
\theta\in [(4-n)/2-\mu, 1]\cap  [0, 1]$,  and $h^{jk}$ are spherically symmetric,
then 
 we have
\beeq\label{eb22} 
\|\pa D^{-\theta} u\|_{X_T}
\les
\| \pa u(0)\|_{\dot H^{-\theta}}+
T^\mu\| h\|_{L^\infty_{t} \dot H^{\frac{n-2}2+\mu}}\| \pa u(0)\|_{\dot H^{1-\theta}}
+ \|D^{-\theta} F\|_{X_T^*}
\eneq
for any spherically symmetric weak solutions to  \eqref{eq-3.25}.
\end{prop}
\begin{rem}\label{thm-rq-s=2}
Here, as is clear from the local energy estimate \eqref{eq-3.5}, when $\theta=0$, 
the second term
on the right of \eqref{eb22} is not necessary. We do not know, however, if it is necessary to have such a term for general $\theta$.
\end{rem}
\begin{prf}
At first, we observe that the local energy estimate \eqref{eq-3.5} applies also for the wave operator in the divergence form
$$  \pa_\al h^{\alpha \beta}(t,x)
  \partial_\beta \ ,$$
  as the difference of these two operators are just a term like $(\pa_\al h^{\alpha \beta})\pa_\be$, which could be absorbed to the left hand by  
\eqref{eq-3.4''} with small $\de_3$, and gives us \eqref{eb21}  with $\theta=0$, which particularly give us
$$
\|  D u\|_{X_T}\les \| F\|_{X_T^*}\  .$$
By duality, we obtain
\beeq\label{eb21-int}
\| \nabla D^{-1}  u\|_{X_T}\les\|   u\|_{X_T}\les \| D^{-1}F\|_{X_T^*}\  .\eneq

By interpolation,
for the proof of \eqref{eb21},
it remains 
 to give the estimate for $\pt u$ with $\theta=1$, for which we shall also argue by duality. 
Observe that the difference
between $(\Box +  \pa_j \tilde h^{jk}
  \partial_k)u$ and $(\Box + \pa_j\pa_k \tilde h^{jk}) u
$ is given by $\pa_j ((\pa_k \tilde h^{jk}) u)=\pa_j ((\pa_k   h^{jk}) u)$, which is an admissible error term thanks to \eqref{eq-3.4}
and
 \eqref{eb21-int}, as we have
$$\|D^{-1} 
\pa_j ((\pa_k h^{jk}) u)\|_{X_T^*}\les
\|(\pa_k h^{jk}) u)\|_{X_T^*}\les
T^\mu \|r^{1-\mu} \pa h\|_{L^\infty_{t,x}}\|u\|_{X_T}
\les \|D^{-1}F\|_{X_T^*}
\ .$$
It is then reduced to the proof of the estimate for $(\Box +  \pa_j\pa_k \tilde h^{jk}) u=F$.
Recall that, for any $G\in C_t^\infty C_0^\infty$ with $\|DG\|_{X^*_T}\le 1$, we have
$$
\|\pa_j\pa_k w\|_{X_T}+\|D\pt w\|_{X_T}\les
\|D\pa w\|_{X_T}\les \|DG\|_{X^*_T}\le 1$$
for any solutions to $(\Box +  \tilde h^{jk}\pa_j\pa_k ) w=G$
with vanishing data on time $t=T$, which follows directly from the estimate 
\eqref{eq-pf1}.
Now, for the purpose of duality, we observe the fact that
$$
\frac{d}{dt}\int_{\R^n} \left(w_t u_t+\nabla w\cdot \nabla u-
 u \tilde h^{jk} \pa_j\pa_k  w\right) dx
+\int_{\R^n}
w_t F+u_t G- u (\pt h^{jk}) \pa_j\pa_k  w 
dx=0\ ,
$$
and so
$$\int_{S_T}  G\pt u dt dx=
\int_{S_T }  ( 
(\pt h^{jk})u\pa_j\pa_k w
-F\pt w
)dt dx\ .
$$ Then, thanks to   \eqref{eq-3.4''}
and
 \eqref{eb21} with $\nabla D^{-\theta} u$,
we obtain that
\begin{eqnarray*}
\left|\int_{S_T }  G\pt u dt dx\right| & \le & \|D\pt w\|_{X_T}\|D^{-1}F\|_{X_T^*}+
T^\mu \|r^{1-\mu} \pa h\|_{L^\infty_{t,x}}\|u\|_{X_T}\|\pa_j\pa_k w\|_{X_T} \\
&\les & \|D^{-1}F\|_{X_T^*}+
T^\mu \|r^{1-\mu} \pa h\|_{L^\infty_{t,x}}\|\nabla D^{-1} u\|_{X_T}
\\
 &  \les &  \|D^{-1}F\|_{X_T^*}\ ,
\end{eqnarray*}
 which, by duality, gives us the desired estimate:
$$\|D^{-1}\pt u\|_{X_T}\les \|D^{-1}F\|_{X_T^*}\ ,$$
and completes the proof of \eqref{eb21}.

Finally, we give the proof of the homogeneous estimates, 
for $( \pa_\al h^{\alpha \beta}
  \partial_\beta)u=0$. For this purpose, we  introduce the homogeneous solution for the standard d'Alembertian $\Box w=0$, $w(0)=u(0)$, $w_t(0)=u_t(0)$. Then it follows from
\eqref{eq-3.5} and $[\Box, D^{\theta}]=0$ that
$$\|D^{\theta} \pa w\|_{X_T}\le  C_0 \|\pa u(0)\|_{\dot H^{\theta}}\ ,$$
for any $\theta\in \R$.
Next, we want to estimate the difference  $v=u-w$, for which we observe that it satisfies $v(0)=\pt v(0)=0$ and
$$ ( \pa_\al h^{\alpha \beta}
  \partial_\beta)v=
  -\pa_j (\tilde h^{jk}
  \partial_k w)\ .$$
Applying \eqref{eb21} for $v$, we obtain that
\begin{eqnarray*}
\|D^{-\theta} \pa v\|_{X_T}
&\les&
T^{\frac{\mu}2}\|r^{\frac{1-\mu}2} D^{-\theta} \pa_j (\tilde h^{jk}
  \partial_k w) \|_{L^2_{t,x}}\\
&\les&
T^{\frac{\mu}2}\|r^{\frac{1-\mu}2} D^{1-\theta} (\tilde h\nabla w)\|_{L^2_{t,x}}
\\
&\les&
T^\mu\| \tilde h\|_{L^\infty_{t} \dot H^{\frac{n-2}2+\mu}}
\|D^{1-\theta} \nabla w\|_{X_T}\\
&\les& T^\mu\| h\|_{L^\infty_{t} \dot H^{\frac{n-2}2+\mu}}
\|\pa u(0)\|_{\dot H^{1-\theta}}
\ ,
\end{eqnarray*}
provided that 
$\theta\in [(4-n)/2-\mu, 1]\cap  [0, 1]$ so that
$|1-\theta|\le (n-2)/2+\mu$ and $\theta\in [0,1]$, where
we have used 
Lemma \ref{thm-wtrace}
in the third inequality, since $h$ and $u$ are assumed to be spherically symmetric. Gluing all 
 these estimates together, we obtain
\eqref{eb22} and this completes the proof.
\end{prf}

\subsection{Local energy estimates with high regularity, radial case}
Considering spherically symmetric equations and solutions, we have
the following version of
local energy estimates with high regularity.

\begin{prop}[Local energy estimates with high regularity]
\label{th-LEkey4}
Let $n\ge 3$, $\mu\in (0,1)$, $h(t,x)=h(t,|x|)\in (\de_0-1, \de^{-1}_0-1)$
and consider radial solutions for
\beeq\label{eq-radial-wave}(-\pt^2+\Delta +h\Delta)\phi=F, \phi(0)=u_0, \phi(0)=u_1\ .\eneq
Then there exists $\de>0$, such that
we have
\beeq\label{eq-LE-high-k}
\|  \tilde\pa D^\theta \phi\|_{X_T}\les
\|D^\theta(\nabla u_0, u_1)\|_{L^2} +
\|D^\theta F\|_{X_T^*}\  , \theta\in[0, [n/2]], 
\eneq
\beeq
\label{eq-LE-high-k-critical}
\|\pa
D^{n/2-1}
\phi\|_{X_{T,1}} 
\les 
\|(\nabla u_0, u_1)\|_{\dot B^{{n/2-1}}_{2,1}}     
+T^{\frac{\mu}2}\|r^{\frac{1-\mu}2}2^{j({n/2-1})}P_j F\|_{\ell_j^1 L^2_{t,x}}
\ ,\eneq
for any classical solutions to \eqref{eq-radial-wave},  provided that
\beeq\label{eq-LE-frac-h}T^\mu\|\pa h\|_{L^\infty_t \dot H^{n/2-1+\mu}}\le \de\ .\eneq
In addition,
when
\eqref{eq-LE-frac-h} is satisfied, for $k=1+[n/2]$, we have
\beeq\label{eq-LE-high-k2}
\|\tilde \pa \nabla^k \phi\|_{X_T}\les
\|\nabla^k(\nabla u_0, u_1)\|_{L^2} +
\|\nabla^k F\|_{X_T^*}
+T^\mu
\|\nabla^k \phi\|_{X_T}
\| h\|_{L^\infty_t \dot H^{\frac{n+2}2+\mu}}
\ , 
\eneq
 if $n$ is odd or $\mu>1/2$, and
\beeq\label{eq-LE-high-k3}
\|\tilde \pa \nabla^k \phi\|_{X_T}\les
\|\nabla^k(\nabla u_0, u_1)\|_{L^2} +
\|\nabla^k F\|_{X_T^*}
+T^{\mu}
 \|h\|_{L^\infty_t \dot H^{k+1}}\|D^{\frac n 2+\mu}\phi\|_{X_T}
\ , 
\eneq
 if 
  $\mu>1/2$.
Similarly, when $n\ge 4$, $\mu\in (0, 1/2)$, $\mu_1\in (0,\mu]$, there exists $\de'>0$, such that
we have
\beeq\label{62eq-LE-high-k}
\| D^\theta \phi\|_{LE_{T}}\les
\|D^\theta(\nabla u_0, u_1)\|_{L^2} +
 \|r^{\frac{1-\mu}2}\<r\>^{\frac{\mu+\mu_1}2} D^\theta F\|_{L^2_{t,x}}
\  , \theta\in[0, [(n-1)/2]], 
\eneq
for any classical solutions to \eqref{eq-radial-wave},  provided that
\beeq\label{62eq-LE-frac-h}
\|h\|_{L^\infty_{t,x}}+
\|\pa h\|_{L^\infty_t \dot H^{n/2-1+\mu}}+\|\pa h\|_{L^\infty_t \dot H^{n/2-1-\mu_1}}\le \de'\ .\eneq
\end{prop}

\begin{prf}
As in Proposition \ref{th-LEkey1},
the estimates
\eqref{eq-LE-high-k} and \eqref{eq-LE-high-k-critical} could be reduced to the proof of
\eqref{eq-LE-high-k}
with $D^\theta$ replaced by $\nabla^k$, with $k\in \N$.
When $\theta=0, 1$, it has been proven from 
\eqref{eb13} of Proposition \ref{th-LEkey1},
by recalling the 
trace
estimates:
$$
\|r^{1-\mu}\pa h\|_{L^\infty_{t,x}}\les \|\pa h\|_{L^\infty_t \dot H^{n/2-1+\mu}}\ .$$

The general case follows from the similar strategy.
By applying $\nabla^\al$ with $|\al|=k\ge 2$, we
have
$$(-\pt^2+\Delta +h\Delta)\nabla^\al\phi=\nabla^\al F
+
[h, \nabla^\al]\Delta \phi
=\nabla^\al F
+\sum_{j=1}^{k}
\mathcal{O}(|\nabla^{j} h| | \nabla^{k-j} \Delta \phi|)
\ .$$
At first, if $1\le j\le k-1$, we have
$n/2+\mu-j\in (1/2, n/2)$, and
\begin{eqnarray*}
\|(\nabla^{j} h)  \nabla^{k-j} \Delta \phi\|_{X_T^*} & \les & T^{\mu/2}
\|r^{j-\mu}\nabla^{j} h\|_{L^\infty}\|r^{\mu/2-j+1/2}\nabla^{k-j} \Delta \phi\|_{L^2} \\
 &\les & T^{\mu/2}
 \|h\|_{L^\infty_t \dot H^{n/2+\mu}}\|r^{-(1-\mu)/2}\nabla^{k-1} \Delta \phi\|_{L^2}\ ,
\end{eqnarray*}
where we have used
\eqref{eq-SteinWeiss} and trace estimate.
For the term with $j=k$, we could proceed similarly, if we have
$n/2+\mu-k\in (1/2, n/2)$, that is,
$k<(n-1)/2+\mu$.
Thus all these commutator terms could be absorbed to the left,
in view of \eqref{eq-LE-frac-h}, which proves
\eqref{eq-LE-high-k} with $0\le k<(n-1)/2+\mu$.

For the remaining case,
$(n-1)/2+\mu\le k\le [n/2]$, it only happens when $n$ is even and $k=n/2$ with $\mu\in (0, 1/2]$. 
Then we have $k-1\ge 1$,
\begin{eqnarray*}
\|(\nabla^{k} h)  \Delta \phi\|_{X_T^*} & \les & T^{ \mu/ 2}
\|r^{-\mu}\nabla^{k} h\|_{L^\infty_t L^2}\|r^{(1+\mu)/2}  \Delta \phi\|_{L^2_t L^{\infty}} \\
 &\les & T^{\mu /2}
 \|h\|_{L^\infty_t \dot H^{k+\mu}}\|r^{-(1-\mu)/2}\nabla^{k-1} \Delta \phi\|_{L^2}\ ,
\end{eqnarray*}
where we have used
\eqref{eq-trace-w'} in Proposition \ref{thm-trace-w} and
Hardy's inequality.
This gives us
\eqref{eq-LE-high-k} with $k=n/2$.

Turning to the proof of \eqref{eq-LE-high-k2}
and \eqref{eq-LE-high-k3}, in which case we have
$(n-1)/2+\mu\le k<(n+1)/2+\mu$.
Notice that 
we still have $n/2+\mu-j\in (1/2, n/2)$ for $1\le j\le k-1$ and, as before, these commutator terms are good terms. For the case $j=k$, we have
$(n+2)/2+\mu-k>1/2$ and so
\begin{eqnarray*}
\|(\nabla^{k} h)  \Delta \phi\|_{X_T^*} & \les & T^{\mu/2}
\|r^{1-\mu+k-2}\nabla^{k} h\|_{L^\infty_{t,x} }\|r^{-(1-\mu)/2+2-k}  \Delta \phi\|_{L^2_{t,x}} \\
 &\les & T^{\mu/2}
 \|h\|_{L^\infty_t \dot H^{(n+2)/2+\mu}}\|r^{-(1-\mu)/2}\nabla^{k}\phi\|_{L^2}
 \\
 &\les & T^{\mu}
 \|h\|_{L^\infty_t \dot H^{(n+2)/2+\mu}}\|\nabla^{k}\phi\|_{X_T}\ ,
\end{eqnarray*}
which gives us  \eqref{eq-LE-high-k2}.
Similarly, for \eqref{eq-LE-high-k3},
 the term with $j=k$  could be controlled as follows 
\begin{eqnarray*}
\|(\nabla^{k} h)  \Delta \phi\|_{X_T^*} & \les & T^{\mu/2}
\|r^{n/2-1}\nabla^{k} h\|_{L^\infty_{t,x} }\|r^{(1-\mu)/2+1-n/2
}  \Delta \phi\|_{L^2_{t,x}} \\
 &\les & T^{\mu/2}
 \|h\|_{L^\infty_t \dot H^{k+1}}\|r^{-(1-\mu)/2}D^{n/2+\mu}\phi\|_{L^2}
 \\
 &\les & T^{\mu}
 \|h\|_{L^\infty_t \dot H^{k+1}}\|D^{n/2+\mu}\phi\|_{X_T}\ ,
\end{eqnarray*}
where we used
$n/2+\mu\ge 2$ and 
\eqref{eq-SteinWeiss},
which completes the proof of  \eqref{eq-LE-high-k3}.

Finally, we treat
 \eqref{62eq-LE-high-k}, for which we follow the similar strategy, by reducing it to $\nabla^k$ with $k= [(n-1)/2]$.
At first, 
for $1\le j\le [n/2]-1$,
 we notice that
$
n/2-\mu_1-j,
n/2+\mu-j\in (1/2, n/2)$, and so
\begin{eqnarray*}
&&\|r^{(1-\mu)/2}\<r\>^{(\mu+\mu_1)/2} (\nabla^{j} h)  \nabla^{k-j} \Delta \phi\|_{L^2_{t,x}}\\
& \les & 
\|r^{j-\mu}\<r\>^{\mu+\mu_1}\nabla^{j} h\|_{L^\infty}\|r^{-(1-\mu)/2-(j-1)}\<r\>^{-(\mu+\mu_1)/2}\nabla^{k-j} \Delta \phi\|_{L^2} \\
 &\les &
 \|h\|_{L^\infty_t \dot H^{n/2+\mu}\cap L^\infty_t \dot H^{n/2-\mu_1}}\|r^{-(1-\mu)/2}\<r\>^{-(\mu+\mu_1)/2}\nabla^{k-1} \Delta \phi\|_{L^2}\ ,
\end{eqnarray*}
where we have used
Lemma 
\ref{thm-w-Hardy}.
For the remaining terms with $j>[n/2]-1$, we see that $n$ is odd, $j=k=(n-1)/2$ and so
$n/2+ \mu-k=1/2+ \mu$,
$n/2- \mu_1-k=1/2- \mu_1$.
Notice that
\eqref{eq-trace2-sd} gives us that
$$\|r^{(n-1)(\frac 1 2-\mu_1)}h\|_{L^{\frac 1 {\mu_1}}}\les
\|h\|_{\dot H^{\frac1 2-\mu_1}},
\|r^{(n-1)(\frac 1 2-\mu_1)-\mu-\mu_1}h\|_{L^{\frac 1{\mu_1}}}\les
\|h\|_{\dot H^{\frac1 2+\mu}},
$$
that is,
\beeq\label{eq-tracemu1}
\|r^{(n-1)(\frac 1 2-\mu_1)-\mu-\mu_1}\<r\>^{\mu+\mu_1}h\|_{L^{\frac1{\mu_1}}}\les
\|h\|_{\dot H^{\frac 1 2+\mu}}
+\|h\|_{\dot H^{\frac 1 2-\mu_1}}.
\eneq
With the help of \eqref{eq-tracemu1}, we obtain, for $1/q=1/2-\mu_1$,
\begin{eqnarray*}
&&\|r^{\frac{1-\mu}2}\<r\>^{\frac{\mu+\mu_1}2} (\nabla^{k} h) \Delta \phi\|_{L^2_{t,x}}\\
& \les & 
\|r^{(n-1)(\frac 1 2-\mu_1)-\mu-\mu_1}\<r\>^{\mu+\mu_1}\nabla^{k} h\|_{L^\infty_t L^{\frac 1{\mu_1}}}\|
r^{-\frac{1-\mu}2 -\frac{n-3}2+ n \mu_1}
\<r\>^{-\frac{\mu+\mu_1}2} \Delta \phi\|_{L^q} \\
 &\les &
 \|h\|_{L^\infty_t \dot H^{\frac n 2+\mu}\cap L^\infty_t \dot H^{\frac n 2-\mu_1}}\|r^{-\frac{1-\mu}2
 -\frac{n-3}2+  \mu_1}\<r\>^{-\frac{\mu+\mu_1}2}
D^{\mu_1}
\Delta \phi\|_{L^2}\\
 &\les &
 \|h\|_{L^\infty_t \dot H^{\frac n 2+\mu}\cap L^\infty_t \dot H^{\frac n 2-\mu_1}}\|r^{-\frac{1-\mu}2}\<r\>^{-\frac{\mu+\mu_1}2}D^{\frac{n-3}2}
\Delta \phi\|_{L^2}\ ,
\end{eqnarray*}
where we have used
 Lemma \ref{thm-w-Hardy-trace}, Lemma 
\ref{thm-w-Hardy} and the assumption $n\ge 4$ so that we have $(n-3)/2\ge \mu_1$.
This gives us \eqref{62eq-LE-high-k}. 
\end{prf}

\section{Local existence and uniqueness for dimension three}\label{sec-4}
With the help of Propositions \ref{th-LEkey1} and \ref{th-LEkey3}, we are able to prove the local existence and uniqueness part of
Theorem \ref{th-1}.

\subsection{Approximate solutions}\label{sec-approx}
Firstly, we fix a 
spherically symmetric function
$\rho \in C_0^{\infty}({\mathbb R}^n)$ which equals $1$ near the origin and $\int_{{\mathbb R}^3} \rho(x)dx = 1$. and set
$\rho_{k} (x) = 2^{3k} \rho (2^{k}x)$. 
Based on $\rho$, we define standard sequence of $C^\infty$,  spherically symmetric, approximate functions to $(u_0, u_1)$,
\begin{equation}
 u_0^{(k)}(x) = \rho_{k} * u_0(x),
 u_1^{(k)}(x) = \rho_{k} * u_1(x), k\ge 3\ .
\label{eq-4.1}
\end{equation}
As is clear, we know that
$$\|(\nabla u_0^{(k)}, u_1^{(k)})\|_{\Hsd{\theta }}\le
\|(\nabla u_0, u_1)\|_{\Hsd{\theta }}, \forall \theta \in \R\ ,$$
$$\|(\nabla u_0^{(k)},u_1^{(k)})\|_{\dot B^{\theta }_{2,1}}\le
\|(\nabla u_0,u_1)\|_{\dot B^{\theta }_{2,1}}, \forall \theta \in \R\ .$$
Since $(u_0, u_1)\in \Hs{s} \times (\Hs{s-1}
\cap \Hsd{s_0-1})$ with $s\in (3/2, 2]$ and
$s_0\in [2-s, s-1]$,
 we have
$$
\lim_{k\to\infty}(\| u_0^{(k)} - u_0\|_{\Hs{s}} +\| u_1^{(k)} -u_1\|_{\Hs{s-1}\cap \Hsd{s_0-1}}) = 0\ .$$
In addition, for any $\theta\in [s_0-3, s-1)$,
we know that
\begin{equation}
 \sum_{k = 3}^{\infty} \left( \|\nabla  u_0^{(k)}
 - \nabla  u_0^{(k+1)}\|_{\dot H^{\theta}({\mathbb R}^3)}
+ \| u_1^{(k)}
 -  u_1^{(k+1)}\|_{\dot H^{\theta}({\mathbb R}^3)}\right)
< \infty.
\label{eq-4.2}
\end{equation}
Indeed, we can easily check this property by using the fact that
$\|\rho_{k} * \varphi - \varphi\|_{L^2} \leq C 2^{-\theta k}\|\varphi\|_{\dot H^\theta}$ for any $\theta\in [0, 2]$.
Moreover, 
there exists subsequence $\{j_k\}$, so that
\beeq \label{eq-4.2'} \|  u_0^{(j_k)}
 -   u_0^{(j_{k+1})}\|_{\dot H^{s}({\mathbb R}^3)}
+ \| u_1^{(j_k)}
 -  u_1^{(j_{k+1})}\|_{\dot H^{s-1}({\mathbb R}^3)} \le
 2^{-k}\ ,\eneq
and we also have \eqref{eq-4.2} for $(u_0^{(j_k)}, u_1^{(j_k)})$ with $\theta=s-1$. 
Furthermore, we could cut-off the data so that they are compactly supported smooth functions, while all of these properties remain valid (with possible augment of the constants). We 
still denote the sequence (after cut-off) as $(u_0^{(j_k)}, u_1^{(j_k)})$.

With $( u_0^{(j_k)},  u_1^{(j_k)})$ as data, we use a standard iteration to define the sequence of approximate solutions. Let 
$F( u)=  a(u) u_t^2+b(u) |\nabla u|^2$,
$u_{2} \equiv 0$ and define $u_k$ ($k\ge 3$) recursively by solving
\beeq
\left\{\begin{array}{l}\Box u_k  + g(u_{k-1})
  \Delta u_k =  F( u_{k-1}), (t,x)\in (0,T)\times \R^3,
 \\
u_k(0,\cdot) =  u_0^{(j_k)},\quad
 \partial_t u_k(0,\cdot) =  u_1^{(j_k)}.\end{array}
 \right.
 \label{eq-4.3}
\eneq

By  Proposition \ref{th-LEkey1}, 
together with a standard existence, uniqueness and regularity theorem,
 we will see that, 
there exists some uniform $T(u_0, u_1)\in (0, \infty)$ to be determined, so that,
 for all $k\ge 2$,
$u_k$ is well defined,
spherically  symmetric, and satisfies $u_k \in C^{\infty}(S_T)$
\begin{equation}\|\pa u_k\|_{L^\infty \dot H^{\theta}}
\le
\|\pa D^\theta u_k\|_{X_{T}}\le 2C_1 \|(\nabla u_0^{(j_k)}, u_1^{(j_k)})\|_{\dot H^{\theta}}\ , \forall \theta\in  [0,1]
\ ,
  \label{eq-4.4'}
\end{equation}
\begin{equation}  \label{eq-4.4Besovbd}
\| \pa u_k\|_{L^\infty_t \dot B^{1/2}_{2,1}}
\le\|\tilde\pa D^{1/2} u_k\|_{X_{T, 1}}
\le 2C_1 \|(\nabla u_0, u_1)\|_{\dot B^{1/2}_{2,1}}=2C_1\ep_c\ .
\end{equation}

\subsection{Uniform boundedness of $u_k$}
In this subsection, we prove the uniform boundedness of the sequence, \eqref{eq-4.4'} and \eqref{eq-4.4Besovbd}.

\begin{lem}
\label{thm-bdd0}
Let $s\in (3/2, 2]$,
$\ep_s:= \|(\nabla_x u_0,u_1 )\|_{\dot H^{s-1}}$ and set $s=3/2+\mu$.
Then there exists
$c=c(g,a,b,\ep_c)$ such that
 the spherically symmetric functions $u_k\in C^\infty\cap C H^\theta\cap C^1 H^{\theta-1}$ are well-defined on $S_T$ for any $k\ge 2$, $\theta\ge 3$ and enjoy the  uniform bounds \eqref{eq-4.4'} and \eqref{eq-4.4Besovbd}, for any $T\in (0, T_0]$
 with
 $T_0=\min (\de_1,
c(g,a,b,\ep_c)\ep_s^{-1/\mu}))$.
\end{lem}
\begin{prf}
The proof proceeds by induction. At first, the result is trivial for $k=2$. Then
 we make the inductive assumption that for some $m\ge 3$, we have
for any $2\le k\le m-1$, 
$u_k\in C^\infty\cap C H^{\theta}\cap C^1 H^{{\theta}-1}$ for any ${\theta}\ge 3$ 
with the bounds \eqref{eq-4.4'}-\eqref{eq-4.4Besovbd} satisfied.

Recall the Sobolev inequality 
$$\|\phi\|_{L^\infty}\le C\|\phi\|_{ \dot B^{3/2}_{2,1}}
\ ,$$
and in view of \eqref{eq-4.4Besovbd} for $u_{m-1}$,
we see that
\beeq\label{eq-Sob}\|u_{m-1}\|_{L^\infty_{t,x}}\le C\|u_{m-1}\|_{ L^\infty_t\dot B^{3/2}_{2,1}}
\le 2CC_1
\|(\nabla u_0, u_1)\|_{\dot B^{1/2}_{2,1}}
=2CC_1\ep_c\ .
\eneq

As
$u_{m-1}\in C^\infty\cap C H^{\theta}\cap C^1 H^{{\theta}-1}$ for any ${\theta}\ge 3$ 
with the bounds \eqref{eq-4.4'} and \eqref{eq-4.4Besovbd},
we see that $F(u_{m-1})\in L^1([0,T]; H^{\theta-1})$ and $g(u_{m-1})\in C^\infty$.
Based on this information, we see from the classical local existence throrem that 
the equation \eqref{eq-4.3} is
 solvable with solution $u_m$ well-defined, smooth in $[0,T]\times \R^n$ and
  $u_m\in C H^{\theta}\cap C^1 H^{{\theta}-1}$ for any ${\theta}\ge 3$.

To
apply Proposition \ref{th-LEkey1} for $u_m$,
we need to check \eqref{eq-3.4'}for $h^{jj}=g(u_{m-1})$ and $h^{\al\be}=0$ with $\al\neq \be$.
 As $\mu\in (0, 1/2]$ and $u_{m-1}$ is spherically symmetric,
by the inequality \eqref{eq-2.1},
we have
\beeq\label{eq-4.5'}\|r^{1-\mu}\pa g(u_{m-1})\|_{L^{\infty}}\le
\|r^{1-\mu} g'(u_{m-1}) \pa u_{m-1}\|_{L^{\infty}}\les
 \|\pa u_{m-1}\|_{L^{\infty} \Hsd{1/2+\mu}}
\les
\ep_s
\eneq
where we have used 
\eqref{eq-4.4'} and \eqref{eq-Sob} for $u_{m-1}$. Here we notice that the implicit constant may depend on $g$ and $\ep_c$ through $\|g'(u_{m-1})\|_{L^\infty}$.
Thus,
with
\beeq\label{eq-life-poly}T_0=
c(g,\ep_c)\ep_s^{-1/\mu}
\ ,\eneq for some small constant $c$,
which may depend on $\ep_c$ and $g$,
 we have
\eqref{eq-3.4'} for $g(u_{m-1})$ and could 
apply Proposition \ref{th-LEkey1} for $u_m$. In conclusion,
we get for $\theta\in [0,1]$ and $T\in (0,\de_1]$,
\beeq
\label{} 
\|\pa D^\theta u_m\|_{X_{T}}\le C_1(
\|(\nabla u_0^{(j_m)}, u_1^{(j_m)})\|_{\dot H^{\theta}} +
T^{\frac \mu 2}
\|r^{\frac{1-\mu}2}
D^{\theta} F(u_{m-1})\|_{L^2_{t,x}}) \  , 
\eneq
and 
\beeq
\label{} 
\|\pa D^{1/2} u_m\|_{X_{T, 1}}
\le C_1(
\|(\nabla u_0, u_1)\|_{\dot B^{1/2}_{2,1}}
+T^{\frac{\mu}2}\|r^{\frac{1-\mu}2}2^{j/2}P_j F(u_{m-1})\|_{\ell_j^1 L^2_{t,x}})
 \  .
\eneq

To control the right hand side, we will exploit Lemma \ref{thm-wtrace}, the weighted fractional chain rule, Theorem \ref{thm-wLeib0}, as well as the
weighted fractional Leibniz rule, Theorem \ref{thm-wLeib-fromChain}. Without loss of generality, we assume $b=0$ and write \beeq\label{eq-Fuexp}F(u)=  a(u) u_t^2=
\tilde a(u) u_t^2+a(0)u_t^2:=F_1(u)+F_2(u)\ .\eneq
We first give the estimate for $F_2(u)$ for any $\theta\in [0, 1]$:
\begin{eqnarray*}
\|r^{\frac{1-\mu}2}2^{j\theta}P_j F_2(u)\|_{\ell_j^q L^2_{t,x}} & \les & 
\|r^{-\frac{1-\mu}2}2^{j\theta}P_j u_t\|_{\ell_j^q L^2_{t,x}}
\|r^{1-\mu} u_t\|_{L^\infty_{t,x}} \\
 & \les & 
 \|r^{-\frac{1-\mu}2}2^{j\theta}P_j u_t\|_{\ell_j^q L^2_{t,x}}
\| u_t\|_{L^\infty_{t} \dot H^{1/2+\mu}}
\end{eqnarray*}
by
 Lemma 
\ref{thm-Ap0}, \eqref{eq-wLeib10} and \eqref{eq-2.1}.
For the first term $F_1(u)$, we use similar argument, followed by
\eqref{eq-wLeib10}, 
\eqref{eq-wLeib6}, and  \eqref{eq-2.1}
 to obtain for $\theta\in [0,1]$,
\begin{eqnarray*}
&&\|r^{\frac{1-\mu}2}2^{j\theta}P_j F_1(u)\|_{\ell_j^q L^2_{t,x}} \\
& \les & 
\|\tilde a(u) \|_{L^\infty}\|r^{\frac{1-\mu}2}2^{j\theta}P_j (u_t^2)\|_{\ell_j^q L^2_{t,x}} +\|r^{-\frac{3-\mu}2}2^{j\theta}P_j  \tilde a(u) \|_{\ell_j^q L^2_{t,x}}
\|r^{2-\mu} u_t^2\|_{L^\infty}
\\ & \les  & C( \|\pa u\|_{L^\infty_t \dot B^{\frac 1 2}_{2,1}})\| u_t\|_{L^\infty_t \dot H^{\frac 1 2+\mu}}  ( \|r^{-\frac{1-\mu}2}2^{j\theta}P_j u_t\|_{\ell_j^q L^2_{t,x}}+\|r^{-\frac{3-\mu}2}2^{j\theta}P_j u\|_{\ell_j^q L^2_{t,x}})\ .
\end{eqnarray*}
For the last term occurred in the last inequality, we recall that 
we have
$$\|r^{-\frac{3-\mu}2}2^{j\theta}P_j u\|_{\ell_j^2 L^2_{t,x}}
\simeq\|r^{-\frac{3-\mu}2}D^{\theta}  u\|_{L^2_{t,x}}
\les\|r^{-\frac{1-\mu}2}D^{1+\theta}  u\|_{L^2_{t,x}}
\simeq\|r^{-\frac{1-\mu}2}D^{\theta} \nabla u\|_{L^2_{t,x}}
$$
in
the case of $q=2$ and $\theta\in [0,1]$, which
 follows directly from the weighted Hardy-Littlewood-Sobolev inequality. The general result for $q$ with non-endpoint $\theta\in (0,1)$ follows then from the real interpolation, that is, we have
$$\|r^{-\frac{3-\mu}2}2^{j\theta}P_j u\|_{\ell_j^q L^2_{t,x}}
\les
\|r^{-\frac{1-\mu}2}2^{j\theta}P_j \nabla u\|_{\ell_j^q L^2_{t,x}}\ .$$
In conclusion, we have proved that, for $q=2$, $\theta\in [0,1]$ or
 $q=1$, $\theta=1/2$
\beeq\label{eq-Funl}
\|r^{\frac{1-\mu}2}2^{j\theta}P_j F(u)\|_{\ell_j^q L^2_{t,x}}  \le 
 C(a, \|\pa u\|_{L^\infty_t \dot B^{\frac 1 2}_{2,1}})
T^{\frac \mu 2} \|\pa D^\theta u\|_{X_{T,q}}
\| \pa u\|_{L^\infty_{t} \dot H^{\frac 1 2+\mu}}\ ,
\eneq
and from which we get
$$
\|\pa D^\theta u_m\|_{X_{T,q}}
\le 2 C_1
\|(\nabla u_0^{(j_m)}, u_1^{(j_m)})\|_{\dot B^{\theta}_{2,q}}\ ,
$$ provided that we  set $T\le \min(\de_1, T_0)$ where $T_0$ is given in \eqref{eq-life-poly} with possibly smaller $c=c(g,a,b,\ep_c)>0$ such that $
4 C_1^2 C(a, 2C_1 \ep_c) T^\mu\le 1$.
This completes the proof by induction.
\end{prf}

To prove convergence of the approximate solutions, when $s<2$, we will also require bounds for  the  solutions in Sobolev space of negative order.
\begin{lem}
\label{thm-bdd1}
Under the same assumption as in Lemma \ref{thm-bdd0}.
Let $(u_0, u_1) \in \Hs{s}({\mathbb R}^3) \times (\Hs{s-1}
\cap \Hsd{s_0-1})({\mathbb R}^3)$
with $s_0\in [2-s, s-1]$.
Then there exist some $c=c(g, a, b,\ep_c)\in (0,1)$ 
and $C>0$ such that
for any $\theta\in [s_0, s-1]$, we have
 \begin{equation}  \label{eq-4.4'''}
 \|D^{\theta-1}\pa u_k\|_{X_T}\le C\ep_{\theta}
+C T^\mu
\ep_{\theta+1}
\ep_{s-1}
 \ , \forall k\ge 2\ ,
\end{equation}
provided that
\beeq\label{eq-life-poly2}T\le \min (\de_1, c \ep_s^{-1/\mu})\ .\eneq
In particular, with $\theta=s-1$, we have
 \beeq  \label{eq-4.4''}
 \|\pa u_k\|_{L^\infty \dot H^{s-2}}
\le
\|D^{s-2}\pa u_k\|_{X_T}\le 2C\|(\nabla u_0, u_1)\|_{\dot H^{s-2}}\ .
 \eneq
\end{lem}
\begin{prf}
As in the proof of Lemma \ref{thm-bdd0},
we proceed by induction. At first, the result is trivial for $k=2$. Then
 we make the inductive assumption that for some $m\ge 3$, we have
 \eqref{eq-4.4'''} satisfied by $u_k$, for any $2\le k\le m-1$. 

To apply Proposition \ref{th-LEkey3}, we write the equation \eqref{eq-4.3} of $u_m$ in the equivalent divergence form for $(t,x)\in (0,T)\times \R^3$:
\beeq
\left\{\begin{array}{l}\Box u_k + \nabla \cdot (g(u_{k-1})
  \nabla u_k) =
  \nabla ( g(u_{k-1}) )\cdot
  \nabla u_k
+    F( u_{k-1}),
 \\
u_k(0,\cdot) =  u_0^{(j_k)},\quad
 \partial_t u_k(0,\cdot) =  u_1^{(j_k)}.\end{array}
 \right.
 \label{eq-4.3'}
\eneq
As we see from the proof of Lemma \ref{thm-bdd0},
\eqref{eq-Sob} and \eqref{eq-4.5'}, we know that
\eqref{eq-3.4''}  is satisfied and we could apply Proposition \ref{th-LEkey3} to obtain for $
\theta\in [1/2-\mu, 1/2+\mu]$:
\begin{eqnarray*}
&&\|D^{-\theta} \pa u_k\|_{X_T}\\
&\les&
\| \pa u_k(0)\|_{\dot H^{-\theta}}+ T^\mu\|  g(u_{k-1}) \|_{L^\infty_{t} \dot H^{\frac 1 2+\mu}}\| \pa u_k(0)\|_{\dot H^{1-\theta}}
\\
&&
+
T^{\frac{\mu}2}\|r^{\frac{1-\mu}2} D^{-\theta} ( \nabla ( g(u_{k-1})) \cdot
  \nabla u_k)\|_{L^2_{t,x}}
+
T^{\frac{\mu}2}\|r^{\frac{1-\mu}2} D^{-\theta} F( u_{k-1})\|_{L^2_{t,x}}
\\
 &\les&
\| \pa u_k(0)\|_{\dot H^{-\theta}}+C(g, \|u_{k-1}\|_{L^\infty})
T^\mu\| u_{k-1}\|_{L^\infty_{t} \dot H^{\frac 1 2+\mu}}\| \pa u_k(0)\|_{\dot H^{1-\theta}}
\\
&&
+
  T^{{\mu}}\|D^{-\theta}\pa (u_{k}, u_{k-1})\|_{X_T}
\|(\nabla  (g(u_{k-1})),
(a(u_{k-1}), b(u_{k-1})) \pa u_{k-1})
\|_{L^\infty_{t} \dot H^{\frac 1 2+\mu}} 
\  ,  
\end{eqnarray*}
where, in the last inequality, we have used
Proposition \ref{thm-wtrace} and fractional chain rule based on the fact that $g(0)=0$.
 To control the last term, as 
$\nabla  (g(u_{k-1}))=g'(u_{k-1}) \nabla u_{k-1}$, we see that all terms are of the  form of $f(u) \pa u$, for which we could
 use the classical 
  fractional Leibniz and chain rule to conclude, 
\begin{eqnarray}\|f(u) \pa u\|_{\dot H^{\frac 1 2+\mu}}
 & \les & 
|f(0)|\|\pa u\|_{ \dot H^{\frac 12+\mu}}
+\| \tilde f(u)\|_{L^\infty}
\|\pa u\|_{ \dot H^{\frac 12+\mu}}
+\| \tilde f(u)\|_{\dot W^{\frac 1 2+\mu, 6}}
\|\pa u\|_{ L^{3}}\nonumber
\\
& \les & 
C(f, \|u\|_{\dot B^{\frac 12}_{2,1}})
\|\pa u\|_{ \dot H^{\frac 12+\mu}} \label{eq-Leib-cont}\ ,
\end{eqnarray}  
where $\tilde f(u)=  f(u)-f(0)$.

In view of the boundedness \eqref{eq-4.4'} and \eqref{eq-4.4Besovbd}, for $
\theta\in [1/2-\mu, 1/2+\mu]$, we see that
$$
 \|D^{-\theta} \pa u_k\|_{X_T}\le \frac C 2
\ep_{1-\theta}+C(g, a,b, \ep_c)
T^\mu\ep_s(
\ep_{2-\theta}+
 \|D^{-\theta}\pa (u_k, u_{k-1})\|_{X_T})
 \ .$$
Thus, for $T$ satisfying \eqref{eq-life-poly2} with sufficiently small $c$, 
 we get by the inductive assumption that
$$
 \|D^{-\theta} \pa u_k\|_{X_T}\le  C
\ep_{1-\theta}+C(g, a,b, \ep_c)
T^\mu\ep_{s-1}
\ep_{2-\theta}
 \ ,$$
which completes the proof.\end{prf}

\subsection{Convergence in $C\Hsd{s_0}$}
In this subsection, we show that the approximate solutions are convergent, in the weaker topology $C\dot H^{s_0}$, so that the desired solution of the quasilinear problem is given by the limit.

\begin{lem}
\label{thm-conv0}
Under the same assumption as in Lemma \ref{thm-bdd0}.
Let $(u_0, u_1) \in \Hs{s} \times (\Hs{s-1}\cap \Hsd{s_0-1})$ with $s_0\in [2-s, s-1]$, and
 $\{u_{k}\}_{k\ge 2}$ be the approximate solutions defined in \eqref{eq-4.3}, or equivalently, \eqref{eq-4.3'}, which satisfy
 the bounds 
 \eqref{eq-4.4'}, \eqref{eq-4.4Besovbd} and
 \eqref{eq-4.4'''}, for any $k\ge 2$. 
Then there exists $c=c(g, a, b,\ep_c)\in (0, 1)$,  such that
for any $T$ with
\beeq\label{eq-life-poly3}T\le \min(\de_1, c (\ep_s+\ep_{s-1})^{-1/\mu}),\eneq
we have
$u_k$ is Cauchy in the space $C([0,T]; \dot H^{s_0})\cap C^{0,1}([0,T]; \dot H^{s_0-1})$, with
\beeq\label{eq-seq-conv}
\sum_{k\ge 3} \|D^{s_0-1}\pa ( u_{k+1} - u_{k})\|_{X_{T}}
\les\ep_{s_0}+\ep_{s_0+1}
+\sum_{k\ge 3} \| \pa 
( u_{k+1} - u_{k})
(0)\|_{\dot H^{s_0-1}\cap\dot H^{s_0}} \ .
\eneq
Here the right hand side is bounded because of
\eqref{eq-4.2}-
\eqref{eq-4.2'}.
\end{lem}
\begin{prf}
If we set $w_{k} =  u_{k+1} - u_{k}$, it satisfies
\begin{align*}
 \Box w_k +& \nabla \cdot (g(u_{k})
  \nabla w_k)    =
   \nabla \cdot ((g(u_{k-1})-g(u_{k}))
  \nabla u_k) \\
    & +  
  \nabla ( g(u_{k}) )\cdot
  \nabla u_{k+1}-
    \nabla ( g(u_{k-1}) )\cdot
  \nabla u_k
 + F( u_{k})- F( u_{k-1}),
\end{align*}
for which we denote the right hand side by $G$.

As we see from the proof of Lemma \ref{thm-bdd0}, we know that
\eqref{eq-3.4}  is satisfied for $h=g(u_{k})$ and we could apply Proposition \ref{th-LEkey3} 
with $\theta\ge 1/2-\mu=2-s$
to obtain 
\beeq\label{eq-seq-conv0}
\|D^{-\theta}\pa w_{k}\|_{X_{T}}\les
\| \pa w_{k}(0)\|_{\dot H^{-\theta}}+
T^\mu\| g(u_{k})\|_{L^\infty_{t} \dot H^{s-1}}\| \pa w_{k}(0)\|_{\dot H^{1-\theta}}
+ \|D^{-\theta} G\|_{X_T^*}\ .
\eneq
For the term involving $g(u_{k})$, we know from
Theorem \ref{thm-wLeib0}, $g(0)=0$,  \eqref{eq-4.4Besovbd} and \eqref{eq-4.4''}
that
\beeq\label{eq-bd-guk}\| g(u_{k})\|_{L^\infty_{t} \dot H^{s-1}}\le C(g, \ep_{c})
\|  u_{k}\|_{L^\infty_{t} \dot H^{s-1}}
\les
\|(\nabla u_0, u_1)\|_{\dot H^{s-2}}=\ep_{s-1}\ .
\eneq

The main part of the proof is to deal with $G$. 
We will write it into a combination of  favorable terms and
deal with each term separately. For this purpose,
we set
$G_{1}=\nabla \cdot ((g(u_{k-1})-g(u_{k}))
  \nabla u_k)$,
  $G_{2}=F( u_{k})- F( u_{k-1})$ and then
\begin{align*}
G-G_{1}&-G_{2}=     g'(u_{k})  \nabla u_{k}\cdot
  \nabla u_{k+1}-
g'(u_{k-1}) \nabla u_{k-1}\cdot
  \nabla u_k\\
    =&  g'(u_{k}) 
    \nabla u_{k} \cdot  \nabla w_{k}+ 
     g'(u_{k}) 
         \nabla w_{k-1} \cdot  \nabla u_{k}+
 (g'(u_{k})- g'(u_{k-1}) )\nabla u_{k-1}\cdot
  \nabla u_k\\
  &=G_{3}+G_{4}+G_{5}\ .
\end{align*}
For $G_j$ with $j\ge 2$, we observe that they fall into the following two categories:
$$\tilde G_2=(f(u_{k})- f(u_{k-1}) )\pa(u_{k-1}, u_k)
 \pa (u_{k-1}, u_k)\ ,$$
 $$\tilde G_3= f(u_{k}) \pa (u_{k-1}, u_k)
 \pa (w_{k-1}, w_k)\ .$$
For all these terms, we claim that we have for any $\theta\in [2-s, s-1]$,
\beeq \label{eq-G}\|D^{-\theta} G\|_{X_T^*} \le
C(g,a,b,\ep_c)T^\mu \ep_s \|D^{-\theta}\pa (w_{k-1}, w_{k})\|_{X_T} \ .
\eneq

Before presenting the proof of  \eqref{eq-G}, we apply it 
to prove \eqref{eq-seq-conv}.
Actually, by \eqref{eq-seq-conv0} and \eqref{eq-bd-guk}, we have
$$\|D^{-\theta}\pa w_{k}\|_{X_{T}}\les
\| \pa w_{k}(0)\|_{\dot H^{-\theta}}+
T^\mu \ep_{s-1} \| \pa w_{k}(0)\|_{\dot H^{1-\theta}}
+ T^\mu \ep_s \|D^{-\theta}\pa (w_{k-1}, w_{k})\|_{X_T}\ ,
$$
for any $\theta\in [2-s, s-1]$,
 where the implicit constant may depend on $g,a,b,\ep_c$.
Let $\theta=1-s_0$,
then for any $T$ satisfying
\eqref{eq-life-poly3} with sufficiently small $c$, we have
$$\|D^{s_0-1}\pa w_{k}\|_{X_{T}}\le
C(\| \pa w_{k}(0)\|_{\dot H^{s_0-1}}+
 \| \pa w_{k}(0)\|_{\dot H^{s_0}})
+ \frac 14 \|D^{s_0-1}\pa (w_{k-1}, w_{k})\|_{X_T}
\ ,
$$
and so
$$\|D^{s_0-1}\pa w_{k}\|_{X_{T}}\le
2C(\| \pa w_{k}(0)\|_{\dot H^{s_0-1}}+
 \| \pa w_{k}(0)\|_{\dot H^{s_0}})
+ \frac 12 \|D^{s_0-1}\pa  w_{k-1}\|_{X_T}
\ ,
$$
for any $k\ge 3$.
However,
recall that $w_{2}=u_3-u_2=u_3$, we know from
 \eqref{eq-4.4'''} and \eqref{eq-life-poly3} that
$$\|D^{s_0-1}\pa  w_2\|_{X_T}\le 
C( \ep_{s_0}
+ T^\mu
\ep_{s_0+1}
\ep_{s-1})\le
C  (\ep_{s_0}+
\ep_{s_0+1})
\ .$$
Thus an iteration argument gives us
that
$\sum\|D^{s_0-1}\pa w_{k}\|_{X_{T}}$ is convergent and  we have \eqref{eq-seq-conv}.
Actually, for any $j\in [3,\infty)$, for finite summation from $3$ to $j$, we have
$$\sum_{k=3}^j \|D^{s_0-1}\pa w_{k}\|_{X_{T}}\le
2C\sum_{k=3}^j \| \pa w_{k}(0)\|_{\dot H^{s_0-1} \cap\dot H^{s_0}}
+\sum_{k=2}^{j-1} \frac 12 \|D^{s_0-1}\pa  w_{k}\|_{X_T}
\ ,
$$
and so
\begin{eqnarray*}
\sum_{k=3}^j \|D^{s_0-1}\pa w_{k}\|_{X_{T}} & \le & 4C\sum_{k=3}^j \| \pa w_{k}(0)\|_{\dot H^{s_0-1} \cap\dot H^{s_0}}
+\|D^{s_0-1}\pa  w_{2}\|_{X_T} \\
 & \le & 4C\sum_{k=3}^j \| \pa w_{k}(0)\|_{\dot H^{s_0-1} \cap\dot H^{s_0}}
+C  (\ep_{s_0}+
\ep_{s_0+1})
\ .
\end{eqnarray*}
Letting $j$ goes to $\infty$, we obtain
 \eqref{eq-seq-conv}.

 It remains to prove  \eqref{eq-G},  for which we divide it into three terms, $G_1$, $\tilde G_2$ and $\tilde G_3$.

{\noindent \bf i) first term: $G_1=\nabla \cdot ((g(u_{k-1})-g(u_{k}))
  \nabla u_k)$}.
For the first term $G_1$, we see that
\begin{eqnarray*}
 \|D^{-\theta} G_1\|_{X_T^*} & \le & 
 T^{\frac{\mu}2}\|r^{\frac{1-\mu}2} D^{-\theta} G_1\|_{L^2_{t,x}}\\
 &\les&
 T^{\frac{\mu}2}\|r^{\frac{1-\mu}2} D^{1-\theta}((g(u_{k-1})-g(u_{k}))
  \nabla u_k) \|_{L^2_{t,x}}\\
 &\les&
 T^{\frac{\mu}2}\|r^{-\frac{1-\mu}2} D^{1-\theta}(g(u_{k-1})-g(u_{k}))
 \|_{L^2_{t,x}}
\|   \nabla u_k\|_{L^\infty_{t} \dot H^{s-1}}
\end{eqnarray*}
where, as $\theta\in [2-s, s]$, in the last inequality, we have used
 Proposition \ref{thm-wtrace} with $|1-\theta|\le s-1$. To control the term involving
$g(u_{k-1})-g(u_{k})$, we observe that
$$g(u)-g(v)
=\int_0^1 g'(v+\la(u-v)) (u-v) d\la
$$
and so
\begin{eqnarray*}
&&\|r^{-\frac{1-\mu}2} D^{1-\theta}(g'(v+\la(u-v)) (u-v)) \|_{L^2_{t,x}} \\
& \les &  \|r^{-\frac{1-\mu}2} D^{1-\theta} (u-v) \|_{L^2_{t,x}} 
\|g'(v+\la(u-v))\|_{L^\infty_{t,x}}
\\&&+
\|r^{-\frac{1-\mu}2}  (u-v) \|_{L^2_{t}
L^{6/(1+2\theta)}_{x}
} 
\|D^{1-\theta} g'(v+\la(u-v))\|_{L^\infty_t L^{3/(1-\theta)}_x}
\\
 & \les & 
 \| g'(v+\la(u-v))\|_{L^\infty_t \dot W^{1-\theta,3/(1-\theta)}_x
\cap L^\infty_{t,x} 
 }
 \|r^{-\frac{1-\mu}2} D^{1-\theta} (u-v) \|_{L^2_{t,x}} 
 \\
 & \les & C(g, \|(u,v)\|_{L^\infty_{t}
\dot B^{3/2}_{2,1}}) T^{\mu/2} \|D^{1-\theta}(u-v)\|_{X_T}\ ,
\end{eqnarray*} where we have used
Theorem \ref{thm-wLeib-fromChain} and
 \eqref{eq-trace-w2} in the first and second inequalities, for $\theta\in [0, 1]$.
In summary, we have proved that
\beeq
\|r^{-\frac{1-\mu}2} D^{1-\theta}(g(u)-g(v))
 \|_{L^2_{t,x}}\le C(g, \|(u,v)\|_{L^\infty_{t}\dot B^{\frac 3 2}_{2,1}}) T^{\frac \mu 2} \|D^{1-\theta}(u-v)\|_{X_T}\ ,\eneq
which gives us
\beeq \label{eq-G1}\|D^{-\theta} G_1\|_{X_T^*} \le
 C(g, \ep_c) T^{\mu} \|D^{1-\theta} w_{k-1}\|_{X_T}
\|   \nabla u_k\|_{L^\infty_{t} \dot H^{s-1}}\ .
\eneq
{\noindent \bf ii) second category of terms: $\tilde G_2$}.
Recall that $$\tilde G_2=(f(u_{k})- f(u_{k-1}) )\pa(u_{k-1}, u_k)
 \pa (u_{k-1}, u_k)\ .$$
  Let us present the proof for the
typical term $\tilde G_2=(f(u_k)- f(u_{k-1}) )\pa y
 \pa z$, for which we know that, as $\theta\in [0, s-1]$,
$\|D^{-\theta} \tilde G_2\|_{X_T^*}$ is bounded by
 \begin{eqnarray*}
 &  & 
 T^{\frac{\mu}2}
 \|r^{-\frac{1-\mu}2} D^{-\theta} ( (f(u_{k})- f(u_{k-1}) )\pa y) \|_{L^2_{t,x}}
 \|\pa z\|_{L^\infty \dot H^{s-1}}
 \\
 & \les & 
 T^{\frac{\mu}2}
 \|r^{-\frac{1-\mu}2+\theta} ( (f(u_{k})- f(u_{k-1}) )\pa y) \|_{L^2_{t,x}}
 \|\pa z\|_{L^\infty \dot H^{s-1}}
 \\
 &\les&
 T^{\frac{\mu}2}\|r^{-\frac{1-\mu}2-1+\theta} 
 (f(u_{k})- f(u_{k-1}) ) \|_{L^2_{t,x}}
\|r \pa y\|_{L^\infty_{t,x}}
 \|\pa z\|_{L^\infty \dot H^{s-1}}\\
 &\les&
    C(f, \|(u_{k-1}, u_{k})\|_{L^\infty_{t,x}})
 T^{\frac{\mu}2}\|r^{-\frac{1-\mu}2-1+\theta}   w_{k-1} \|_{L^2_{t,x}}
\|\pa y\|_{L^\infty_{t}
\dot B^{1/2}_{2,1}}
 \|\pa z\|_{L^\infty \dot H^{s-1}}
 \\
 &\les&
    C(f, \|(u_{k-1}, u_{k})\|_{L^\infty_{t,x}})
 T^{\frac{\mu}2}\|r^{-\frac{1-\mu}2}D^{1-\theta}   w_{k-1} \|_{L^2_{t,x}}
\|\pa y\|_{L^\infty_{t}
\dot B^{1/2}_{2,1}}
 \|\pa z\|_{L^\infty \dot H^{s-1}}\ .
\end{eqnarray*}
That is, we have
\beeq\label{eq-G2} \|D^{-\theta} \tilde G_2\|_{X_T^*} \le C(f, \ep_c) T^{\mu} \|D^{1-\theta} w_{k-1}\|_{X_T}\|   \pa (u_k, u_{k-1})\|_{L^\infty_{t} \dot H^{s-1}}\ .\eneq

 {\noindent \bf iii)  third category of terms: $\tilde G_3= f(u_{k}) \pa (u_{k-1}, u_k)
 \pa (w_{k-1}, w_k)$}. In this case, 
with the help of Proposition \ref{thm-wtrace},
 we see that
 $\|D^{-\theta} \tilde G_3\|_{X_T^*}$ is bounded by
$$
 T^{\frac{\mu}2}
 \|r^{-\frac{1-\mu}2} D^{-\theta}  \pa (w_{k-1}, w_k)
  \|_{L^2_{t,x}}
 \| f(u_{k}) \pa (u_{k-1}, u_k)\|_{L^\infty \dot H^{s-1}}\ .
$$
Similar to the proof of \eqref{eq-Leib-cont} in Lemma \ref{thm-bdd1}, we know that
$$ \| f(u_{k}) \pa (u_{k-1}, u_k)\|_{L^\infty_t \dot H^{s-1}}
\les C(f, \|\pa (u_{k-1}, u_k)\|_{L^\infty_t\dot B^{1/2}_{2,1}})
\| \pa (u_{k-1}, u_k)\|_{L^\infty_t \dot H^{s-1}}\ ,
$$
and so
\beeq\label{eq-G3} \|D^{-\theta} \tilde G_3\|_{X_T^*} \le
 C(f, \ep_c
 ) T^{\mu} \|D^{-\theta}\pa (w_{k-1}, w_{k})\|_{X_T}
\|   \pa (u_k, u_{k-1})\|_{L^\infty_{t} \dot H^{s-1}}\ .
\eneq

In summary, in view of  \eqref{eq-G1},  \eqref{eq-G2} and \eqref{eq-G3},
as well as the  uniform bounds \eqref{eq-4.4'} and \eqref{eq-4.4Besovbd}, 
we complete the proof of \eqref{eq-G} and 
Lemma \ref{thm-conv0}.
\end{prf}

\subsection{Local 
wellposedness in $H^s$}
Equipped with 
Lemma \ref{thm-bdd0},
Lemma \ref{thm-bdd1},
and  Lemma \ref{thm-conv0}, we are ready to prove the (unconditionally) local wellposedness.
\begin{lem}
\label{thm-lwp0}
Let $n=3$, $s=3/2+\mu\in (3/2, 2]$ and 
$s_0\in [2-s, s-1]$.
Considering the initial value problem \eqref{ea1}-\eqref{ea2}, 
with $(u_0, u_1) \in \Hs{s}  \times (\Hs{s-1}\cap \Hsd{s_0-1})$.
Then, for any $T$ satisfying
\eqref{eq-life-poly3},
there exists a unique weak solution 
\beeq\label{eq-19-0}u \in L^{\infty}_t \Hs{s}  \cap C^{0,1}_t\Hs{s-1}\cap C_t \dot H^{s_0} \cap C^{0,1}_t\dot H^{s_0-1}  \eneq
in $[0,T]\times\R^3$ for
the initial value problem \eqref{ea1}-\eqref{ea2}. 
Moreover, there exists $C_2>0$ such that the solution satisfies
 $\pa u \in C([0,T]; \dot H^{\theta-1})$ for any $\theta\in [s_0,s)$, 
\begin{equation}\label{eq-19-1}\|\pa u\|_{L^\infty \dot H^{\theta}}
\le
\|\pa D^\theta u\|_{X_{T}} 
\le C_2 \|(\nabla u_0, u_1)\|_{\dot H^{\theta}}\ , \forall \theta\in [0, s-1] 
\ ,
\end{equation}
\begin{equation} \label{eq-19-2}
\| \pa u\|_{L^\infty_t \dot B^{1/2}_{2,1}}
\le 
\|\pa D^{1/2} u\|_{X_{T,1}}
\le C_2 \|(\nabla u_0, u_1)\|_{\dot B^{1/2}_{2,1}}\ ,
\end{equation}
 \beeq  \label{eq-19-3}
 \|\pa u\|_{L^\infty \dot H^{\theta-1}}
\le
\|D^{\theta-1}\pa u\|_{X_T}\le C_2(\ep_\theta+\ep_{s-1})
\ ,
\forall \theta\in [s_0,s-1]
\ .
 \eneq
 \end{lem}

\begin{prf}
By Lemma \ref{thm-bdd0}
and Lemma \ref{thm-bdd1},
the approximate solutions ${u_{k}}$ are well defined and satisfy the bounds
 \eqref{eq-4.4'}, \eqref{eq-4.4Besovbd} and
 \eqref{eq-4.4'''}. Moreover,
  Lemma \ref{thm-conv0} tells us that  ${u_{k}}$
 is Cauchy in the space $C([0,T]; \dot H^{s_0})\cap C^{0,1}([0,T]; \dot H^{s_0-1})$, for which we denote the limit by $u
\in
C([0,T]; \dot H^{s_0})\cap C^{0,1}([0,T]; \dot H^{s_0-1})$.
 By 
 Helly's selection theorem,  we see that there is a subsequence of ${u_{k}}$, which is weak star convergent to $u$,
 in  $L^{\infty}([0,T]; \Hsd{s}) \cap C^{0,1}([0,T]; \Hsd{s-1})$, and so 
we have \eqref{eq-19-0}.
Then, it is clear that $\pa u_k$ is convergent to $\pa u$ in $C([0,T]; \dot H^{\theta-1})$ for any $\theta\in [s_0,s)$,  which follows directly, by interpolation, from the boundedness for $\theta=s$ and continuity for $\theta=s_0$.
Consequently, in view of the definition of $u_k$, 
\eqref{eq-4.3},
$u$ is the desired weak solution for 
 the initial value problem \eqref{ea1}-\eqref{ea2}, as well as the bounds
 \eqref{eq-19-1}- 
\eqref{eq-19-3}.
 
 It remains to prove the unconditional uniqueness.
Suppose there is a solution
\beeq\label{eq-19-4}v \in L^{\infty}_t \Hs{s}  \cap C^{0,1}_t\Hs{s-1}\cap C_t \dot H^{s_0},
\pt v \in C_t\dot H^{s_0-1} \ , \eneq
in $[0, T_1]\times \R^3$
for
the initial value problem \eqref{ea1}-\eqref{ea2},
 for some $T_1\in (0, T]$. 
The key observation here is that
by \eqref{eq-19-4}, we have
$$D^{s_0-1}\pa v\in L^\infty H^1$$
and so by Hardy's inequality,
$$\|r^{-\frac{1-\mu}{2}}D^{s_0-1} \pa v\|_{L^2([0,{T_2}]\times \R^3)}
\les_T \|D^{s_0-1} \pa v\|_{L^\infty([0,{T_2}]; \dot H^{\frac{1-\mu}2})}
\les \|D^{s_0-1} \pa v\|_{L^\infty([0,{T_2}]; H^1)}
$$
for any ${T_2}\in (0, T_1]$. In other words, we see that
$D^{s_0-1}\pa v\in X_{T_2}$
 for any ${T_2}\in (0, T_1]$.

Similar to the proof of Lemma \ref{thm-conv0}, we 
set $w  =  u - v$ with $D^{s_0-1}\pa w\in X_{T_2}$, and write the equation for $w$ as follows
\begin{align*}
-\partial_t^2 w + \Delta w +& \nabla \cdot (g(u)
  \nabla w)    =
   \nabla \cdot ((g(v)-g(u))
  \nabla v) \\
    & +  
  \nabla ( g(u) )\cdot
  \nabla u-
    \nabla ( g(v) )\cdot
  \nabla v
 + F( u)- F( v):=G(u, w)\ ,
\end{align*}
together with $w(0,x)=0$, $\pt w(0,x)=0$.

As $u$ is constructed as limit of $u_k$, 
 we could apply Proposition \ref{th-LEkey3} 
with $\theta=1-s_0$ for 
the wave operator $-\partial_t^2 + \Delta + \nabla \cdot g(u)
  \nabla$. That is, we have
\beeq\label{eq-Wave-perturb1}
\|D^{s_0-1}\pa w\|_{X_{T}}\les
\|D^{s_0-1} G\|_{X_T^*})\ .
\eneq
With the help of \eqref{eq-Wave-perturb1},
applied to  $w  =  u - v$,
 together with
 the similar proof as
\eqref{eq-G1},  \eqref{eq-G2} and \eqref{eq-G3}, we get that
\begin{eqnarray*}
\|D^{s_0-1}\pa w\|_{X_{T}} &  \les &
 \|D^{s_0-1} G(u,w)\|_{X_T^*}
 \\
& \les & 
 C(g,a,b, \|\pa (u, v)\|_{L^\infty_{
T_1 
}\dot B^{1/2}_{2,1}}) T^{\mu} \|D^{s_0-1}\pa w\|_{X_T}
\|   \pa (u, v)\|_{L^\infty_{T_1
} \dot H^{s-1}}
 \\
& \les & 
T^{\mu} \|D^{s_0-1}\pa w\|_{X_T}\ .
\end{eqnarray*}
Thus, with $T_2\in (0, T_1]$ sufficiently small, we see that
$\|D^{s_0-1}\pa w\|_{X_{T_2}}=0$ and so
$w\equiv 0$ in $[0,{T_2}]\times\R^3$,
in view of $w(0,x)=0$. After a simple iteration argument, this proves that
$u\equiv v$ in $[0,T_1]\times\R^3$, which completes the proof of unconditional uniqueness.\end{prf}

\section{High dimensional well-posedness}\label{sec-5-high}

Let $n\ge 4$, 
 $s=n/2+\mu$ with 
$\mu$ as in \eqref{eq-mu}, and
$\ep_s$, $\ep_c$ be as in
\eqref{eq-1-7}.
In this section, we prove the existence and uniqueness part of
Theorem \ref{th-high}, following the similar approach  as in Section \ref{sec-4}.

\subsection{Approximate solutions}
As in subsection \ref{sec-approx}, we could construct a sequence of 
spherically symmetric, compactly supported, smooth functions
$(u_0^{(k)}, u_1^{(k)})\to (u_0, u_1)$ in $\Hs{s}\times\Hs{s-1}$, such that
\begin{equation}
\|(\nabla u_0^{(k)},u_1^{(k)})\|_{\dot B^{\theta }_{2,q}}\le C_{\theta, q}
\|(\nabla u_0,u_1)\|_{\dot B^{\theta }_{2,q}}, \ \forall \theta \in \R,\ q\in [1,\infty]\ ,
\label{eq-6.1}
\end{equation}
\begin{equation}
\|\nabla  u_0^{(k)}
 - \nabla  u_0^{(k+1)}\|_{\Hs{s}({\mathbb R}^n)}
+ \| u_1^{(k)}
 -  u_1^{(k+1)}\|_{\Hs{s-1}({\mathbb R}^n)}
\le 2^{-k}.
\label{eq-6.2}
\end{equation}
 Let 
$F( u)=  a(u) u_t^2+b(u) |\nabla u|^2$,
$u_{2} \equiv 0$ and define $u_k$ ($k\ge 3$) recursively by solving
\beeq
\left\{\begin{array}{l}\Box u_k  + g(u_{k-1})
  \Delta u_k =  F( u_{k-1}), (t,x)\in (0,T)\times \R^n,
 \\
u_k(0,\cdot) =  u_0^{(k)},\quad
 \partial_t u_k(0,\cdot) =  u_1^{(k)}.\end{array}
 \right.
 \label{eq-6.3}
\eneq

\subsection{Uniform boundedness of $u_k$}
Let $C$ be the implicit constant in the estimates of Proposition \ref{th-LEkey4}.
We claim that we have
 the uniform bounds
\beeq\label{6eq-LE-high-k}
\|\tilde \pa D^\theta u_k\|_{X_T}\le 2 C
\|D^\theta(\nabla u_0, u_1)\|_{L^2} , \ \theta\in [0, s-1]\ ,
\eneq
\beeq\label{6eq-LE-high-k-critical}
\|\pa D^{\frac n 2-1}u_k\|_{X_{T,1}} 
\le 2 C
\|(\nabla u_0, u_1)\|_{\dot B^{{\frac n 2-1}}_{2,1}}     \ ,
\eneq
for any $T>0$ satisfying
\beeq\label{6eq-life}
T^{\mu} 
f(C \ep_c)
\ep_s\le c\ ,\eneq
for some increasing function $f$, and constants $c\ll 1\ll C$.

We prove the bounds by induction. It is trivially true when $k=2$. Assuming for some $m\ge 2$, it is true for  any $k\le m$, then,
for $h=g(u_{m})$, we have
\beeq\label{6eq-LE-frac-h}T^\mu\|\pa h\|_{L^\infty_t \dot H^{\frac n2-1+\mu}} \le C(\ep_c)
T^{\mu} 
\ep_s 
\le \de\ ,\eneq
and so
is the requirement 
\eqref{eq-LE-frac-h}
of Proposition \ref{th-LEkey4} satisfied.

Based on Proposition \ref{th-LEkey4}, we know that
$$
\|\tilde \pa D^\theta u_{m+1}\|_{X_T}\le 
C\|D^\theta(\nabla u_0, u_1)\|_{L^2} +C
\|D^\theta   F( u_{m})\|_{X_T^*}\  ,\ \theta\in[0, [\frac n2]], 
$$$$
\|\pa  D^{\frac n 2-1}u_{m+1}\|_{X_{T,1}} 
\le C
\|(\nabla u_0, u_1)\|_{\dot B^{{\frac n 2-1}}_{2,1}}     
+C T^{\frac{\mu}2}\|r^{\frac{1-\mu}2}2^{j({\frac n 2-1})}P_j  F( u_{m})\|_{\ell_j^1 L^2_{t,x}}
\ .$$
To control the nonlinear term, we apply Proposition \ref{thm-wtrace}
to obtain, for the sample term $F(u)=a(u)u_t^2$,
\begin{eqnarray*}
\|D^{\theta} F( u)\|_{X_T^*}
 & \les & T^{\frac\mu 2}\|r^{-\frac{1-\mu}2}D^{\theta} \pa u\|_{ L^2_{t,x}} 
 \|a(u)\pa u\|_{L^\infty_t \dot H^{s-1}}
 \\
& \les &T^{\mu} \|\tilde \pa D^\theta u\|_{X_T}
\ \tilde{a}(\|u\|_{L^\infty \dot B^{\frac n2}_{2,1}})
\|\pa u\|_{L^\infty_t \dot H^{s-1}}
\end{eqnarray*}
whenever $\theta\in [0, s-1]$, where, since $s-1<n/2$,
we have used the following well-known consequence of the fractional Leibniz rule
and chain rule
$$
\|(a(u)-a(0))v\|_{\dot H^{s-1}}
\les
\|a(u)-a(0)\|_{\dot B^{\frac n2}_{2,1}}\|
v\|_{ \dot H^{s-1}}
\les C(\|u\|_{L^\infty})
\|u\|_{\dot B^{\frac n2}_{2,1}}\|v\|_{ \dot H^{s-1}}\ .
$$
Similarly, by \eqref{eq-wtrace2}, we have
\begin{eqnarray*}
\|r^{\frac{1-\mu}2}2^{j({\frac n 2-1})}P_j  F( u)\|_{\ell_j^1 L^2_{t,x}}
 & \les & 
 \|r^{-\frac{1-\mu}2}2^{j({\frac n 2-1})}P_j  \pa u\|_{\ell_j^1 L^2_{t,x}}
 \|a(u)\pa u\|_{L^\infty_t \dot H^{s-1}}  \\
 & \les & T^{\frac \mu 2} \|\pa D^{\frac n2-1}u\|_{X_{T,1}}  
\tilde{a}(\|u\|_{L^\infty \dot B^{\frac n2}_{2,1}})
\|\pa u\|_{L^\infty_t \dot H^{s-1}}\ .
\end{eqnarray*}
Based on the induction assumption and
\eqref{6eq-life}, we 
get
\eqref{6eq-LE-high-k} and
\eqref{6eq-LE-high-k-critical} for $k=m+1$,
  if we set $c>0$ to be sufficiently small.
This completes the proof by induction.

\subsection{Convergence in  $C\dot H^1\cap C^{0,1} L^2$}
Let $w_{k} =  u_{k+1} - u_{k}$, it satisfies
$$
 \Box w_k + g(u_{k})
  \Delta w_k   =
(g(u_{k-1})-g(u_{k}))
  \Delta u_k 
 + F( u_{k})- F( u_{k-1}).
$$
Thus, by Theorem \ref{th-LE0}, we have
$$\|\pa w_k\|_{X_T}
\les 
\|\pa w_k(0)\|_{L^2}+
\|
(g(u_{k-1})-g(u_{k}))
  \Delta u_k \|_{X_T^*}
 +\| F( u_{k})- F( u_{k-1})\|_{X_T^*}.
$$
Notice that
$\|
(g(u_{k-1})-g(u_{k}))
  \Delta u_k \|_{X_T^*}$ is controlled by
\begin{eqnarray*}
 &  &
    T^{\frac \mu 2}
C(\ep_c)\|r^{\frac n 2-1}w_{k-1} \|_{L^\infty}
\|r^{-\frac{ 1-\mu}2-(\frac n 2-2+\mu)}  \Delta u_k\|_{L^2}
 \\
 & \les & T^{\frac\mu 2}
C(\ep_c)\|\pa w_{k-1} \|_{L^\infty L^2}
\|r^{-\frac{1-\mu}2} D^{s}  u_k\|_{L^2}
 \\
 & \les & T^{\mu}\ep_s
C(\ep_c)\|\pa w_{k-1} \|_{X_{T}}\ .
\end{eqnarray*}
Similarly, for the sample term $F(u)=a(u)u_t^2$, we have
\begin{eqnarray*}
&&
\| F( u_{k})- F( u_{k-1})\|_{X_T^*}\\
&\les &T^{\frac\mu 2}
\|r^{\frac n2-1}(a(u_{k})-a(u_{k-1}))\|_{L^\infty} \|r \pa u_{k}\|_{L^\infty}
\|r^{\frac{1-\mu}2-\frac n 2} \pa u_{k}\|_{L^2}\\
&&+C_1(\ep_c) T^{\frac\mu 2} \|r^{1-\mu}\pa (u_{k-1}, u_{k})\|_{L^\infty}
\|r^{-\frac{1-\mu}2}\pa w_{k-1}\|_{L^2}\\
&\les&
T^{\mu}\ep_s\tilde
C(\ep_c)\|\pa w_{k-1} \|_{X_{T}}\ .
\end{eqnarray*}
By letting
$c$ 
in \eqref{6eq-life}
be even smaller, we could conclude
 $$\|\pa w_k\|_{X_{T}}
\le C 
\|\pa w_k(0)\|_{L^2}+
\frac12 \|\pa w_{k-1} \|_{X_{T}}\ ,$$
which yields convergence in $C\dot H^1\cap C^{0,1} L^2$, thanks to \eqref{eq-6.2}.

\subsection{Local well-posedness}
Let $u\in C  H^1\cap C^{0,1} L^2$ be the limit of $u_k$, 
by weak star compactness,
we have $\pa u\in L^\infty H^{s-1}$ and
$\pa u_k$ is convergent to $\pa u$ in $C([0,T]; \dot H^{\theta-1})$ for any $\theta\in [1,s)$.
Then, in view of the definition of $u_k$, 
\eqref{eq-6.3},
it is clear that $u$ is a weak solution for 
 the initial value problem \eqref{ea1}-\eqref{ea2}.
 
The unconditional uniqueness follows from the similar argument as that in Lemma \ref{thm-lwp0}, and we omit the proof.

\section{Persistence of regularity}\label{sec-6persi} 
In this section, we show that 
persistence of regularity for the weak solutions, when the initial data have higher regularity,
as well as the continuous dependence on the data.

In Sections \ref{sec-4} and \ref{sec-5-high}, 
for data in $H^s\times H^{s-1}$, with additional requirement in $\dot H^{s_0-1}$ for initial velocity when $n=3$,
we have constructed solutions in $H^s$, when
$$
s\in
\left\{
\begin{array}{ll }
  (n/2, (n+1)/2]  ,  &   n \textrm{ odd} \ ,\\
  (n/2, (n+2)/2) ,      &    n \textrm{ even}\ .
\end{array}
\right.
$$
Recall the classical energy argument shows local well-posedness in $H^{s_1}$ for any $s_1>n/2+1$, together with 
the persistence of higher regularity.
Keeping this fact in mind, we need only to prove
 the persistence of regularity in $H^{s_1}$, with
$s_1=[(n+4)/2]$.

\subsection{Persistence of regularity: a weaker version}\label{sec:persist-slight}
At first, we prove a weaker version of persistence of regularity, that is, when the data has slightly better regularity
$s_2=n/2+\mu_2$, $(u_0, u_1)\in \dot H^{s_2}\times \dot H^{s_2-1}$, with $$\mu_2\in 
\left\{
\begin{array}{ll }
  (\mu, 1/2]  ,  &   n \textrm{ odd} \ ,\\
  (\mu, (\mu+1)/2) ,      &    n \textrm{ even}\ ,
\end{array}
\right.$$
where
$\mu=s-n/2$.

Fix $T<T_*$,
we have uniform bound on $\|\pa u\|_{L^\infty_t \dot H^{s-1}\cap
L^\infty_t \dot B^{n/2-1}_{2,1} 
([0,T]\times\R^n)}$.
With $\de_3>0$ to be determined,
by dividing $[0,T]$ into finitely many small, disjoint, adjacent intervals $I_j$, we have 
$|I_j|^{\mu}\|\pa u\|_{L^\infty_t \dot H^{s-1}(I_j\times\R^n)}
\le \de_3$,
so that
\beeq\label{5eq-LE-frac-h}
|I_j|^{\mu}\|\pa g(u)\|_{L^\infty_t \dot H^{s-1}(I_j\times\R^n)} \le \de\eneq
for each $I_j=[T_j, T_{j+1}]$,
with $\Delta T_j=|I_j|$,
where
 $\delta$ is that occurred in
\eqref{eq-LE-frac-h}
of Proposition \ref{th-LEkey4}.
 In addition, we could possibly shrink $I_j$, so that we could apply the iteration argument in $I_j$ to obtain uniform bound in $H^s$, for data $(u(T_j), \pt u(T_j))$ at $t=T_j$.

By recasting the  iteration argument  for local well-posedness on $I_j$, we obtain for the iterative $C^\infty$ sequence $u_k$ on $I_j$, with
$$\|\pa u_k\|_{L^\infty_t \dot H^{s-1}(I_j\times\R^n)}
\le C_j \|\pa u(T_j)\|_{\dot H^{s-1}},
\|\pa u_k \|_{L^\infty_t \dot B^{n/2-1}_{2,1} }
\le C_j \|\pa u(T_j)\|_{\dot B^{n/2-1}_{2,1} },
$$
$$\lim_{k\to \infty}\|\pa (u_k-u)\|_{L^\infty_t L^2(I_j\times\R^n)}=0\ .$$
Assuming, by induction in $j$, that
\beeq\label{5eq-induction}
\|\pa u_k (T_j)\|_{\dot H^{s_2-1}}\le
C \|\pa u(T_j)\|_{\dot H^{s_2-1}}
\le \tilde C_j 
\|\pa u(0)\|_{\dot H^{s_2-1}}
\ .
\eneq
Applying Proposition \ref{th-LEkey4} with $\theta=s_2-1$, we have
\beeq\label{5thm-linea}
\|D^{\theta} u_{k+1}\|_{\widetilde{LE}_{{I_j},{\mu}}}:=
\|\tilde \pa D^{\theta} u_{k+1}\|_{X_{\Delta T_j}(I_j)}\les
\| \pa u(T_j)\|_{\dot H^{\theta}} +
\|D^\theta F(u_{k})\|_{X_{\Delta T_j}^*(I_j)}\  .
\eneq

As for the nonlinear term, we have
\begin{lem}\label{5thm-nl}
Let $n$ be odd or $\mu_2<(\mu+1)/2$, $F(u)=a(u) u_t^2+b(u)|\nabla u|^2$, then for radial functions $u$,
\beeq\label{5eq-nl}\| D^{s_2-1} F(u)
\|_{X_{\Delta T_j}^*(I_j)}
\les 
  C(
\|\pa u \|_{L^\infty_t \dot B^{n/2-1}_{2,1} })
|I_j|^{ {\mu}} \|D^{s_2-1} u\|_{\widetilde{LE}_{{I_j},{\mu}}}
\| \pa u\|_{L^\infty_{t} \dot H^{s-1}}\ .\eneq
\end{lem}
\begin{prf}
As in \eqref{eq-Fuexp}, without loss of generality, we deal with $F_1(u)$ and $F_2(u)=u_t^2$.
For $F_2(u)=u_t^2$, we have
\begin{eqnarray*}
\|r^{\frac{1-{\mu}}2}D^{\theta} F_2(u)\|_{ L^2_{t,x}} & \les & 
\|r^{-\frac{1-{\mu}}2}D^{\theta}u_t\|_{  L^2_{t,x}}
\|r^{1-{\mu}} u_t\|_{L^\infty_{t,x}} \\
 & \les & |I_j|^{\frac {\mu} 2} \|D^\theta u\|_{\widetilde{LE}_{{I_j},{\mu}}}
\| \pa u\|_{L^\infty_{t} \dot H^{n/2-1+{\mu}}}
\end{eqnarray*}
by
Theorem \ref{thm-wLeib-fromChain} and \eqref{eq-2.1}.
Concerning the other term
$F_1(u)=\tilde a(u) u_t^2=\tilde a(u) F_2(u)$, with
$\tilde a(0)=0$, we get from 
Theorem \ref{thm-wLeib-fromChain}
that
\begin{eqnarray*}
&&
\|r^{\frac{1-{\mu}}2}D^{\theta} F_1(u)\|_{ L^2_{t,x}}
\\
& \les & 
\|\tilde a(u) \|_{L^\infty_{t,x}}
\|r^{\frac{1-{\mu}}2}D^{\theta} F_2(u)\|_{ L^2_{t,x}}
+
\|r^{-\frac{3-{\mu}}2}D^{\theta}\tilde a(u)\|_{ L^2_{t,x}}
\|r^{2-{\mu}} F_2(u)\|_{L^\infty_{t,x}}
\\
 & \les & 
\|r^{1-{\mu}} \pa u\|_{L^\infty_{t,x}}
( \|r^{-\frac{1-{\mu}}2}D^{\theta}\pa u\|_{  L^2_{t,x}}
C(\|u\|_{L^\infty_{t,x}})+{ \|r^{-\frac{3-{\mu}}2}D^{\theta}\tilde a(u)\|_{ L^2_{t,x}}}
\|r \pa u\|_{L^\infty_{t,x}}
)
\\ & \les
  &
  \tilde C(
\|\pa u \|_{L^\infty_t \dot B^{\frac n 2-1}_{2,1} })
|I_j|^{\frac {\mu} 2} \|D^\theta u\|_{\widetilde{LE}_{{I_j},{\mu}}}
\| \pa u\|_{L^\infty_{t} \dot H^{\frac n 2-1+{\mu}}}
\ ,
\end{eqnarray*}
where we have used
Theorem \ref{thm-wLeib0}
with $\theta=s_2-1\in (0,1]$ when $n=3$,
in the last inequality.

For odd $n\ge 5$ , the inequality still holds for $\mu_2<1/2$.
Actually, as $\theta=s_2-1$ with $k=[\theta]=(n-3)/2\ge 1$.
With $\al=-\frac{3-{\mu}}2$, we see that
$$\al<\frac n 2, k-\al<\frac n 2,$$
and
so we could apply Proposition \ref{thm-wchain-high}, together
with Lemma \ref{thm-trace}, to obtain
\begin{eqnarray}
\|r^{-\frac{3-{\mu}}2}D^{\theta}\tilde a(u)\|_{ L^2_{t,x}} &
 \les  &   C(\max_{j\le k}
 \|r^j \nabla^j u \|_{L^\infty_{t,x}} )
\|r^{-\frac{3-\mu}2}D^{\theta} u\|_{ L^2_{t,x}}\ . 
\nonumber\\
& \les &  C(\|\pa u \|_{L^\infty_t \dot B^{\frac n 2-1}_{2,1} })
\|r^{-\frac{3-\mu}2}D^{\theta} u\|_{ L^2_{t,x}}\ .
 \label{5eq-weightedchain}
\end{eqnarray}
Alternatively, when $\mu_2=1/2$ and so $\theta=(n-1)/2\ge 2$, it could be estimated directly, as follows
\begin{eqnarray*}
 \|r^{-\frac{3-{\mu}}2}D^{\theta}\tilde a(u)\|_{ L^2_{t,x}}
&\les&
\sum_{|\sum \be_l|=\theta, |\be_1|\ge |\be_j|\ge 1}
 \|r^{-\frac{3-{\mu}}2} \Pi_{l=1}^j \nabla^{\be_l} u\|_{ L^2_{t,x}}\\
&\les&
\sum_{1\le |\be_1|\le \theta}
 \|r^{-\frac{3-{\mu}}2+|\be_1|-\theta}
 \nabla^{\be_1} u\|_{ L^2_{t,x}}
   \Pi_{l=2}^j\|r^{|\be_l|} \nabla^{\be_l} u\|_{ L^\infty_{t,x}}
\\
 &\les& C(\|\pa u \|_{L^\infty_t \dot B^{n/2-1}_{2,1} })
\|r^{-\frac{3-\mu}2}D^{\theta} u\|_{ L^2_{t,x}}\ .
\end{eqnarray*}

 For the case of even $n$, 
we have 
$k=[\theta]=(n-2)/2$, 
$\tau=\theta-k$,
$n/2< k-\al< n/2+1$. Let
$q, p\in (2,\infty)$ to be determined, such that $1/q+1/p=1/2$,
the similar argument in Proposition \ref{thm-wchain-high} gives us the desired bound, except the following term
$$\sum_{|\sum \be_l|=k, |\be_l|\ge 1}\|r^{\tau-\frac n q} D^{\tau} (\tilde a^{(j)}(u)-\tilde a^{(j)}(0))\|_{ L^q_{x}}
\|r^{\al-\tau+\frac n q}\Pi_{l=1}^j \nabla^{\be_l} u
\|_{ L^p_{x}}
\ .
$$
As  $-n<\tau q-n<n(q-1)$, we have
$r^{\tau q -n}\in A_q$ and so
\begin{eqnarray*}
\|r^{\tau-\frac n q} D^{\tau} (\tilde a^{(j)}(u)-\tilde a^{(j)}(0))\|_{ L^q_{x}} 
& \les & C(\|\pa u \|_{\dot B^{\frac n 2-1}_{2,1} })
\|r^{\tau -\frac n q} D^{\tau} u\|_{ L^q_{x}} \\
 & \les & C(\|\pa u \|_{\dot B^{\frac n 2-1}_{2,1} })
\| u\|_{\dot B^{\frac n 2}_{2,1}} 
\end{eqnarray*}
where we have used   Theorem \ref{thm-wLeib0} and
Lemma \ref{thm-trace}.
For another term, we
let $p$ be sufficiently close to $2$ such that
$\theta-\be_1
\in (1/2-1/p, n/2)
$.
Because of the assumption that
$\al-\mu_2+2=\frac{\mu+1}2-\mu_2
>0
$, we also have 
$$\al<\frac n 2,\ \al-(\theta-|\be_1|)
>-\frac n 2\ ,
$$
and thus we could apply \eqref{eq-trace-w2} to obtain
\begin{eqnarray*}
\|r^{\al-\tau +\frac n q}\Pi_{l=1}^j \nabla^{\be_l} u
\|_{ L^p_{x}} & \le &
\|r^{\al-\theta+\frac n q+|\be_1|} \nabla^{\be_1} u
\|_{ L^p_{x}}
\|r^{k-|\be_1|}\Pi_{l=2}^j \nabla^{\be_l} u\|_{ L^\infty_{x}}
 \\
& \les & \|r^\al D^\theta u\|_{L^2}\| u\|_{\dot B^{\frac n 2}_{2,1}} ^{j-1}\ .
\end{eqnarray*}
Thus, we still have \eqref{5eq-weightedchain},
which completes the proof.
\end{prf}

In view of
\eqref{5thm-linea} and
Lemma \ref{5thm-nl}, we have
$$\|D^{s_2-1} u_{k+1}\|_{\widetilde{LE}_{{I_j},{\mu}}}\les
\|\pa u(T_j)\|_{\dot H^{s_2-1}} +
|I_j|^{{\mu} } \|D^{s_2-1} u_{k}\|_{\widetilde{LE}_{{I_j},{\mu}}}
\| \pa u_{k}\|_{L^\infty_{t} \dot H^{s-1}}
\  ,
$$
for any $k\ge 2$.
Then, with $\de_3>0$ sufficiently small, we could conclude with the uniform bound
$$\|D^{s_2-1} u_{k+1}\|_{\widetilde{LE}_{{I_j},{\mu}}}\les
\| \pa u(T_j)\|_{\dot H^{s_2-1}}\ ,$$
for any $k\ge 2$, which, combined with
the induction assumption
\eqref{5eq-induction},
 gives us the desired bound
$$\|D^{s_2-1} u\|_{\widetilde{LE}_{{I_j},{\mu}}}\les
\| \pa u(T_j)\|_{\dot H^{s_2-1}} \les
\|\pa u(0)\|_{\dot H^{s_2-1}}\ .$$
As \eqref{5eq-induction} is trivial
when $T_j=0$, by induction, we see that
\eqref{5eq-induction} holds for any $j$ and thus
$$\|D^{s_2-1} u\|_{\widetilde{LE}_{T,{\mu}}}\les
\|  \pa u(0)\|_{\dot H^{s_2-1}}\ .$$
This completes the proof of $\pa u\in L^\infty_t \dot H^{s_2-1}([0,T]\times \R^n)$. As it is true for any $T<T_*$, we see that 
$\pa u\in L^\infty_{loc} \dot H^{s_2-1}([0,T_*)\times \R^n)$.

Notice also that in the case of even $n$, the result could be iterated to show that for any $s_2\in (s, n/2+1)$, we have persistence of regularity.

\subsection{Persistence of regularity for odd $n$}
Now we could prove persistence of regularity to $H^{s_1}$ with $s_1=[(n+4)/2]$.
Let us begin with the case of odd $n$, when
$s_1=(n+3)/2$. 

As we see from Subsection \ref{sec:persist-slight}, we could assume we have 
$H^k$ solution, where $k=(n+1)/2=[(n+2)/2]$ and $\mu=1/2$. Also, it suffices for us to prove
\beeq\label{5eq-induction2}
\|\pa u\|_{L^\infty \dot H^{k}([0,T]\times\R^n)}\les
 1+\|\pa u(0)\|_{\dot H^{k}}
\ ,
\eneq
for any $T$ such that
$$
T^{1/2}
\|\nabla^{k} u\|_{X_{T}}\les
T^{1/2}
\|D^{k-1} u\|_{\widetilde{LE}_{T,{1/2}}}
\ll 1,\ $$
$$
\|u\|_{L^\infty_t H^{k}([0,T]\times\R^n)}+
\|\pt u\|_{L^\infty_t H^{k-1}([0,T]\times\R^n)}
\les 1
\ .$$ For simplicity, we will just illustrate the proof for solutions, instead for the approximate solutions.

By
 \eqref{eq-LE-high-k2}, we have 
\beeq\label{5eq-LE-high-k2}
\|\tilde \pa \nabla^k u\|_{X_T}\les
 \|\pa u(0)\|_{\dot H^{k}}
 +\|\nabla^k F\|_{X_T^*}
+T^{1/2}
\|\nabla^k u\|_{X_T}
\| g(u)\|_{L^\infty_t \dot H^{k+1}}
\ . 
\eneq
Recall the classical Schauder estimates yield
$$
\| g(u)\|_{L^\infty_t  H^{k+1}} 
\les 
C(g, \|u\|_{L^\infty})
\| u\|_{L^\infty_t H^{k+1}} ,$$
which shows that the last term on the right of
 \eqref{5eq-LE-high-k2} is admissible.
 
Then, to finish the proof of
\eqref{5eq-induction2}, we need only to prove 
a nonlinear estimate, for the nonlinear term $\|\nabla^k F\|_{X_T^*}$, which is provided by the following
\begin{lem}\label{5thm-nl-odd}
Let $n$ be odd and $k=(n+1)/2$, $F(u)=a(u) u_t^2+b(u)|\nabla u|^2$, then for radial functions $u$,
\beeq\label{5eq-nl-1}\|r^{1/4}\nabla ^{k} F (u)\|_{ L^2_{t,x}}
\les 
  C(\| \pa u\|_{L^\infty_{t}  H^{k-1}})
 \|r^{-1/4} \nabla^{k}\pa u\|_{L^2} \| \pa u\|_{L^\infty_{t}  H^{k-1}}\ .
\eneq
\end{lem}
\begin{prf} 
At first, when there are no derivatives acting on $a(u)$ or $b(u)$, we need only to control
\begin{eqnarray*}\|r^{1/4}\nabla ^{k} (\pa u)^2\|_{ L^2_{t,x}}&
\les&
\|r^{-1/4}\nabla ^{k} \pa u\|_{ L^2_{t,x}}
\|r^{1/2} \pa u\|_{ L^\infty_{t,x}}\\
&\les&
\|r^{-1/4}\nabla ^{k} \pa u\|_{ L^2_{t,x}}
\|\pa u\|_{ L^\infty_{t}\dot H^{k-1}}\ ,
\end{eqnarray*}
by Theorem \ref{thm-wLeib-fromChain} and \eqref{eq-2.1}.

For the remaining case, thanks to the uniform boundedness of $u$, we are reduced to control
$$
\|r^{1/4} 
\Pi_{j=1}^l \nabla^{\al_j}\pa u
\|_{ L^2_{x}}$$
where $l\ge 3$,
 $\sum |\al_j|=k+2-l$. Without loss of generality, we assume $|\al_j|$ is non-increasing. Notice then that
 $$
k+2-l= \sum |\al_j|\ge 2|\al_2|\Rightarrow
|\al_2|\le \frac{n-1}4\Rightarrow
|\al_2|\le \frac{n-3}2\ ,$$
where we used the fact that $|\al_2|$ must be integer. Then
we see from \eqref{eq-2.1} that
$$\|r^{|\al_j|+1/2}\nabla^{\al_j}\pa u\|_{ L^\infty_{t,x}}
+\|r^{|\al_j|+1} \nabla^{\al_j}\pa u\|_{ L^\infty_{t,x}}
\les
\|\pa u\|_{ L^\infty_{t} H^{k-1}}\ ,$$
 for any $j\ge 2$.
 
 When $|\al_1|\ge 1$, we have 
$-1/4-k+|\al_1|>-n/2$, and 
\begin{eqnarray*}\|r^{1/4} 
\Pi_{j=1}^l \nabla^{\al_j}\pa u
\|_{ L^2_{x}}&
\les&
\|r^{1/4+1/2-\sum_{j\ge 2} (|\al_j|+1)} \nabla^{\al_1}\pa u\|_{L^2}\\&&
\times
\|r^{1/2+|\al_2|}\nabla^{\al_2}\pa u\|_{L^\infty}
\Pi_{j=3}^l \|r^{1+|\al_j|}\nabla^{\al_j}\pa u\|_{L^\infty}
\\
&\les&\|\pa u\|_{ L^\infty_{t} H^{k-1}}^{l-1}
\|r^{-1/4-k+|\al_1|} \nabla^{\al_1}\pa u\|_{L^2}\\
&\les&\|\pa u\|_{ L^\infty_{t} H^{k-1}}^{l-1}
\|r^{-1/4} \nabla^{k}\pa u\|_{L^2}\ ,
\end{eqnarray*}
by \eqref{eq-SteinWeiss}.
 While in the case $|\al_1|=0$, we have
 $l=k+2$, $|\al_j|=0$, and 
\begin{eqnarray*}\|r^{1/4} 
\Pi_{j=1}^l \nabla^{\al_j}\pa u
\|_{ L^2_{x}}&
\les&
\|r^{1/4} 
\Pi_{j=1}^{k+2} \pa u
\|_{ L^2_{x}}\\&\les&
\|\pa u\|_{ L^\infty_{t} H^{k-1}}^{k+1}
\|r^{1/4-(k+1)/2}
\<r\>^{-(k+1)/2}
  \pa u\|_{L^2}\\&\les&
\|\pa u\|_{ L^\infty_{t} H^{k-1}}^{k+1}
\|r^{1/4-k}
  \pa u\|_{L^2}\\
&\les&\|\pa u\|_{ L^\infty_{t} H^{k-1}}^{k+1}
\|r^{-1/4}D^{k}  \pa u\|_{L^2}
\ .
\end{eqnarray*}
This completes the proof.
\end{prf}

\subsection{Persistence of regularity for even $n$}
When $n$ is even, we use similar argument.
Here  $s_1=[(n+4)/2]=n/2+2$, $k=n/2+1$ and $\mu\in (1/2, 1)$.
We need only to prove
\beeq\label{5eq-induction3}
\|\pa u\|_{L^\infty \dot H^{k}([0,T]\times\R^n)}\les
 1+\|\pa u(0)\|_{\dot H^{k}}
\ ,
\eneq
for any $T\ll 1$ such that
$$
T^{\mu}
\|D^{n/2-1+\mu} \pa u\|_{X_{T}}
\ll 1\ ,
\|(u,\pa u)\|_{ L^\infty H^{n/2-1+\mu}}+
\|\pa u\|_{ X_T}\les 1\ .
 $$

By \eqref{eq-LE-high-k3}, we have
\beeq\label{5eq-LE-high-k3}
\|\tilde \pa \nabla^k u\|_{X_T}\les
 \|\pa u(0)\|_{\dot H^{k}}+
\|\nabla^k F\|_{X_T^*}
+T^{\mu}
 \|g(u)\|_{L^\infty_t \dot H^{k+1}}\|D^{n/2+\mu}u\|_{X_T}
\ , 
\eneq
where, as before, the last term is admissible, thanks to Schauder estimates.
 
 As for the nonlinear term, we have the following estimates, which is sufficient to conclude the proof of
\eqref{5eq-induction3}. 
\begin{lem}\label{5thm-nl-even}
Let $n\ge 4$ be even, $k=n/2+1$, $\mu=2/3$, 
and 
$F(u)=a(u) u_t^2+b(u)|\nabla u|^2$, then for radial functions $u$,
\beeq\label{5eq-nl-2}
\|\nabla^k F\|_{X_T^*}
\les 
  C(
\|\pa u \|_{L^\infty_t H^{\frac n 2-1+\mu}})
T^{ {\mu}} 
(\|\pa u\|_{X_T}+
\|D^k\pa  u\|_{X_T})
\| \pa u\|_{L^\infty_{t} H^{\frac n 2-1+{\mu}}}\ .\eneq
\end{lem}

\begin{prf} The proof proceeds similar as
that of Lemma \ref{5thm-nl-odd}.
At first, when there are no derivatives acting on $a(u)$ or $b(u)$, we need only to control
\begin{eqnarray*}
\|\nabla^k (\pa u)^2\|_{X_T^*}
&
\les&
T^{\mu/2} \|r^{-(1-{\mu})/2}\nabla ^{k} \pa u\|_{ L^2_{t,x}}
\|r^{1-{\mu}} \pa u\|_{ L^\infty_{t,x}}\\
&\les&
T^{\mu}\|\nabla ^{k} \pa u\|_{X_T}
\|\pa u\|_{ L^\infty_{t}  \dot H^{n/2-1+{\mu}}}\ ,
\end{eqnarray*}
by Theorem \ref{thm-wLeib-fromChain} and \eqref{eq-2.1}.

For the remaining case, thanks to the uniform boundedness of $u$, we are reduced to control
$$
\|r^{(1-\mu)/2} 
\Pi_{j=1}^l \nabla^{\al_j}\pa u
\|_{ L^2_{x}}$$
where $l\ge 3$,
 $\sum |\al_j|=k+2-l$ and $|\al_j|$ is non-increasing.

 When $|\al_2|=0$,
 we see from \eqref{eq-2.1} that
$$\|r^{1-\mu}\pa u\|_{ L^\infty_{t,x}}
\les
\|\pa u\|_{ L^\infty_{t} \dot H^{n/2-1+\mu}}
\ ,$$
 for any $j\ge 2$.
 As $3\le l\le k+2$, $
|\al_1|=k+2-l$, 
 $
k-\mu(l-2)\in [(1-\mu)k, k]$,
and $(1-\mu)/2-(l-1)(1-\mu)>-n/2$,
we have
\begin{eqnarray*}
\|r^{(1-\mu)/2} 
\Pi_{j=1}^l \nabla^{\al_j}\pa u
\|_{ L^2_{x}}
 & \les & 
 \|r^{(1-\mu)/2-(l-1)(1-\mu)} \nabla^{\al_1}\pa u
\|_{ L^2_{x}}
\Pi_{j=2}^l\|r^{1-\mu} \pa u
\|_{ L^\infty_{x}}
 \\
 & \les & 
 \|r^{-(1-\mu)/2}D^{(l-2)(1-\mu)} \nabla^{\al_1}\pa u
\|_{ L^2_{x}}
\|\pa u\|_{  \dot H^{n/2-1+\mu}}^{l-1}
 \\
 & \les & 
 \|r^{-(1-\mu)/2}D^{k-\mu(l-2)} \pa u
\|_{ L^2_{x}}
\|\pa u\|_{  \dot H^{n/2-1+\mu}}^{l-1}\ .
\end{eqnarray*}

On the other hand, if 
 $|\al_1|= |\al_2|=1$, 
as
 $\mu>1/2$, we have
 $n/2+\mu-1-|\al_2|>1/2$,
  and then
 $$
\|r^{|\al_j|+1-\mu} \nabla^{\al_j}\pa u
\|_{ L^\infty_{x}}
\les
\|\pa u\|_{  \dot H^{n/2-1+\mu}}\ ,
$$ for any $j\ge 2$.
 Thus
 \begin{eqnarray*}
&&\|r^{(1-\mu)/2} 
\Pi_{j=1}^l \nabla^{\al_j}\pa u
\|_{ L^2_{x}}\\
 & \les & \|r^{(1-\mu)/2-(1-\mu)(l-1)-(k+2-l-|\al_1|)}\nabla^{\al_1}\pa u
\|_{ L^2_{x}}
\Pi_{j=2}^l\|r^{|\al_j|+1-\mu} \nabla^{\al_j}\pa u
\|_{ L^\infty_{x}}
 \\
 & \les & 
 \|r^{-(1-\mu)/2+\mu(l-2)-n/2}
 \nabla \pa u
\|_{ L^2_{x}} 
\|\pa u\|_{  \dot H^{n/2-1+\mu}}^{l-1}
 \\
 & \les & 
 \|r^{-(1-\mu)/2}D^{k-\mu(l-2)} \pa u
\|_{ L^2_{x}}
\|\pa u\|_{\dot  H^{n/2-1+\mu}}^{l-1}\ ,
\end{eqnarray*} where in the last inequality, we have used the fact
$\mu(l-2)\ge \mu>(1-\mu)/2$, thanks to $\mu\in (1/2, 1)$, so that we could apply
\eqref{eq-SteinWeiss}.

 It remains to consider the case
 $|\al_1|\ge 2$,
then
$$-(1-\mu)/2-(k- |\al_1|)>-n/2,
$$
and so
$$\|r^{-(1-\mu)/2-(k-|\al_1|)}\nabla^{\al_1}\pa u
\|_{ L^2_{x}}\les
\|r^{-(1-\mu)/2}\nabla^{k}\pa u
\|_{ L^2_{x}}\ .
$$
Also,
notice that
 $$|\al_2|\le
k+2-l-|\al_1|\le k+2-3-2=\frac{n}2-2
\ ,$$
we see from \eqref{eq-2.1} that
$$
\|r^{|\al_2|+1-\mu} \nabla^{\al_2}\pa u
\|_{ L^\infty_{x}}
\les
\|\pa u\|_{  \dot H^{n/2-1+\mu}},\ 
\|r^{|\al_j|+1}\nabla^{\al_j}\pa u
\|_{ L^\infty_{x}}
\les
\|\pa u\|_{  \dot H^{n/2-1}}\ ,
$$
 for any $j\ge 3$.
 Thus
 \begin{eqnarray*}
&&\|r^{(1-\mu)/2} 
\Pi_{j=1}^l \nabla^{\al_j}\pa u
\|_{ L^2_{x}}\\
 & \les & \|r^{-(1-\mu)/2-(k-|\al_1|)}\nabla^{\al_1}\pa u
\|_{ L^2_{x}}
\|r^{|\al_2|+1-\mu} \nabla^{\al_2}\pa u
\|_{ L^\infty_{x}}
\Pi_{j=3}^l\|r^{|\al_j|+1} \nabla^{\al_j}\pa u
\|_{ L^\infty_{x}}
 \\
 & \les & 
 \|r^{-(1-\mu)/2}\nabla^{k}\pa u
\|_{ L^2_{x}} 
\|\pa u\|_{  \dot H^{n/2-1+\mu}}
\|\pa u\|_{  \dot H^{n/2-1}}^{l-2}
 \\
 & \les & 
 \|r^{-(1-\mu)/2}D^{k} \pa u
\|_{ L^2_{x}}
\|\pa u\|_{  H^{n/2-1+\mu}}^{l-1}\ .
\end{eqnarray*}
This completes the proof.
\end{prf}

\subsection{Continuous dependence}
The continuous dependence property is essentially included in the proofs of convergence of the approximate solutions,   Lemma \ref{thm-conv0}, and the unconditional uniqueness.

Let $T_*>0$ be the lifespan of the solution $u$, with data $(u_0, u_1)$.
Fix $T<T_*$ and
$s_1\in (s_c, s)$,
we have uniform bound on $\|\pa u\|_{L^\infty_t  H^{s-1}
([0,T]\times\R^n)}$.
When $n=3$, as $s_0<s-1$, without loss of generality, we could assume 
$s_0=s_1-1$ and also the
uniform bound
$\|\pa u\|_{L^\infty_t  \dot H^{s_1-2}}$.
As the proof for $n\ge 4$ is relatively easier, we present only the proof for $n=3$.

\subsubsection{Short time continuity}
Before proving the full continuous dependence property, we present a result of short time continuous dependence, for data with regularity
$\tau_0$ to solution with regularity
$\tau_1$, with $s\ge \tau_0>\tau_1\ge s_1$.
Suppose
$\|\pa u(0)  \|_{\dot H^{s-1}\cap \dot H^{s_1-2}}\le M<\infty$, 
$\tau_0-\tau_1\ge \ep>0$,
we would like to find $T>0$, with the following property:
for any $\ep>0$,
there exists $\de(\ep)>0$,
such that 
whenever
$\|(\nabla (u_0-v_0), u_1-v_1)\|_{\dot H^{\tau_0-1}\cap \dot H^{s_1-2}}\le \de
$, the corresponding solution $v\in L^\infty H^{\tau_0}\times C_t^{0,1} H^{{\tau_0}-1}$ is well-defined in $[0, T]\times \R^3$ and 
 $$
\|\pa (u-v)\|_{L^\infty (\dot H^{\tau_1-1}
\cap \dot H^{s_1-2})}\le\ep\ .$$
Here, $T$ could be chosen to be independent of the specific choices of $\tau_0$, $\tau_1$.

At first, by assuming $\de\le 1$, we could always assume
$\|\pa v(0)  \|_{\dot H^{\tau_0-1}\cap \dot H^{s_1-2}}\le M+1<\infty$.
Based on this, we know that
$$(\ep_{\tau_0}+\ep_{\tau_0-1})^{-1/(\tau_0-s_c)}
\ge (2(M+1))^{-1/(s_1-s_c)}>0\ ,
$$
and thus, in view of Lemma \ref{thm-lwp0} and
\eqref{eq-life-poly3}, the
corresponding solution $v\in L^\infty H^{\tau_0}\times C_t^{0,1} H^{{\tau_0}-1}\cap
C_t \dot H^{s_1-1} \cap C^{0,1}_t\dot H^{s_1-2}$ is well-defined in $[0, T_4]\times \R^3$,
together with a uniform bound in $L^\infty \dot H^{\tau_0}\times C_t^{0,1} \dot H^{{\tau_0}-1}$,
where
\beeq T_4:= \min(\de_4, c (2(M+1))^{-1/(s_1-s_c)}, T_0).\eneq

We need to give the estimate of 
 $$
\|\pa (u-v)\|_{L^\infty (\dot H^{\tau_1-1}
\cap \dot H^{s_1-2})}\ ,$$
in terms of the norm  of $\pa(u-v)(0)$. 
For this purpose, we give firstly the estimate of
 $
\|\pa (u-v)\|_{L^\infty  \dot H^{s_1-2}}$.
Let $w  =  u- v$ and $\mu_0=\tau_0-s_c$, it satisfies
\begin{align*}
\Box w +& \nabla \cdot (g(u)
  \nabla w)    =
   \nabla \cdot ((g(v)-g(u))
  \nabla v) \\
    & +  
  \nabla ( g(u) )\cdot
  \nabla u-
    \nabla ( g(v) )\cdot
  \nabla v
 + F( u)- F( v):=G(u, w)\ ,
\end{align*}
together with $w (0,x)=u_0-v_0 $, $\pt w (0,x)=u_1-v_1$.
Notice that $T^{\mu_0}\| g(u)\|_{L^\infty_{t} \dot H^{\tau_0}}\les 1$,
then by  Proposition \ref{th-LEkey3}, and arguing as in
 Lemma \ref{thm-conv0}, we obtain that for any $R\in (0, T_4]$, 
\begin{eqnarray*}&&\|D^{s_1-2}\pa w\|_{X_{R}}\\  &\les &
\| \pa w(0)\|_{\dot H^{s_1-2}}+
R^{\mu_0}\| g(u)\|_{L^\infty_{R} \dot H^{\tau_0-1}}\| \pa w_{j}(0)\|_{\dot H^{s_1-1}}+ \|D^{s_1-2} G(u,w )\|_{X_R^*}
 \\
& \les & \| \pa w (0)\|_{\dot H^{s_1-2}\cap \dot H^{s_1-1}}\\&&+
 C(g,a,b, \|\pa (u, v)\|_{L^\infty_{t\in [0,T]}\dot B^{s_c-1}_{2,1}}) R^{\mu_0} \|D^{s_1-2}\pa w\|_{X_R}
\|   \pa (u, v)\|_{L^\infty_{t\in [0, T]} \dot H^{\tau_0-1}}
 \\
& \les & 
\| \pa w(0)\|_{\dot H^{s_1-2}\cap \dot H^{s_1-1}}+
 R^{\mu_0} \|D^{s_1-2}\pa w\|_{X_R}
\ ,
\end{eqnarray*}
where the implicit constant is independent of $R\in (0,T_4]$.
Thus by choosing $R$ small enough, we obtain
$$
\|\pa w \|_{L^\infty( [0, R] ; \dot H^{s_1-2})}\les
 \|D^{s_1-2}\pa w \|_{X_R}\les \| \pa w (0)\|_{\dot H^{s_1-2}\cap \dot H^{s_1-1}}\ .$$
 Iterating this argument finite many times ($\sim T_4/R$), we obtain
\beeq\label{eq-cont-1}\|\pa w \|_{L^\infty( [0, T_4] ; \dot H^{s_1-2})}\les
 \| \pa w (0)\|_{\dot H^{s_1-2}\cap \dot H^{s_1-1}}\ .\eneq

Combined with the uniform bounds, as that in 
 Lemma \ref{thm-lwp0},
 $$\|\pa w \|_{L^\infty( [0, T_4] ; \dot H^{\tau_0-1})}\le
 \|\pa u\|_{L^\infty( [0, T_4] ; \dot H^{\tau_0-1})}+\|\pa v \|_{L^\infty( [0, T_4] ; \dot H^{\tau_0-1})}\le
  2 C_2 (M+1)
 \ ,$$
 we obtain, for any $t\in [0, T_4]$,
\beeq\label{eq-cont-2}\|\pa w(t) \|_{ \dot H^{\tau_1-1}}\le
\|\pa w(t) \|_{ \dot H^{\tau_0-1}}^{1-\theta}
\|\pa w(t) \|_{ \dot H^{s_1-2} }^{\theta}
\les 
 \| \pa w (0)\|_{\dot H^{s_1-2}\cap \dot H^{s_1-1}}^{\theta}\ ,\eneq
 where
 $\tau_0(1-\theta)+ (s_1-1)\theta=\tau_1$.

\subsubsection{Long time continuity}
With short time continuity available, it is easy to conclude long time continuity. Actually,
as $T<T_*$, there exists $M<\infty$ such that
$$\|\pa u  \|_{L^\infty([0,T]; \dot H^{s-1}\cap \dot H^{s_1-2})}\le M<\infty\ .$$
For fixed $s_1\in (s_c, s)$, we have uniform $T_4$ so that we have
short time continuity, in any interval with length less than $T_4$, around the solution $u$. Thus, we could divide $[0, T]$ into finite, say, $N$, adjacent intervals $\{I_j\}_{j=1}^N$, with
$|I_j|<T_4$,
$I_j=[t_{j-1}, t_{j}]$,
$t_0=0$, $t_N=T$.

 Let 
 $\tau_j=s-j(s-s_1)/N$, we could apply
 short time continuity, from $t_{j-1}$ to the interval $I_j$. Gluing together, we obtain the 
 long time continuity.

\section{Three dimensional almost global existence with small data}\label{sec-7-n=3}

In this section, when $n=3$, we show that the lower bound of the lifespan, available from local results, could be improved to almost global existence, 
Theorem \ref{th-2-al-global}. Without loss of generality, we assume
$s=3/2+\mu$ with $\mu\in (0, 1/2]$ and the solution lies in $CH^s\cap C^1 H^{s-1}$.

Let $I\subset J:= [0, T_*(u_0,u_1))$ be the subset such that for any $T\in I$, we have
\beeq\label{eq-continuity}\|D^{s-1} u\|_{LE_{T}}\le 10C_3 \ep_s,\ 
\| u\|_{LE_{T}}\leq 10 C_3 \ep_1 
\ ,\eneq
where $C_3$ denotes the constant in \eqref{eb16}.
It is clear that $I$ is non-empty and closed set in $J$. By bootstrap argument,
to show existence up to
$\exp(c_1/(\ep_1+\ep_s))$,
 we need only to show that \eqref{eq-continuity} holds for $5C_3$ instead of $10 C_3$, 
for any $T\in I\cap [0, 
\exp(c_1/(\ep_1+\ep_s))]$, provided that
$\ep_1+\ep_s<\de$ for some sufficiently small $\de>0$.

By Sobolev embedding, we see that
$$
\|u\|_{L^\infty_{t,x}(S_T)}\le C \|\nabla u\|_{L^\infty H^{s-1}(S_T)}\le 10 CC_3
(\ep_1+\ep_s)\les 1
\ ,$$
$g'(u)=\O(1)$,
and  so
$$\|g(u)\|_{L^\infty(S_T)}\le
\|u\int_0^1 g'(\la u)d\la \|_{L^\infty(S_T)}\les \ep_1+\ep_s\ .$$
Moreover, we have
\beeq\label{eq-4.5}\|r^{1-\mu}\pa g(u)\|_{L^{\infty}}\les
\|r^{1-\mu}\pa u\|_{L^{\infty}}\les
 \|\pa u\|_{L^{\infty} \dot H^{s-1}}
\les
\ep_s\ ,\eneq
\beeq\label{eq-4.6}\|r\pa g(u)\|_{L^{\infty}}\les
\|r \pa u\|_{L^{\infty}}\les
\|\pa u\|_{L^{\infty}_t H^{s-1}}\les
\ep_1+\ep_s\ .
\eneq
From these estimates, we see that \eqref{eb15} is satisfied when $T\le
\exp(c/\ep_c)$ with
$c+\ep_1+\ep_s\ll 1$.

Recall that $u$ is constructed through approximation of $C^\infty_t C_c^\infty$ solutions of approximate equations,
Proposition \ref{th-LEkey2} applies for $u$ as well, which gives us
\beeq\label{eq-continuity2}\|D^{s-1} u\|_{LE_{T}}\leq C_3 
\ep_s+
C_3
(\ln\<T\>)^{\frac 12}
\|r^{\frac{1-\mu}2}
\<r\>^{\frac{\mu}2}
D^{s-1} F(u)\|_{L^2_{t,x}} 
\ ,\eneq
\beeq\label{eq-continuity3}
\| u\|_{LE_{T}}\le C_3\ep_1
+C_3
(\ln\<T\>)^{\frac 12}
\|r^{\frac{1-\mu}2}
\<r\>^{\frac{\mu}2}
F(u)\|_{L^2_{t,x}} 
\ .\eneq

When $F(u)=a(u)u_t^2$, in view of \eqref{eq-4.5} and \eqref{eq-4.6},
we have
\begin{eqnarray*}
\|r^{\frac{1-\mu}2}
\<r\>^{\frac{\mu}2}
F(u)\|_{ L^2_{t,x}} 
 & \les &
\|r^{-\frac{1-\mu}2}
\<r\>^{-\frac{\mu}2}
u_t\|_{   L^2_{t,x}} 
\|r^{1-\mu}
\<r\>^{\mu}
a(u) u_t\|_{L^\infty_{t,x}} 
 \\
&\les&
(\ln \<T\>)^{1/2} \ep_1 (\ep_1+\ep_s)\ .
\end{eqnarray*}
For the term with $D^{s-1}=D^{1/2+\mu}$,
 by Theorem \ref{thm-wLeib-fromChain} and
 Theorem \ref{thm-wLeib0}, together with Lemma \ref{thm-Ap0}, we have
\begin{eqnarray*}
&&\|r^{\frac{1-\mu}2}
\<r\>^{\frac{\mu}2}
D^{s-1} F(u)\|_{  L^2_{t,x}} \\
 & \les &
\|r^{-\frac{1-\mu}2}
\<r\>^{-\frac{\mu}2}
D^{s-1} u_t\|_{ L^2_{t,x}} 
\|r^{1-\mu}
\<r\>^{\mu}
(|a(u)|+|a(0)|)u_t\|_{L^\infty_{t,x}}  \\
&  & +
\|r^{-\frac{3}2(1-\mu)}
\<r\>^{-\frac{\mu}2}
D^{s-1}  (a(u)-a(0))\|_{L^2_{t,x}} 
\|r^{2(1-\mu)}
\<r\>^{\mu}
u_t^2\|_{L^\infty_{t,x}} 
 \\
 & \les &
\|r^{-\frac{1-\mu}2}
\<r\>^{-\frac{\mu}2} D^{s-1} u_t\|_{L^2_{t,x}} 
\|\pa u\|_{L^\infty_t  H^{s-1}}
\\&&+
\|r^{-\frac{1-\mu}2}
\<r\>^{-\frac{\mu}2} 
r^{\mu-1}D^{s-1} u\|_{L^2_{t,x}} 
\|\pa u\|_{L^\infty_t  \dot H^{s-1}}
\|\pa u\|_{L^\infty_t  H^{s-1}}
 \\
 & \les &
(\ln \<T\>)^{1/2} \ep_s (\ep_1+\ep_s)+
\|r^{-\frac{1-\mu}2}
\<r\>^{-\frac{\mu}2} 
D^{\frac 3 2} u \|_{L^2_{t,x}} 
\|\pa u\|_{L^\infty_t  \dot H^{s-1}}\\
&\les&
(\ln \<T\>)^{1/2} \ep_s (\ep_1+\ep_s) \ ,
\end{eqnarray*}
where in the second to the last inequality, we have used
Lemma \ref{thm-w-Hardy}.

Then, combined with 
\eqref{eq-continuity2} and
\eqref{eq-continuity3}, we arrive  at
$$\|D^{\la-1} u\|_{LE_{T}}
\leq C_3 
 \ep_\la +
C
(\ep_1+\ep_s)\ep_\la\ln\<T\>,\ 
\la=1,s
\ ,
$$
and so
$$\|D^{s-1} u\|_{LE_{T}}\le  5 C_3 \ep_s,\ 
\| u\|_{LE_{T}}\le 5 C_3  \ep_1  \ ,$$
for any
$T\in I\cap [0,
\exp(c_1/(\ep_1+\ep_s))]$,
where $c_1=\min( c, 1/(4C))$.

\section{High dimensional global  well-posedness}
\label{sec-8}

In this section, we show that 
when $\ep_s+\ep_{1}$ is small enough,
the lower bound of the lifespan could be improved to global existence, 
when $n\ge 4$.

For any $s>s_c$,
there exists $\mu\in (0, 1/3)$ such that $s> s_c+\mu$.
Without loss of generality, we assume
$s= s_c+\mu$ and the solution lies in $CH^s\cap C^1 H^{s-1}$.
Let $I\subset J:= [0, T_*(u_0,u_1))$ be the subset such that for any $T\in I$, we have
\beeq\label{6eq-continuity}\|D^{s-1} u\|_{LE_{T}}\le 10 C\ep_s,\ 
\| u\|_{LE_{T}}\leq 10 C  \ep_1  
\ ,\eneq
where $C$ is the constant occurred in Proposition \ref{th-LEkey4} \eqref{62eq-LE-high-k}. 
It is clear that $I$ is non-empty and closed set in $J$. By bootstrap argument,
to show global existence,
 we need only to show that \eqref{6eq-continuity} holds for $5C$ instead of $10 C$, 
for any $T\in I$, provided that
$\ep_1+\ep_s<\de$ for some sufficiently small $\de>0$.

By Sobolev embedding, we see that
$$
\|u\|_{L^\infty_{t,x}(S_T)}\les \|\pa u\|_{L^\infty  H^{s-1} (S_T)}\les
\ep_1+\ep_s
\ ,$$
$g'(u)=\O(1)$,
and  so
$$\|\<r\>^{\mu_1} g(u)\|_{L^\infty(S_T)}\les
\|\<r\>^{\mu_1} u  \|_{L^\infty(S_T)}\les
\|u\|_{L^\infty (\dot H^{\frac n 2-\mu_1}\cap \dot H^{\frac n 2+\mu_1})(S_T)}
\les \ep_1+\ep_s\ ,$$
provided that $\mu_1\le \mu$.
Moreover, we have
\beeq\label{6eq-4.5}\|r^{1-\mu}\pa g(u)\|_{L^{\infty}}\les
\|r^{1-\mu}\pa u\|_{L^{\infty}}\les
 \|\pa u\|_{L^{\infty} \Hsd{s-1}}
\les
\ep_s \ ,\eneq
\beeq\label{6eq-4.6}\|r^{1+\mu_1}\pa g(u)\|_{L^{\infty}}\les
\|r^{1+\mu_1} \pa u\|_{L^{\infty}}\les
\|\pa u\|_{L^{\infty}_t \dot H^{n/2-1-\mu_1}}\les
 \ep_1+\ep_s \ ,
\eneq
and for
any $0\le j\le [(n-1)/2]$,
$$\|r^j \nabla^j u\|_{L^\infty_x}\les \|u\|_{\dot B^{n/2}_{2,1}}\les \ep_1+\ep_s\ 
.$$
From these estimates, we see that \eqref{62eq-LE-frac-h} is satisfied when 
$\ep_1+\ep_s\ll 1$.

Recall that $u$ is constructed through approximation of $C^\infty_t C_c^\infty$ solutions of approximate equations,
Proposition \ref{th-LEkey4} \eqref{62eq-LE-high-k} applies for $u$ as well, which gives us
\beeq\label{6eq-continuity2}\| u\|_{LE_{T}}\le C
\ep_1+C
\|w 
 F(u)\|_{L^2_{t,x}} 
\ ,\eneq
and in the case of odd $n$,
\beeq\label{6eq-continuity3}
\|D^{s-1} u\|_{LE_{T}}\le C\ep_s
+C
\|w 
D^{s-1}  F(u)\|_{L^2_{t,x}} 
\ ,\eneq
where we set $w=r^{\frac{1-\mu}2}
\<r\>^{ \frac{\mu+\mu_1}2 }$.
We claim that the following variant of
\eqref{6eq-continuity3}
\beeq\label{6eq-continuity3'}
\|D^{s-1} u\|_{LE_{T}}\le C\ep_s
+C
\|w 
D^{s-1}  F(u)\|_{L^2_{t,x}} +\tilde C\ep_s(\ep_1+\ep_s)
\ ,\eneq
for some $\tilde C$,
 applies for even $n$ as well. Before presenting the proof of \eqref{6eq-continuity3'},  let us use it to conclude the global existence.

At first, we have
\begin{eqnarray*}
\|w
 F(u)\|_{ L^2_{t,x}} 
 & \les &
\|w^{-1}  \pa u\|_{ L^2_{t,x}} 
\|w^2
(|a(u)|+|b(u)|)\pa u\|_{L^\infty_{t,x}}   \\
 & \les &
\|w^{-1}
\pa  u\|_{  L^2_{t,x}} 
C(\| u\|_{L^\infty_{t,x}})
\|\pa u\|_{L^\infty_t  H^{s-1}}  \\
&\les&
 (\ep_1+\ep_s)\ep_1 \ .
\end{eqnarray*}
Concerning the part with $D^{s-1}$,
when $F(u)=a(u)u_t^2=(a(0)+\tilde a(u))u_t^2$ and $n$ is odd,
by Theorem \ref{thm-wLeib-fromChain},
Proposition \ref{thm-wchain-high} with $[s-1]=k=(n-3)/2$
and $k+ (1-\mu)+ (1+\mu_1)/2<n/2$,
 together with Lemma \ref{thm-Ap0}, 
 we see that
 $\|w
D^{s-1} F(u)\|_{ L^2_{t,x}}$ is controlled by
\begin{eqnarray*}
 &  &
\|w^{-1}
D^{s-1} u_t\|_{ L^2_{t,x}} 
\|w^2
u_t\|_{L^\infty_{t,x}}  +
\|r^{-(1-\mu)}
w^{-1}
D^{s-1} \tilde a(u)\|_{ L^2_{t,x}} 
\|r^{1-\mu}
w
u_t^2\|_{L^\infty_{t,x}} 
 \\
 & \les &
\|D^{s-1} u\|_{LE_T} 
\|\pa u\|_{L^\infty_t   H^{s-1}}  +
\|r^{-(1-\mu)}
w^{-1}
D^{s-1} u\|_{ L^2_{t,x}} \|\pa u\|_{L^\infty_t \dot H^{s-1} }
 \\
 & \les &
\|D^{s-1} u\|_{LE_T} 
\|\pa u\|_{L^\infty_t   H^{s-1}}  +
\|w^{-1}
D^{s-\mu} u\|_{ L^2_{t,x}}  \|\pa u\|_{L^\infty_t \dot H^{s-1} }
\\
&\les&(\ep_1+\ep_s)\ep_s\ ,
\end{eqnarray*}
where in the second  inequality, we have used
Lemma \ref{thm-w-Hardy}.

When $n$ is even, 
we have $[s-1]=n/2-1$, 
and we could apply Proposition \ref{thm-wchain-high} only if 
$\mu>(1+\mu_1)/2$.
For the remaining case $0<\mu\le (1+\mu_1)/2$, we notice that
 $1-\mu+\frac{1+\mu_1}2<n/2$, and we could apply Lemma 
\ref{thm-w-Hardy} to obtain
$$\|w^{-1}r^{ \mu-1} D^{s-1} \tilde a(u)\|_{L^2_{t,x}}
\les\|w^{-1} D^{s-\mu} \tilde a(u)\|_{L^2_{t,x}}
\les\|w^{-1} \nabla^{\frac n 2}\tilde a(u)\|_{L^2_{t,x}}\ .
 $$ Notice that
$$|\nabla^{\frac n 2} \tilde a(u)|\les \sum_{|\sum \be_l|=\frac n2, |\be_1|\ge |\be_l|\ge 1}
\Pi_{l=1}^j |\nabla^{\be_l} u|
\les
\sum_{|\be_l|<n/2, l\ge 2}
r^{|\be_1|-\frac n2}|\nabla^{\be_1} u|
\Pi_{l=2}^j |r^{|\be_l|}\nabla^{\be_l} u|\ ,
$$
we get
\beeq\label{eq-even-fractional}
\|w^{-1}\nabla^{\frac n 2} \tilde a(u)\|_{L^2_{t,x}}\les
\sum_{1\le j\le \frac n2}
\|w^{-1}r^{j-\frac n2}\nabla^{j} u\|_{L^2_{t,x}}
\les\|w^{-1}D^{\frac n2} u\|_{L^2_{t,x}}
\les \ep_1+\ep_s\ ,
\eneq
and thus we have the same estimate as for the odd spatial dimension. 
 
Then, combined with 
\eqref{6eq-continuity2} and
\eqref{6eq-continuity3}, we arrived at
$$\|D^{s-1} u\|_{LE_{T}}\le C
\ep_s+\tilde C
\ep_s(\ep_1+\ep_s)
\ ,
\|  u\|_{LE_{T}}\le 
C\ep_1+\tilde C
\ep_1(\ep_1+\ep_s)
\ .
$$
Consequently, with $\ep_1+\ep_s\ll 1$, we have
$$\|D^{s-1} u\|_{LE_{T}}\le 2C
\ep_s
\ ,
\| u\|_{LE_{T}}\le 2
C\ep_1
\ ,
$$ for any
$T\in I$,
which gives us
 \eqref{6eq-continuity} holds for $2C$.
 By continuity, we see that $T_*=\infty$ and this completes the proof.

\subsection{
(8.6) for even spatial dimension}
In the case of even $n$ with $s=n/2+\mu\ge 2$, we could apply
\eqref{62eq-LE-high-k} 
with $\theta=s-2$, for the equation of $\nabla u$,
$$(\Box + g(u)\Delta)\nabla u = \nabla F(u)-(\nabla g(u))\Delta u\ ,$$
which gives us
\beeq\label{6eq-1}\|D^{s-1} u\|_{LE_{T}}
  \les 
\ep_s+\|w
D^{s-1}  F(u)\|_{L^2_{t,x}} 
+\|w
D^{s-2}(
(\nabla g(u))\Delta u)\|_{L^2_{t,x}} \ .\eneq

When $n\ge 6$, we have
$n/2-2-\mu_1>1/2$ so that
$$\|w^2 r\Delta u\|_{L^\infty_x}\les \|\Delta u\|_{\Hs{s-2}}\les \ep_1+\ep_s,\ 
\|r^{1-\mu}\nabla u\|_{L^\infty_{t,x}} \les \ep_s\ .
$$
 Moreover, by Lemma 
\ref{thm-w-Hardy}, as $2-\mu+(1+\mu_1)/2<n/2$, we have
the following similar estimate as that of
\eqref{eq-even-fractional},
$$\|w^{-1}r^{\mu-2} D^{s-2} g'(u)\|_{L^2_{t,x}}
\les\|w^{-1} D^{\frac n 2} g'(u)\|_{L^2_{t,x}}
\les\|w^{-1}D^{\frac n2} u\|_{L^2_{t,x}}
\les \ep_1+\ep_s\ ,
$$ which gives us
\beeq\label{eq-8.6-claim}\|w^{-1}r^{\mu-2} D^{s-2} g'(u)\|_{L^2_{t,x}}
\les \ep_1+\ep_s\ .
\eneq
Then, by  Theorem \ref{thm-wLeib-fromChain}, we get
\begin{eqnarray*}
\|w D^{s-2}((\nabla g(u))\Delta u)\|_{L^2_{t,x}} 
& \les & 
\|w^{-1} D^{s-2}\Delta u\|_{L^2_{t,x}} 
\|w^2  g'(u)\nabla u\|_{L^\infty_{t,x}}\\&& +
\|w^{-1}r^{-1} D^{s-2}\nabla u\|_{L^2_{t,x}} 
\|w^2 r g'(u)\Delta u\|_{L^\infty_{t,x}} \\
&&+
\|w^{-1}r^{\mu-2} D^{s-2}g'(u)\|_{L^2_{t,x}} 
\|w^2 r^{2-\mu} \nabla u\Delta u\|_{L^\infty_{t,x}}
 \\
 & \les & \ep_s (\ep_1+\ep_s)\ .
\end{eqnarray*}

Turning to the case of  $n=4$, for which we have $s=2+\mu$. Let $q=2/(1-2\mu)$ so that $1/q+\mu=1/2$, and
$$\|  r^{3\mu} \Delta u\|_{L^\infty_t L^q_x} \les
\|  D^{\mu} \Delta u\|_{L^\infty_t L^2_x} \les \ep_s\ .$$
Moreover, we claim that
\beeq\label{eq-8.6-claim2}\|w r^{-3\mu} D^{\mu}\nabla g(u)\|_{L^2_{t}L^{1/\mu}_{x}} \les \ep_1+\ep_s\ .\eneq
Thus we have
\begin{eqnarray*}
\|w D^{\mu}((\nabla g(u))\Delta u)\|_{L^2_{t,x}} 
& \les & 
\|w^{-1} D^{\mu}\Delta u\|_{L^2_{t,x}} 
\|w^2 \nabla g(u)\|_{L^\infty_{t,x}}\\&& +
\|w r^{-3\mu} D^{\mu}\nabla g(u)\|_{L^2_{t}L^{1/\mu}_{x}} 
\|  r^{3\mu} \Delta u\|_{L^\infty_t L^q_x} 
 \\
 & \les & \ep_s (\ep_1+\ep_s)\ .
\end{eqnarray*}
It remains to give the proof of the claim \eqref{eq-8.6-claim2}.
Actually, notice that
$$w r^{-3\mu}=r^{ 1-4\mu }
\<r\>^{  \mu+\mu_1} w^{-1}
\les
r^{3(\frac 1 2-\mu)}( w^{-1} r^{-\frac 1 2-\mu}
+w^{-1}r^{\mu_1-\frac 1 2}) \ ,
$$ 
an application of  Lemma \ref{thm-w-Hardy-trace} gives us that
$$\|w r^{-3\mu} D^\mu \nabla g(u)\|_{L^2_{t}L^{1/\mu}_{x}}
\les
\|w^{-1}r^{-\frac 12-\mu}D^{\frac 12}\nabla g(u)\|_{L^2_{t,x}}
+\|w^{-1}r^{\mu_1-\frac 12}D^{\frac 12}\nabla g(u)\|_{L^2_{t,x}}\ ,$$
where we have used the assumption $\mu\le 1/3$ to ensure 
$-(1+\mu_1)/2-1/2-\mu\ge -(n-1)/2$.
The second term on the right could be 
controlled by using Proposition \ref{thm-wchain-high} and Lemma 
\ref{thm-w-Hardy}, as follows
$$
\|w^{-1}r^{-\frac 12+\mu_1}D^{\frac 12}\nabla g(u)\|_{L^2_{t,x}}
\les
\|w^{-1}r^{-\frac 12+\mu_1}D^{\frac 32} u\|_{L^2_{t,x}}
\les
\|w^{-1}D^{2-\mu_1} u\|_{L^2_{t,x}}\les \ep_1+\ep_s\ .
$$
Instead, concerning the first term on the right, we use Lemma 
\ref{thm-w-Hardy} to obtain
\begin{eqnarray*}
\|w^{-1}r^{-\frac 12-\mu}D^{\frac 12}\nabla g(u)\|_{L^2_{t,x}} & \les & \|w^{-1}r^{-\mu}\Delta g(u)\|_{L^2_{t,x}} \\
& \les & \|w^{-1}r^{-\mu}\Delta u\|_{L^2_{t,x}}
+\|w^{-1}r^{-1}\nabla u\|_{L^2_{t,x}}
\|r^{1-\mu}\nabla u\|_{L^\infty_{t,x}}\\
&\les&\ep_1+\ep_s\ .
\end{eqnarray*}

\section{Appendix: Proof of Morawetz type estimates, Lemma
\ref{th-LE}}\label{sec-LEprf}

In this section, we are interested in proving the fundamental Morawetz type estimates, Theorem
\ref{th-LE}.
Let $S_T=[0,T)\times \R^n$ with $n\ge 3$, we consider the  linear wave equations 
\eqref{eq-3.2}, that is
\beeq\label{eq-9.1}
\Box_h u:=
(-\pt^2+\Delta +\tilde h^{\al\be}(t,x)\pa_\al\pa_\be)u=F\ ,\eneq 
where we assume
$\tilde h^{\al\be}= h^{\al\be}-m^{\al\be}$,
 $h^{\al\be}=h^{\be\al}$, 
 $\tilde h^{00}=0$
 and 
$\Box_h$ is uniform hyperbolic, in the sense of
\eqref{eq-3.1}.

\begin{lem}[Morawetz type estimates]\label{thm-KSS-0}
Let $f=f(r)$ be any fixed differential function.
For any solutions $u\in C^\infty([0,T], C_0^\infty(\R^n))$ to the equation \eqref{eq-9.1}
 in $S_T$ with  \eqref{eq-3.1} and $n\ge 3$,
we have
\beeq\label{eq-9.2} (f X u)(
h^{\al\be}\pa_\al\pa_\be   u)
=\pa_\ga   P^\ga_h-
Q\ ,\eneq
where 
$Xu=\left(\pa_r+\frac{n-1}{2r}\right) u$,
$P^0_h=f   h^{0\beta}  \partial_\beta u  Xu
$, $P^j_h=\CO(( |f|+|rf'|+|f \tilde h|)   |\tilde\pa u|^2)$,
$$Q=Q_0+\CO(
( \frac{|f\tilde h|} r+|\pa (f\tilde h)|) |\pa u| |\tilde\pa u|
)\ ,$$
  \beeq\label{eq-5.5LE}Q_0=  \frac{ 2f-r f' }r \frac{|\ang u|^2}2+f' \frac{|\pa_r u|^2+
|\pt u|^2}{2}
-\frac{n-1}4\Delta \left(\frac{f }{r}\right) u^2\ ,\eneq
and $|\ang u|^2=|\nabla u|^2-|\pa_r u|^2$.
\end{lem}
This is essentially coming from multiplying $f(r) \left(\pa_r+\frac{n-1}{2r}\right) u$ to the wave equation and a tedious calculation of integration by parts. See e.g. \cite[P199-200]{MetSo06},
\cite[P273 (2.10)-(2.11)]{HWY1}.
Typically, $f$ is chosen to be differential functions
satisfying
\beeq\label{520-1}f\le 1, 2f \ge r f'(r)\ge 0,  -\Delta (f/r)\ge 0,\eneq
which ensure that $Q_0$ is positive semidefinite.
In literatures, some of the the typical choices are $f=1$ \cite{Mo68}, $1-(3+r)^{-\delta}$ ($\delta>0$) \cite{Strauss75}, $r/(R+r)$ 
\cite{Ster05, MetSo06}, $(r/(R+r))^{\mu}$ ($\mu\in (0,1)$, \cite{HWY1, HWY2}).

\subsection{Details: general case}
Let $\omega^j=\om_j=x^j/r$.
As $\pa_j=\omega_j\pa_r+\ang_j$,
$\pa_r=\omega^j \pa_j$, $2X=\pa_r-\pa_r^*=\omega^j \pa_j 
+\pa_j\omega^j=2\omega^j \pa_j +(n-1)/r$,
we have $[X,\pt]=0$, $$[X,\pa_k]=
[\omega^j,\pa_k] \pa_j+\frac{n-1}{2}[\frac 1 r, \pa_k]
=-\frac{\de^{j}_k-\omega^j\omega_k}{r}\pa_j+\frac{n-1}{2r^2}\omega_k
=\frac{1}r (-\pa_k+\omega_k X)\ .
$$

Notice that
\begin{eqnarray*}
\pa_\al\pa_\be u X u 
 & = & 
 \pa_\al (\pa_\be u X u)-\pa_\be u \pa_\al X u
   \\
 & = & 
 \pa_\al (\pa_\be u X u)-\pa_\be u [\pa_\al, X] u
 - \pa_\be u X \pa_\al u\ ,
\end{eqnarray*}
as $u_{\al\be}=u_{\be\al}$, we obtain
\begin{eqnarray*}
2\pa_\al\pa_\be u X u 
 & = & 
 \pa_\al (\pa_\be u X u)
 + \pa_\be (\pa_\al u X u)
 -\pa_\be u [\pa_\al, X] u
  -\pa_\al u [\pa_\be, X] u\\&&
 - \pa_\be u X \pa_\al u
  - \pa_\al u X \pa_\be u\ .\end{eqnarray*}
Notice also that  $\pa_j (\omega^j F G)=F XG+GXF
$,
we get
$$2\pa_\al\pa_\be u X u 
 = 
 \pa_\al (\pa_\be u X u)
 + \pa_\be (\pa_\al u X u)
 -\pa_\be u [\pa_\al, X] u
  -\pa_\al u [\pa_\be, X] u-\pa_j (\omega^j
 \pa_\be u  \pa_\al u)\ .$$

To be specific, we have
$$\pt^2 u X u 
= \pt (\pt u X u)
 -\pa_j(\omega^j \frac{u_t^2}2)\ ,
$$
$$2\pt \pa_j u X u 
= \pt (\pa_j u X u)+\pa_j  (\pt u X u)
-\pa_k (\omega^k \pt u  \pa_j u)
+ \frac{u_t
(-\pa_k+\omega_k X)u
}r \ ,
$$
\begin{align*}
2\pa_j\pa_k u X u =    \pa_j (\pa_k u X u)
& + \pa_k (\pa_j u X u)
-\pa_m (\omega^m
 \pa_j u  \pa_k u)
 \\
    &  +
\frac{ \pa_j u  (-\pa_k+\omega_k X)u
+\pa_k u(-\pa_j+\omega_j X)u
}r\ .
\end{align*}
 In summary, we have
\beeq\label{eq-multiplier-1} (\pa_\al\pa_\be u ) X u =\pa_\ga P^\ga_{\al\be}+Q_{\al\be}\ ,\eneq
with
$P^\ga_{\al\be}=\CO(|\pa u||\tilde\pa u|)$, $r Q_{\al\be}=\CO(|\pa u||\tilde\pa u|)$.

\subsection{Details: general multiplier} By \eqref{eq-multiplier-1}, when we multiply $f Xu$ against $\Box u$, we get
$$
\pa^\al\pa_\al u f X u =
\pa_\ga (f m^{\al\be}P^\ga_{\al\be})
-f'\omega_k m^{\al\be}P^k_{\al\be}
+fm^{\al\be}Q_{\al\be}\ ,
$$ where
$$m^{\al\be}P^0_{\al\be}=-P^0_{00}+\sum_j P^0_{jj}=-\pt u X u\ ,$$
$$m^{\al\be}P^k_{\al\be}
=-P^k_{00}+\sum_j P^k_{jj}=
\omega^k \frac{u_t^2}2
+\sum_j (\de^k_j u_j Xu-\frac 12\omega^k u_j^2)
=\omega^k \frac{u_t^2-|\nabla u|^2}2+u_k Xu\ ,
$$
$$ m^{\al\be}Q_{\al\be}= -  Q_{00}+  \sum_j  Q_{jj}
=\sum_j \pa_j u 
\frac 1r (-\pa_j+\omega_j X)u
 =\frac{1}r(-|\nabla u|^2+u_r^2)+\frac{n-1}{2r^2} u\pa_r u
 \ .
 $$
Notice that 
$\frac{u\pa_r u}{r^2} =
\frac{1}{r^2} \pa_r u^2=
-\frac{1}{2 }\nabla r^{-1}\cdot \nabla u^2
\ ,
$
$$2f\frac{
u\pa_r u
}{r^2} =
-\nabla \frac fr\cdot \nabla u^2+\frac{f'}r \pa_r u^2 
= 
-\nabla\cdot\left(u^2 \nabla \frac fr\right)+u^2 \Delta\left(\frac fr\right)
+\frac{ f'} r  \pa_r u^2
,
$$
we obtain
$$-f'\omega_k m^{\al\be}P^k_{\al\be}
=-f'
(
\frac{u_t^2-|\nabla u|^2}2+u_r Xu)\ ,
$$
and
 $$
f  m^{\al\be}Q_{\al\be}=
-\frac{f}r |\ang u|^2+
\frac{n-1}{4 }\left(\Delta\left(\frac fr\right)u^2
+\frac{ f'} r  \pa_r u^2
\right)
+\pa_j F^j\ ,
$$ with
$F^j=\CO(( |f|+|rf'|) r^{-2} u^2 )$.

In summary, we have
\begin{eqnarray*}
(\pa^\al\pa_\al u) f X u &=&
\pa_\ga (f m^{\al\be}P^\ga_{\al\be})
-f'\omega_k m^{\al\be}P^k_{\al\be}
+fm^{\al\be}Q_{\al\be}
\\
&=&\pa_\ga \tilde P^\ga-\frac{f}r |\ang u|^2+
\frac{n-1}{4 }\Delta(\frac fr)u^2-f'
(
\frac{u_t^2-|\nabla u|^2}2+u_r^2)\\
&=&
\pa_\ga \tilde P^\ga-f' \frac{u_t^2+u_r^2}{2}
+\frac{n-1}{4}\Delta(\frac{f}{r})  u^2
-\frac{(2f-rf')|\ang u|^2}{2r}
\\
&=&\pa_\ga \tilde P^\ga-Q_0\ ,
\end{eqnarray*}
where $\tilde P^j=\CO(( |f|+|rf'|)   |\tilde\pa u|^2)$,
$\tilde P^0=-f u_t Xu$.

For perburtation, we have
\beeq f \tilde h^{\al\be}(\pa_\al\pa_\be u) X u =\pa_\ga (f \tilde h^{\al\be} P^\ga_{\al\be})
- P^\ga_{\al\be} \pa_\ga (f\tilde h^{\al\be})
+f \tilde h^{\al\be} Q_{\al\be}\ .
\eneq
In summary, we obtained
\eqref{eq-9.2}.

\subsection{Choice of multiplier function $f$}
To prove the Morawetz type estimates, Lemma
\ref{th-LE}, we will choose two kinds of the multiplier functions $f$, with parameter $R>0$,
\beeq\label{eq-mul1}f=\frac{r}{R+r}\ ,\eneq
\beeq\label{eq-mul2}f=\(\frac{r}{R+r}\)^\mu=
\(1-\frac{R}{R+r}\)^\mu
,\ \mu\in (0, 1)\ .\eneq
Of course, \eqref{eq-mul1} could be viewed as the limit case of \eqref{eq-mul2} when $\mu=1$.

Now we do the calculation for $f$ given in \eqref{eq-mul2} with $\mu\in (0,1]$.
 We first notice that
\beeq\label{eq-5.9}
 f'(r) =\mu\(\frac{r}{R+r}\)^{\mu-1}
\frac{R}{(R + r)^2}=\mu \frac{R r^{\mu-1}}{(R+r)^{\mu+1}}\ge 0
\ ,\eneq
\beeq\label{eq-5.10}
 \frac{f(r)}{r} - f'(r) =
  \frac{r^{\mu- 1}}{(R + r)^{\mu}}\left(
1 - \frac{\mu R}{R + r}\right)\ge 0\ ,
\eneq
and
\beeq\label{eq-5.11}
 \frac{2f(r)-rf'(r)}{r} \ge  \frac{f(r)}{r} \ge f'(r)\ .
 \eneq
  In order to compute $-\Delta (f(r)/r)$, we recall that
$$- \Delta \left(\frac{f(r)}{r}\right)
 =- r^{1-n}\partial_r\left(r^{n-1}\partial_r \frac{f(r)}{r}\right) 
= r^{1-n}\partial_r\left(r^{n-2}\left(\frac{f(r)}{r}-f'(r)\right)\right).
$$
Using this identity and \eqref{eq-5.10}, we see that
$-\Delta (f(r)/r)$ equals to
\begin{eqnarray*}
 && r^{1-n}\partial_r\left(
  \frac{r^{n+\mu- 3}}{(R + r)^{\mu}}\left(
1 - \frac{\mu R}{R + r}\right)
\right) \\
 &=& \left(\frac{(n+\mu-3)r^{-3+\mu}}{(R+r)^{\mu}}
-\frac{\mu r^{-2+\mu}}{(R+r)^{\mu+1}}\right)
 \left(1-\frac{\mu R}{R+r}\right) +\frac{\mu R r^{-2+\mu}}{(R+r)^{\mu + 2}}\\
 &= &
\frac{r^{-3+\mu}}{(R+r)^{\mu}}  \left((n+\mu-3)
-\frac{\mu r}{R+r}\right)
 \left(1-\frac{\mu R}{R+r}\right)
+\frac{\mu R r^{-2+\mu}}{(R+r)^{\mu + 2}}\\
 &= &
\frac{r^{-3+\mu}}{(R+r)^{\mu}}  \left(n-3+
\frac{\mu R}{R+r}\right)
 \left(1-\frac{\mu R}{R+r}\right)
+\frac{\mu R r^{-2+\mu}}{(R+r)^{\mu + 2}} ,
\end{eqnarray*}
from which we see that, as $n\ge 3$,
\begin{equation}\label{eq-5.12}
- \Delta \left(\frac{f(r)}{r}\right) \ge
(1-\mu)
\frac{\mu R^2 r^{-3+\mu}}{(R+r)^{\mu+2}} 
+\frac{\mu R r^{-2+\mu}}{(R+r)^{\mu + 2}}\ge 0\ .
\end{equation}

In summary, we see that when $\mu\in (0,1)$, we get $Q_0$ from \eqref{eq-5.5LE} is non-negative and has the following lower bound for $r\le R$
\beeq\label{eq-5.13}Q_0\ge f' \frac{|\pa u|^2}{2}-\frac{n-1}4\Delta \left(\frac{f }{r}\right) u^2
\ges_\mu  \frac{|\tilde\pa u|^2}{R^{\mu}r^{1-\mu}}\ ,
\eneq
 where the implicit constant depends only on $n$ and $\mu\in (0,1)$, and in particular, independent of $R>0$.
On the other hand, when $\mu=1$, 
 $Q_0$ from \eqref{eq-5.5LE} is still non-negative and has the following lower bound for $R/2\le r\le R$
\beeq\label{eq-5.14}Q_0\ge f' \frac{|\pa u|^2}{2}-\frac{n-1}4\Delta \left(\frac{f }{r}\right) u^2
\ge \frac{1}{8R}  |\pa u|^2+
\frac{ n-1 }{32}\frac{1}{R^{2}r} u^2
\ges \frac{|\tilde\pa u|^2}{r}\ .
\eneq

\subsection{Proof of Morawetz type estimates,  Theorem
\ref{th-LE}.}
Equipped with Lemma \ref{thm-ener} and  Lemma \ref{thm-KSS-0}, together with the observations \eqref{eq-5.13}-\eqref{eq-5.14}, we could give the proof of Morawetz type estimates,  Lemma
\ref{th-LE}.

Let us begin with the proof of \eqref{eb4}. At first, applying
\eqref{eq-5.13} with $R=1$, and \eqref{eq-5.14} with $R\ge 1$, that is we use $f=\(\frac{r}{1+r}\)^\mu$ and $\frac{r}{R+r}$ with $R\ge 1$, we get
\begin{eqnarray*}
&  &\int_{r\le 1} \frac{|\tilde\pa u|^2}{r^{1-\mu}} dx dt+\sup_{R\ge 1}
\int_{R/2\le r\le R} \frac{|\tilde\pa u|^2}{r} dx dt
 \\ & \les &
 \sup_{f=\(\frac{r}{1+r}\)^\mu, \frac{r}{R+r}, R\ge 1}\int_{S_T} Q_0 dxdt
 \\
 & \les & \sup_{f=\(\frac{r}{1+r}\)^\mu, \frac{r}{R+r}, R\ge 1}\(
- \int_{S_T} f F \left(\pa_r+\frac{n-1}{2r}\right) u  dx dt+\left. \int_{\R^n} P_h^0(t, \cdot)dx\right|_{t=0}^T\right.
 \\
 &  &+\left. \int_{S_T} Q_0-Q dxdt\)
 \\ & \les &
  \int_{S_T} |F\tilde\pa u| dxdt
 +\int_{\R^n} |\pa u(T)||\tilde \pa u(T)|dx
  +\int_{\R^n} |\pa u(0)||\tilde \pa u(0)|dx\\
  &&+
  \sup_{f=\(\frac{r}{1+r}\)^\mu, \frac{r}{R+r}, R\ge 1}
\int_{S_T} \(|f\pa h|+\frac{|f \tilde h|}r\)  |\pa u||\tilde \pa u| dxdt
 \\ & \les &
  \int_{S_T} |F\tilde\pa u| dxdt
 +\|\tilde \pa u(t)\|_{L^\infty L^2(S_T)}^2+
\int_{S_T} \(|\pa h|+\frac{| \tilde  h|}{r^{1-\mu}\<r\>^\mu}\)  |\pa u||\tilde \pa u| dxdt
\end{eqnarray*}
where we have used 
\eqref{eq-9.2} in the second inequality,
the facts
$|f|\le 1$, $0\le f'\le f/r$,
$$|Q-Q_0|\les
\(|f\pa h|+|f' \tilde  h|+\frac{|f  \tilde h|}r\)  |\pa u||\tilde \pa u|
\les
\(|f\pa h|+\frac{|f  \tilde h|}r\)  |\pa u||\tilde \pa u|
\ ,
$$
and
$|P^0|\les |\pa u||\tilde \pa u|$ in the third inequality.
By Lemma \ref{thm-ener} and Hardy's inequality, we see that
\begin{eqnarray*}
\|u\|^2_{X_1}& :=&\int_{r\le 1} \frac{|\tilde\pa u|^2}{r^{1-\mu}} dx dt+\sup_{R\ge 1}
\int_{R/2\le r\le R} \frac{|\tilde\pa u|^2}{r} dx dt
+\|\pa u(t)\|_{L^\infty L^2(S_T)}^2
 \\ & \les &
  \int_{S_T} |F\tilde\pa u| dxdt
 +\|\pa u(0)\|_{ L^2}^2+
\int_{S_T} \(|\pa h|+\frac{| \tilde  h|}{r^{1-\mu}\<r\>^\mu}\)  |\pa u||\tilde \pa u| dxdt\ .
\end{eqnarray*}
Thus to give \eqref{eb4}, we need only to show that
\beeq\label{eq-KSS}
\|u\|_{LE_{T}}\les \|u\|_{X_1}\ ,
\eneq
which essentially follows from a standard argument of
Keel-Smith-Sogge \cite{KSS}.
Here for completeness, we write down the proof. The first and second terms are trivial to control.
For the remaining two terms, 
with $\al\in [0, \mu]$, we have
\begin{eqnarray*}
&&
\|r^{-\frac{1-\mu}2}
\<r\>^{-\frac{\al}2}
 \tilde\partial u\|_{L^2_{t,x}(S_T)}^2
 \\&
\les&
\|r^{\frac{\mu-1}2}\tilde\partial u\|_{L^2 (r\le 1)}^2+\sum_{0\le j\le \ln \<T\>}\|r^{\frac{\mu-1-\al}2}\tilde\partial u\|_{L^2 (r\simeq 2^j)}^2
+\|r^{\frac{\mu-1-\al}2}\tilde\partial u\|_{L^2 (r\ge \<T\>)}^2
\\
&
\les&\|u\|_{X_1}^2+\sum_{0\le j\le \ln \<T\>}2^{j(\mu-\al)} \|r^{-\frac{1}2}\tilde\partial u\|_{L^2 (r\simeq 2^j)}^2
+\<T\>^{\mu-1-\al}\|\tilde\partial u\|_{L^2_{t,x} (r\ge \<T\>)}^2
\\
&
\les&\sum_{0\le j\le \ln \<T\>} 2^{j (\mu-\al)}\|u\|_{X_1}+\<T\>^{ \mu-\al }
 \|\tilde\partial u\|_{L^\infty_t L^2_x } \les   C_\al (T)  \|u\|_{X_1}^2
 \ ,
\end{eqnarray*}
where
$$C_\al(T)=\left\{
\begin{array}{ll }
\ln \<T\>      &   \al=\mu\ , \\
  \<T\>^{ {\mu-\al}}    &   \al\in [0,\mu)\ .
\end{array} \right.$$
This completes the proof of \eqref{eq-KSS}, and so is  \eqref{eb4}.

Turning to the proof of \eqref{eb4-vari}, we will use
$f=\(\frac{r}{T+r}\)^\mu\le 1$.
Applying
\eqref{eq-5.13}, we get as before,
\begin{eqnarray*}
&  &\int_{r\le T} \frac{|\tilde\pa u|^2}{T^\mu r^{1-\mu}} dx dt
 \\ & \les &
  \int_{S_T} (|F|+|f\tilde\pa \tilde h|  |\pa u|)|\tilde\pa u| dxdt
 +\|\tilde\pa u(t)\|_{L^\infty L^2(S_T)}^2
 \\ & \les &
\int_{S_T} \(|F|+\(|\pa h|+\frac{|  \tilde h|}{r^{1-\mu}(T+r)^\mu}\)  |\pa u|\) |\tilde \pa u| dxdt+\|\tilde\pa u(t)\|_{L^\infty L^2(S_T)}^2
\ .
\end{eqnarray*}
Together with Lemma \ref{thm-ener}, we see that
\begin{eqnarray*}
\|u\|^2_{X_2}& :=&\int_{r\le T} \frac{|\tilde\pa u|^2}{T^\mu r^{1-\mu}} dx dt
+\|\pa u(t)\|_{L^\infty L^2(S_T)}^2
 \\ & \les &
\int_{S_T} \(|F|+\(|\pa h|+\frac{|  \tilde h|}{r^{1-\mu}(T+r)^\mu}\)  |\pa u|\) |\tilde \pa u| dxdt +\|\pa u(0)\|_{ L^2}^2\ .
\end{eqnarray*}
Moreover,
we have
\begin{eqnarray*}
\|r^{-\frac{1-\mu}2}\tilde\partial u\|_{L^2_{t,x}}&
\les&
\|r^{-\frac{1-\mu}2}\tilde\partial u\|_{L^2 (|x|\le T)}+\|r^{-\frac{1-\mu}2}\tilde\partial u\|_{L^2 (|x|\ge T)}
\\
&
\les&
\|r^{-\frac{1-\mu}2}\tilde\partial u\|_{L^2 (|x|\le T)}
+T^{-\frac{1-\mu}2}\|\tilde\partial u\|_{L^2_{t,x} (|x|\ge T)}
\\
&
\les&\|r^{-\frac{1-\mu}2}\tilde\partial u\|_{L^2 (|x|\le T)}+
T^{\frac{\mu}2}\|\tilde\partial u\|_{L^\infty_t L^2_x }\les T^{\frac{\mu}2} \|u\|_{X_2}
\end{eqnarray*}
which gives us \eqref{eb4-vari}.

\subsection*{Acknowledgment}
The author would like to thank Sergiu Klainerman for stimulating
discussions during the author's visits in the Spring of 2017 and Summer of 2018. He would also like to thank the annonymous referee for the careful reading and helpful comments.
The author was supported in part by NSFC 11971428 and  National Support Program for Young Top-Notch Talents.


\begin{thebibliography}{10}

\bibitem{MR2003417}
Serge Alinhac.
\newblock An example of blowup at infinity for a quasilinear wave equation.
\newblock {\em Ast\'{e}risque}
No. 284, pages 1--91. 2003.

\bibitem{MR1728676}
Hajer Bahouri and Jean-Yves Chemin.
\newblock \'{E}quations d'ondes quasilin\'{e}aires et effet dispersif.
\newblock {\em Internat. Math. Res. Notices}, (21):1141--1178, 1999.

\bibitem{MR1719798}
Hajer Bahouri and Jean-Yves Chemin.
\newblock \'{E}quations d'ondes quasilin\'{e}aires et estimations de
  {S}trichartz.
\newblock {\em Amer. J. Math.}, 121(6):1337--1377, 1999.

\bibitem{BL}
J{\"o}ran Bergh and J{\"o}rgen L{\"o}fstr{\"o}m.
\newblock {\em Interpolation spaces. {A}n introduction}.
\newblock Springer-Verlag, Berlin-New York, 1976.
\newblock Grundlehren der Mathematischen Wissenschaften, No. 223.

\bibitem{CN16}
David Cruz-Uribe and Virginia Naibo.
\newblock Kato-{P}once inequalities on weighted and variable {L}ebesgue spaces.
\newblock {\em Differential Integral Equations}, 29(9-10):801--836, 2016.

\bibitem{Da19}
Piero D'Ancona.
\newblock A short proof of commutator estimates.
\newblock {\em J. Fourier Anal. Appl.}, 25(3):1134--1146, 2019.

\bibitem{MR3583356}
Boris Ettinger and Hans Lindblad.
\newblock A sharp counterexample to local existence of low regularity solutions
  to {E}instein equations in wave coordinates.
\newblock {\em Ann. of Math. (2)}, 185(1):311--330, 2017.

\bibitem{FW1}
Daoyuan Fang and Chengbo Wang.
\newblock Local well-posedness and ill-posedness on the equation of type
  {$\square u=u^k(\partial u)^{\alpha}$}.
\newblock {\em Chinese Ann. Math. Ser. B}, 26(3):361--378, 2005.

\bibitem{FW2}
Daoyuan Fang and Chengbo Wang.
\newblock Some remarks on {S}trichartz estimates for homogeneous wave equation.
\newblock {\em Nonlinear Anal.}, 65(3):697--706, 2006.

\bibitem{FWQLW}
Daoyuan Fang and Chengbo Wang.
\newblock Local existence for nonlinear wave equation with radial data in 2+1
  dimensions.
\newblock Preprint. ArXiv:0705.2849 [math.AP], 2007.

\bibitem{FW11}
Daoyuan Fang and Chengbo Wang.
\newblock Weighted {S}trichartz estimates with angular regularity and their
  applications.
\newblock {\em Forum Math.}, 23(1):181--205, 2011.

\bibitem{HJLW20p}
Kunio Hidano, Jin-Cheng Jiang, Sanghyuk Lee, and Chengbo Wang.
\newblock Weighted fractional chain rule and nonlinear wave equations with
  minimal regularity.
\newblock {\em Rev. Mat. Iberoam.}, 36(2):341--356, 2020.

\bibitem{HWY2}
Kunio Hidano, Chengbo Wang, and Kazuyoshi Yokoyama.
\newblock The {G}lassey conjecture with radially symmetric data.
\newblock {\em J. Math. Pures Appl. (9)}, 98(5):518--541, 2012.

\bibitem{HWY1}
Kunio Hidano, Chengbo Wang, and Kazuyoshi Yokoyama.
\newblock On almost global existence and local well posedness for some 3-{D}
  quasi-linear wave equations.
\newblock {\em Adv. Differential Equations}, 17(3-4):267--306, 2012.

\bibitem{HWY3}
Kunio Hidano, Chengbo Wang, and Kazuyoshi Yokoyama.
\newblock Combined effects of two nonlinearities in lifespan of small solutions
  to semi-linear wave equations.
\newblock {\em Math. Ann.}, 366(1-2):667--694, 2016.

\bibitem{HiYo06}
Kunio Hidano and Kazuyoshi Yokoyama.
\newblock Space-time {$L^2$}-estimates and life span of the
  {K}lainerman-{M}achedon radial solutions to some semi-linear wave equations.
\newblock {\em Differential Integral Equations}, 19(9):961--980, 2006.

\bibitem{MR0420024}
Thomas J.~R. Hughes, Tosio Kato, and Jerrold~E. Marsden.
\newblock Well-posed quasi-linear second-order hyperbolic systems with
  applications to nonlinear elastodynamics and general relativity.
\newblock {\em Arch. Rational Mech. Anal.}, 63(3):273--294 (1977), 1976.

\bibitem{Jo81}
Fritz John.
\newblock Blow-up for quasilinear wave equations in three space dimensions.
\newblock {\em Comm. Pure Appl. Math.}, 34(1):29--51, 1981.

\bibitem{MR1031987}
L.~V. Kapitanski\u{\i}.
\newblock Estimates for norms in {B}esov and {L}izorkin-{T}riebel spaces for
  solutions of second-order linear hyperbolic equations.
\newblock {\em Zap. Nauchn. Sem. Leningrad. Otdel. Mat. Inst. Steklov. (LOMI)},
  171(Kraev. Zadachi Mat. Fiz. i Smezh. Voprosy Teor. Funktsi\u{\i}.
  20):106--162, 185--186, 1989.

\bibitem{KSS}
Markus Keel, Hart~F. Smith, and Christopher~D. Sogge.
\newblock Almost global existence for some semilinear wave equations.
\newblock {\em J. Anal. Math.}, 87:265--279, 2002.
\newblock Dedicated to the memory of Thomas H. Wolff.

\bibitem{KM95}
S.~Klainerman and M.~Machedon.
\newblock Smoothing estimates for null forms and applications.
\newblock {\em Duke Math. J.}, 81(1):99--133, 1995.
\newblock A celebration of John F. Nash, Jr.

\bibitem{MR1962783}
S.~Klainerman and I.~Rodnianski.
\newblock Improved local well-posedness for quasilinear wave equations in
  dimension three.
\newblock {\em Duke Math. J.}, 117(1):1--124, 2003.

\bibitem{KlMa93}
Sergiu Klainerman and Matei Machedon.
\newblock Space-time estimates for null forms and the local existence theorem.
\newblock {\em Comm. Pure Appl. Math.}, 46(9):1221--1268, 1993.

\bibitem{KlRo05}
Sergiu Klainerman and Igor Rodnianski.
\newblock Rough solutions of the {E}instein-vacuum equations.
\newblock {\em Ann. of Math. (2)}, 161(3):1143--1193, 2005.

\bibitem{MR3402797}
Sergiu Klainerman, Igor Rodnianski, and Jeremie Szeftel.
\newblock The bounded {$L^2$} curvature conjecture.
\newblock {\em Invent. Math.}, 202(1):91--216, 2015.

\bibitem{KS97}
Sergiu Klainerman and Sigmund Selberg.
\newblock Remark on the optimal regularity for equations of wave maps type.
\newblock {\em Comm. Partial Differential Equations}, 22(5-6):901--918, 1997.

\bibitem{LiZh95}
Ta-Tsien Li and Yi~Zhou.
\newblock A note on the life-span of classical solutions to nonlinear wave
  equations in four space dimensions.
\newblock {\em Indiana Univ. Math. J.}, 44(4):1207--1248, 1995.

\bibitem{Ld93}
Hans Lindblad.
\newblock A sharp counterexample to the local existence of low-regularity
  solutions to nonlinear wave equations.
\newblock {\em Duke Math. J.}, 72(2):503--539, 1993.

\bibitem{Ld96}
Hans Lindblad.
\newblock Counterexamples to local existence for semi-linear wave equations.
\newblock {\em Amer. J. Math.}, 118(1):1--16, 1996.

\bibitem{MR1666844}
Hans Lindblad.
\newblock Counterexamples to local existence for quasilinear wave equations.
\newblock {\em Math. Res. Lett.}, 5(5):605--622, 1998.

\bibitem{Ld08}
Hans Lindblad.
\newblock Global solutions of quasilinear wave equations.
\newblock {\em Amer. J. Math.}, 130(1):115--157, 2008.

\bibitem{LiuW21}
Mengyun Liu and Chengbo Wang.
\newblock Concerning ill-posedness for semilinear wave equations.
\newblock {\em Calc. Var. Partial Differential Equations}, 60(1):19, 2021.

\bibitem{MaNaNaOz05}
Shuji Machihara, Makoto Nakamura, Kenji Nakanishi, and Tohru Ozawa.
\newblock Endpoint {S}trichartz estimates and global solutions for the
  nonlinear {D}irac equation.
\newblock {\em J. Funct. Anal.}, 219(1):1--20, 2005.

\bibitem{MetSo06}
Jason Metcalfe and Christopher~D. Sogge.
\newblock Long-time existence of quasilinear wave equations exterior to
  star-shaped obstacles via energy methods.
\newblock {\em SIAM J. Math. Anal.}, 38(1):188--209, 2006.

\bibitem{MeTa12MA}
Jason Metcalfe and Daniel Tataru.
\newblock Global parametrices and dispersive estimates for variable coefficient
  wave equations.
\newblock {\em Math. Ann.}, 353(4):1183--1237, 2012.

\bibitem{MW17}
Jason Metcalfe and Chengbo Wang.
\newblock The {S}trauss conjecture on asymptotically flat space-times.
\newblock {\em SIAM J. Math. Anal.}, 49(6):4579--4594, 2017.

\bibitem{MR1168960}
Gerd Mockenhaupt, Andreas Seeger, and Christopher~D. Sogge.
\newblock Local smoothing of {F}ourier integral operators and
  {C}arleson-{S}j\"{o}lin estimates.
\newblock {\em J. Amer. Math. Soc.}, 6(1):65--130, 1993.

\bibitem{Mo68}
Cathleen~S. Morawetz.
\newblock Time decay for the nonlinear {K}lein-{G}ordon equations.
\newblock {\em Proc. Roy. Soc. Ser. A}, 306:291--296, 1968.

\bibitem{MuSch13}
Camil Muscalu and Wilhelm Schlag.
\newblock {\em Classical and multilinear harmonic analysis. {V}ol. {I}}, volume
  137 of {\em Cambridge Studies in Advanced Mathematics}.
\newblock Cambridge University Press, Cambridge, 2013.

\bibitem{PS}
Gustavo Ponce and Thomas~C. Sideris.
\newblock Local regularity of nonlinear wave equations in three space
  dimensions.
\newblock {\em Comm. Partial Differential Equations}, 18(1-2):169--177, 1993.

\bibitem{MR1644105}
Hart~F. Smith.
\newblock A parametrix construction for wave equations with {$C^{1,1}$}
  coefficients.
\newblock {\em Ann. Inst. Fourier (Grenoble)}, 48(3):797--835, 1998.

\bibitem{SmTa05}
Hart~F. Smith and Daniel Tataru.
\newblock Sharp local well-posedness results for the nonlinear wave equation.
\newblock {\em Ann. of Math. (2)}, 162(1):291--366, 2005.

\bibitem{Ster05}
Jacob Sterbenz.
\newblock Angular regularity and {S}trichartz estimates for the wave equation.
\newblock {\em Int. Math. Res. Not.}, (4):187--231, 2005.
\newblock With an appendix by Igor Rodnianski.

\bibitem{Ster07}
Jacob Sterbenz.
\newblock Global regularity and scattering for general non-linear wave
  equations. {II}. {$(4+1)$} dimensional {Y}ang-{M}ills equations in the
  {L}orentz gauge.
\newblock {\em Amer. J. Math.}, 129(3):611--664, 2007.

\bibitem{Strauss75}
Walter~A. Strauss.
\newblock Dispersal of waves vanishing on the boundary of an exterior domain.
\newblock {\em Comm. Pure Appl. Math.}, 28:265--278, 1975.

\bibitem{Tataru99}
Daniel Tataru.
\newblock On the equation {$\square u=|\nabla u|^2$} in {$5+1$} dimensions.
\newblock {\em Math. Res. Lett.}, 6(5-6):469--485, 1999.

\bibitem{MR1749052}
Daniel Tataru.
\newblock Strichartz estimates for operators with nonsmooth coefficients and
  the nonlinear wave equation.
\newblock {\em Amer. J. Math.}, 122(2):349--376, 2000.

\bibitem{Ta02}
Daniel Tataru.
\newblock Strichartz estimates for second order hyperbolic operators with
  nonsmooth coefficients. {III}.
\newblock {\em J. Amer. Math. Soc.}, 15(2):419--442, 2002.

\bibitem{Tay}
Michael~E. Taylor.
\newblock {\em Tools for {PDE}}, volume~81 of {\em Mathematical Surveys and
  Monographs}.
\newblock American Mathematical Society, Providence, RI, 2000.
\newblock Pseudodifferential operators, paradifferential operators, and layer
  potentials.

\bibitem{W17}
Chengbo Wang.
\newblock Long-time existence for semilinear wave equations on asymptotically
  flat space-times.
\newblock {\em Comm. Partial Differential Equations}, 42(7):1150--1174, 2017.

\bibitem{WangQ17}
Qian Wang.
\newblock A geometric approach for sharp local well-posedness of quasilinear
  wave equations.
\newblock {\em Ann. PDE}, 3(1):Art. 12, 108, 2017.

\bibitem{YinZhou16}
Silu Yin and Yi~Zhou.
\newblock Global existence of radial solutions for general semilinear
  hyperbolic systems in 3{D}.
\newblock {\em J. Math. Anal. Appl.}, 433(1):49--73, 2016.

\bibitem{Zh03}
Yi~Zhou.
\newblock On the equation {$\square\phi=|\nabla\phi|^2$} in four space
  dimensions.
\newblock {\em Chinese Ann. Math. Ser. B}, 24(3):293--302, 2003.

\bibitem{ZhouLei08}
Yi~Zhou and Zhen Lei.
\newblock Global low regularity solutions of quasi-linear wave equations.
\newblock {\em Adv. Differential Equations}, 13(1-2):55--104, 2008.

\end{thebibliography}

\end{document}